\documentclass[11pt]{article}
\usepackage{graphicx}
\usepackage{verbatim}
\newcommand{\RR}{\mathrm{I\!R\!}}
\newcommand{\ZZ}{{\Bbb Z}}
\newcommand{\QQ}{{\Bbb Q}}
\newcommand{\KK}{{\Bbb K}}
\newcommand{\M}[1]{\mathcal{#1}}
\newcommand{\MU}{\ensuremath{\mathcal{U} }}
\newcommand{\MA}{\ensuremath{\mathcal{A} }}
\newcommand{\MMM}{\ensuremath{\mathcal{MM} }}
\newtheorem{Def}{Definition}
\newtheorem{Prop}[Def]{Proposition}
\newtheorem{Prop+Def}[Def]{Proposition and Definition}
\newtheorem{Th}[Def]{Theorem}

\input amssym.def
\input amssym
\begin{document}
\title{ Some Combinatorial Aspects of Movies and Movie-moves
 in the Theory of Smoothly Knotted Surfaces in $\RR ^4$}
\author{Glenn Lancaster\\DePaul University\\Chicago, Illinois
\and Richard Larson\\University of Illinois at
Chicago\\Chicago,Illinois \and
Jacob Towber\\DePaul
University\\Chicago,Illinois  }
\date{\today}
\maketitle
\section{Introduction}\label{S:intro}
Throughout this paper, by a \emph{knotted surface} (or,
synonymously, \emph{2-knot}) will always be meant a \emph{smooth}
embedding in $\RR^{\;4}$ of a compact 2-manifold $M$ (where $M$ is
smooth, not necessarily connected, and is not assigned an
orientation.) \footnote{While there are also highly non-trivial
questions concerning smooth embeddings of such $M$ in $\RR^{\;5}$,
the present paper will study \emph{only} smooth embeddings in
$\RR^{\;4}$. For $n\geq 6$, any two smoothly knotted surfaces in
$\RR^{\;n}$ are smoothly ambient isotopic, so this third case is
without further interest. The analogous theories of 2-knots in the
PL- and continuous categories, are not studied in the present
paper. } By \emph{isotopy} of two such smoothly knotted surfaces
in $\RR^{\;4}$, will be meant smooth ambient isotopy.

There is a beautiful combinatorial description of the
isotopy-classes of such smooth knotted surfaces in $\RR^{\;4}$,
via the "doubly time-oriented movies" of Carter, Rieger and Saito
(\cite{CRS}), explained below in \S~\ref{S:movies}. Putting all
this in a more general context, the collection of smooth isotopy
classes of knotted surfaces is a portion of the collection of
``\emph{two-tangles}'', whose intricate structure (that of a "free
braided monoidal 2-category with duals")has been elucidated by the
work of Fischer (\cite{F}), Kharlamov and Turaev(\cite{KT}) Baez
and Neuchl (\cite{BN}) and especially in the recent work of Baez
and Langford (\cite {BL}, and \cite{LL}). Indeed, the method of
CRS movies comes into its own, only with the extension to
2-tangles.
\medskip

This theory of movies involves a remarkable interplay between
low-dimensional differential topology, combinatorics, quantum
groups and 2-category theory. In the present paper (except for a
few elementary topological remarks , in \S \ref{SS:Stills},
 reviewing the classical theory of knotted curves in $\RR^3$), we
confine our attention strictly to the {\bf combinatorial} portion
of this theory---i.e., the portion related to such concepts as
elementary transitions, flickers, and Carter-Rieger-Saito movies,
explained in \S \ref{SS:ET} and \S \ref{S:movies}---and also the
combinatorial concept of movie-moves, explained in \S
\ref{SS:movieIsotopy} and \S \ref{SS:moviemoves}.
\smallskip

The material which presents itself for study when we thus focus
our attention, is divided in the present paper into the three
following portions (in Sections 2, 3, 4 respectively).

In Section 2 of the present paper, we review the combinatorial
portion of the material on movies in [CS2], [CRS], [LL] and [BL].
While this represents only a small portion of the full theory of
these movies, these combinatorial movies, and the operations on
them which occur naturally in this paper, form an intricate
structure which it may be worthwhile disentangling from the
existing literature.

In Section 2, this material is analyzed in an extremely detailed
way. This level of detail is necessary, both for the purely
mathematical purposes of the present paper, and also to make
possible the program(s) described below.

More speculatively, the extremely detailed combinatorial
information presented here for knotted surfaces, might possibly be
useful data in the future extension of these results to yet
higher-dimensional generalizations. It seems fascinating to
speculate, about what might be the third term in the sequence of
combinatorial structures which begins:\\
$\mathcal{E}=\mbox{set of elementary events }=\{Cup,Cap,NE,NW\}$,
5 Reidemeister and Yetter moves  \\
--- followed by:\\
$\mathcal{ET}=\mbox{set of elementary transitions, 31 types of
movie-moves.}$

Another advantage to restricting to the combinatorial portion of
the theory of CRS movies, is that it seems much easier to write a
program for directly computing with this portion of the theory.
Actually, we have written two such programs, based on somewhat
different principles---with the idea that, when they both give the
same result, this increases the reliability of the answer. These
two (equivalent) programs--- or rather, class libraries, one in
C++ and one in Java--- will be referred to in this paper, as
C++2KnotsLib and Java2KnotsLib.
 \footnote{These class libraries, together with some driver programs
 relevant to the present paper,
are available (including source code) on
the three following URL's:\\
  http://www.uic.edu/\textasciitilde rgl/javaknots.html ---this contains
   Java2KnotsLib\\
  while C++2KnotsLib may be downloaded from either
  of:\\
  http://condor.depaul.edu/\textasciitilde jtowber/2knots $\;\;\;$
  http://condor.depaul.edu/\textasciitilde glancast/2knots}

 It should be emphasized that Section 2 is entirely an exposition
of a subset of the results in [CRS] and [BL]. In Section 3 we
present constructions, next to be sketched, which we believe have
some novelty. Their explicit computation involves use of the
programs and class-libraries just mentioned.

 Let $\mathcal{MM}$ denote the set of all movie-moves. This will
be explained in more detail in \S \ref{SS:moviemoves} below. For
now, it suffices to know that each movie-move is associated with a
positive integer, which we shall call its \underline{{\bf type}},
which ranges from 1 to 31 inclusive. One application of the
constructions in (III) is to study some aspects of this set
$\mathcal{MM}$, in relation to the theory of 2-knots.
\smallskip

 It is natural to ask, which of these 31 types of movie-move is
really needed, in the sense of the following definition. For
instance, it will be proved in this paper, that (in the sense of
the following definition) movie-move 31 is not a consequence of
the remaining 30 types.
\begin{Def}\label{def:essentialMM}
Let $1\leq i\leq 31$. We shall say that movie-moves of type $i$
are \underline{{\bf essential}} if there exist two CRS movies M
and M' such that (working in the smooth category):
\begin{enumerate}
\item M and M' have 0 strings coming in and 0 going out---ie the
2-tangles they represent are knotted surfaces.
\item The knotted
surfaces represented by M and M' are ambient isotopic.
\item It is impossible to get from M to M', by a sequence of movie-moves,
all of type different from $i$.
\end{enumerate}
When such movies $M,M'$ do not exist, we shall say that
movie-moves of type i are \underline{{\bf non-essential}.}
\end{Def}
Using the constructions in \S 3, we shall prove that movie-moves
of type 31 are essential. This result is proved in \S 3.6 below,
utilizing the machinery next to be explained.
\medskip

Define
$$\underline{31}:=\{i\in \ZZ\; |\; 1\leq i\leq 31\}\;,$$
and let $\mathcal{U}$ denote any proper subset of
$\underline{31}$. Let
$$\mathcal{U}^c:=\{i\in \underline{31}\; |\; i\notin \mathcal{U}\}$$
denote the complement of $\mathcal{U}$ in $\underline{31}$.
\medskip

 We now need the following
higher-dimensional analogue of Kauffman's concept of ``regular
isotopy'' of link-diagrams:
\begin{Def}\label{def:Uregular}
Let $\M{U}$ and $\M{U}^c$ be as above. Let $M$ and $M'$ be two
movies.\\ Then we say $M$ and $M'$ are
 \underline{{\bf $\M{U}$-regularly isotopic}} if it is possible to
go from $M$ to $M'$ by a sequence of movie-moves NOT in $\M{U}$\\
---in more detail, if there is an non-empty ordered sequence
$$\M{M}_1,\cdots,\M{M}_N \;\;\;(N>0)$$
of movie-moves, with these three properties:
\begin{enumerate}
\item Each $\M{M}_i$ (for $1\leq i\leq N)$ is a movie-move of type
in $\M{U}^c$ (i.e., NOT in $\M{U}$.)
\item If $\M{U}_i$ is the movie-move $(M_i==M_i')$, then
$$M_i'=M_{i+1}\mbox{ for }1\leq i\leq N-1$$
\item $M_1=M$ and $M_N'=M'$.
\end{enumerate}
\end{Def}
Another key ingredient in our proof that movie-moves of type 31
are essential, is the construction of certain $\M{U}$-regular
isotopy invariants, which are here called \underline{{\bf
$\M{U}$-amplitude-invariants}}. The construction of these
amplitude-invariants, which is next to be sketched, is based on
$sl_q(2)$ and the Kauffman amplitude $<>$.
\smallskip

Let
$$A=\ZZ[q,q^{-1}]$$
and let $V$ be free over $A$ of rank 2. As explained in \S
\ref{SS:Stills}, to every elementary event or still
$$E:m\rightarrow n\;,$$
the Kauffman amplitude assigns an $A$-linear map
$$<E>:V^{\otimes m} \rightarrow V^{\otimes n}$$
For technical reasons, our constructions become much simpler, if
we replace $A$ by its quotient-field
$$\KK:=\QQ (q)$$
and $V$ by the 2-dimensional vectorspace
$$V_1:=V\otimes _A \KK$$
over $\KK$, at the same time (by a slight abuse of notation)
identifying the Kauffman bracket $<E>$ with
$$<E>\otimes _E \KK:V_1^{\otimes m} \rightarrow V_1^{\otimes n}$$

 We introduce in \S \ref{S:amp} a higher-dimensional analog of the
 Kauffman bracket:
namely, the concept of an \underline{{\bf amplitude-assignment}
$\mathcal{A}$}, which assigns to every elementary transition
$\mathcal{E}:S=>T:m\rightarrow n$, an $\KK$-linear map
$$<\mathcal{E}>_{\mathcal{A}}:V_1^{\otimes m} \rightarrow
V_1^{\otimes n}$$ Such an $\mathcal{A}$ assigns in a natural way
(explained in Def.\ref{Def:ampFlicker}) to every flicker
$$F:S=>T:m\rightarrow n$$ an $\mathcal{A}$-amplitude
$<F>_{\mathcal{A}}$. Finally, the $\mathcal{A}$-amplitude of a
movie $M$ is defined to be given by the following construction:
\smallskip

Let the movie $M$ consist of the stills $S_1,\cdots,S_N$, with
consecutive stills $S_i$, $S_{i+1}$ joined by the flicker $F_i$
for $1\leq i\leq N-1$. It is an important feature of the
Carter-Rieger-Saito concept of movie, that two movies which have
the same stills but different flickers between them, count in
general as non-isotopic movies, and can represent non-isotopic
2-knots. This feature suggests that the $\mathcal{A}$-amplitude of
$M$ should be defined in such a way as to involve the $F$'s as
well as the $S$'s; we here set
$$<M>_{\mathcal{A}}:=<S_1>-<F_1>_{\mathcal{A}}+ <S_2>-<F_2>_{\mathcal{A}}+
<S_3>-\cdots  +<S_N>\eqno(*)$$
\medskip

 We define an amplitude-assignment $\mathcal{A}$ to be
\underline{{\bf $\MU$-balanced}} if the expression (*) just given,
provides an $\M{U}$-regular isotopy-invariant of movies. The main
purpose of the present paper, is to construct (by using computer
programs to be discussed below) $\MU$-balanced
amplitude-assignments $\mathcal{A}$, which then yield
$\M{U}$-regular isotopy invariants $<M>_\mathcal{A}$ of knotted
surfaces. $\M{U}$-regular isotopy invariants obtained in this way,
are here referred to as $sl_q(2)$-based \footnote{This represents
a pious hope that these constructions may be generalized to apply
to other quantum groups. The work presented here is organized
strongly around the Kauffman bracket, which has $sl_q(2)$ as its
quantum symmetry group, acting on its two-dimensional fundamental
representation.} amplitude-invariants.

It is easy to see that an amplitude-assignment $\mathcal{A}$ is
$\M{U}$-balanced, if and only if it `respects' (in a suitable
sense given by Def.\ref{def:respect}) all movie-moves not in
$\M{U}$. \emph{A priori}, there are infinitely many movie-moves to
be considered in this connection. At this point it is helpful to
add a somewhat technical requirement, that $\mathcal{A}$ be
\underline{{\bf semi-normal}}, explained in
Def.\ref{def:seminormalAA} of \S \ref{SS:balanced}. With this
\emph{Ansatz}, we may reduce to the consideration of only a finite
collection of movie-moves, as follows:

Using the numbering of movie-moves furnished both in [CRS] and in
[BL], it is proved that if $\mathcal{A}$ is semi-normal, it
automatically respects all movie-moves of types 15, 16, 22, (Prop.
\ref{def:seminormalAA}) and also all movie-moves of types
17,18,19,20 (Prop. \ref{prop:semiRespects }) Moreover, the
movie-moves of the 23 remaining types can be reduced to a finite
subset, such that an amplitude-assignment $\mathcal{A}$ is
$\MU$-balanced, if and only if all movie-moves in this finite
subset, are respected by $\mathcal{A}$. (Prop.\ref{finiteTests}).
\smallskip

In order to complete our construction of $\MU$-balanced
amplitude-assignments,  we must augment our requirement of
semi-normality, obtaining a stronger condition of \underline{{\bf
strong normality }}, as given by Def.\ref{def:normal}. Now the
search for $\MU$-balanced strongly normal amplitude-assignments,
reduces to a finite number of linear equations over $A$, via an
explicit algorithm furnished by Prop.\ref{prop:Assoc} in \S
\ref{SS:Assoc}.
\smallskip

We apply the machinery thus developed, to prove, in \S 3.6, the
result (asserted above) that movie-move 31 is \emph{essential} in
the sense of Def. \ref{def:essentialMM}, i.e. is not a consequence
of the remaining movie-moves. The proof to be presented here,
involves explicitly solving the system Assoc($\mathcal{U}$) of
equations associated to $\mathcal{U}$, in the two special cases
$\mathcal{U}=\{31\}$ and $\mathcal{U}=\emptyset$.
\smallskip

Finally, \S 4 is devoted to showing how our Program can be used to
perform some miscellaneous computations on 7 sample movies
(representing, respectively: an unknotted sphere, an unknotted
Klein bottle, an unknotted torus, two knotted spheres, and the
`1-twist and 2-twist spun trefoil' described on p.36 of [CS2].)
\medskip

 We should like to express our gratitude to Carter
and Saito for a large number of conversations in which they
patiently explained some of their ideas about knotted surfaces.
Without their help, this paper could not have been written. Also,
the third author would like to express gratitude to the
mathematics departments of Harvard, University of Chicago, and
University of Illinois at Chicago, for extending hospitality
during time in which work was done on the present paper.

\section{Carter-Rieger-Saito Doubly Time-oriented Movies}\label{S:movies}
\subsection{Time-oriented Link-diagrams and their Kauffman Amplitudes}\label{SS:Stills}
By a ``link-diagram in $\RR^{\;2}$ " will be meant an ordered pair
(D,C), where D is a smooth immersion in $\RR^{\;2}$ of a finite
disjoint union of circles, which is ``generic", that is, such that
the only singularities of D are simple transverse intersections;
and where C is the additional structure furnished by the
assignment for each singularity S of D, with branches b1,b2 over
S, of a label "over-crossing" for one of these two branches, and
"under-crossing" for the other one.

 If L is a smooth link in
$\RR^{\;3}$, then for almost every unit vector $v$ in $\RR^{\;3}$, the projection
 $$\RR^{\;3} \rightarrow  \RR^{\;2}$$ along $v$ maps L
to the portion D of such a link-diagram, with the additional
structure C then determined so the direction from each
over-crossing to the corresponding under-crossing is parallel (in
the positive sense) to $v$.

\begin{figure}
\centering
\includegraphics[scale=0.5]{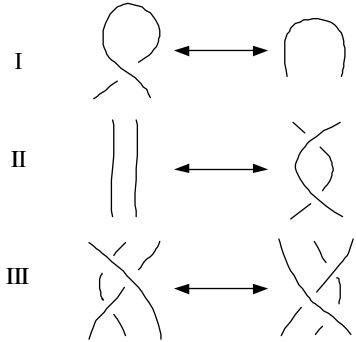}
\caption{\label{fig:reid1_3} The Three Reidemeister Moves }
\end{figure}

The reader is assumed familiar with the fact that two such link
diagrams are projections of isotopic links, if and only if it is
possible to pass from one to the other by planar isotopies,
together with a finite sequence of the three "Reidemeister moves"
(as pictured in Fig. \ref{fig:reid1_3}.) Moreover, certain
mathematical constructions, reminiscent of constructions in
quantum mechanics, have (since the revolutionary advent in 1984 of
the Jones polynomial(\cite{J})) been used to import "amplitudes"
into this situation, and so to obtain isotopy-invariants of
links.The reader will not be assumed familiar with the subset of
these constructions utilized in this paper (Yetter's two additions
to the Reidemeister moves, and the relation of these additions to
the concept of braided monoidal categories and to the Kauffman
amplitudes under ``regular isotopy".) This subset of these
familiar constructions will next be discussed in the remainder of
the present section, from a perspective which (it is hoped) will
make more natural the higher-dimensional analogues of these
constructions in the later portions of this paper.

\medskip

Suppose we introduce an oriented ``time" dimension, into the real
Euclidean plane P containing the link diagrams under study---i.e.,
suppose we have chosen a specific unit vector $\tau$ in P
(pointing from `past'to  `future'). (In all figures drawn in the
present paper, $\tau$ will point from bottom to top. Many papers
on the subject use the opposite convention.) This extra bit of
structure makes possible the following further well-known
modifications in the preceding set-up---and so makes possible the
construction of the
 monoidal category $\Sigma$  (of \underline{'stills'}) next to be described
(a quotient of which is a braided monoidal category with duals,
the category of \underline{'tangles'}, which seems to underly the
theory of the Jones , HOMFLY and Kauffman polynomials.)

\begin{Def}\label{tau.regular}
Let $\tau$ be a unit vector in a Euclidean plane $P$, and let
$$(D:n.S^1\rightarrow P,C)$$
 be a link-diagram in $P$ (where we denote by
$n.S^1$ the disjoint union of n circles.)

Then we shall say that D is \underline{$\tau$-regular} if, for every singular point $\pi$
of D (which, by hypothesis, is a simple transverse self-intersection) neither of the two
tangents at $\pi$ is parallel to $\tau$.

Two such $\tau$-regular link-diagrams $(D,C),(D',C')$ with
$$D,D':n.S^1\rightarrow P$$ will be called
\underline{$\tau$-equivalent} if there is a smooth map $\phi:(n.S^1\times [0,1])\rightarrow
P,$ with these two properties:

i) For each $t\in [0,1]$,
$$\phi_t:n.S^1\rightarrow P, x\mapsto \phi (x,t)$$
is a $\tau$-regular link-diagram.

ii) $\phi_0=D$ and $\phi_1=D'.$

(iii) In the obvious sense, $\phi$ takes $C$ into $C'$.
\smallskip

When such a $\phi$ exists, we shall write $$(D,C)\sim_1 (D',C')\;.$$
\end{Def}
\bigskip

\begin{figure}
\mbox{
\begin{minipage}[b]{5.0cm}
\centering
\includegraphics[scale=0.5]{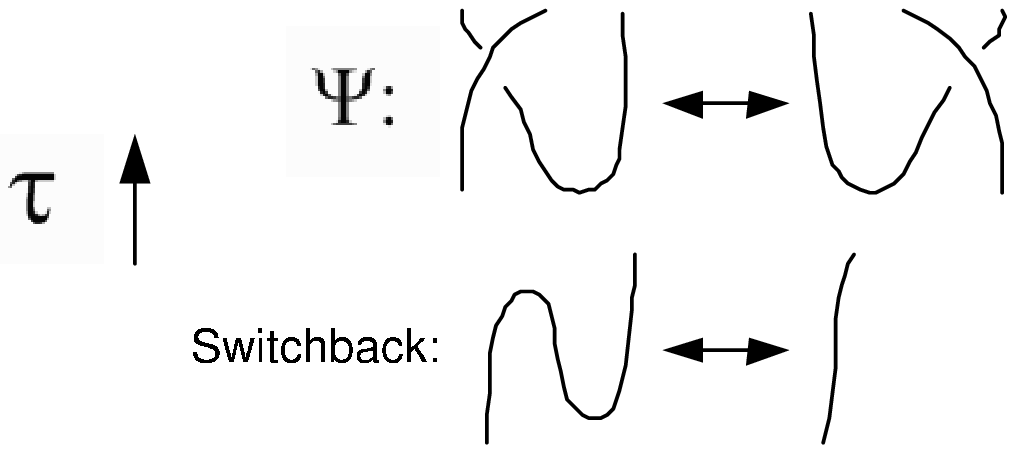}
\caption{\label{fig:yetter} The Yetter Moves  }
\end{minipage}
}
\hfill
\mbox{
\begin{minipage}[b]{6.0cm}
\centering
\mbox{
\includegraphics[scale=0.5]{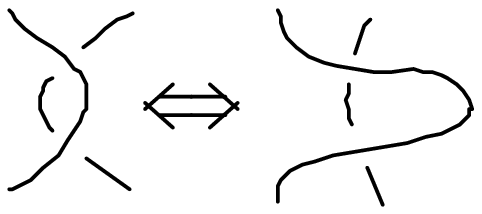} }
\begin{minipage}[b][0.5\height][b]{4.0cm}
\mbox{}\\ 
\end{minipage}

\caption{\label{fig:isotopy_move0}Move 0 = Planar~Isotopy}
\end{minipage}
}
\end{figure}

The point of all these definitions, is the following well-known refinement of Reidemeister's
results (which is precisely made-to-order for the extension below to knotted surfaces):
\begin{Th}\label{prop:Yetter} Let $(D,C),(D',C')$ be
 $\tau$-regular projections of smooth links
$L$ resp. $L'$ in $\RR^{\;3}$; then, $L$ and $L'$ are ambient isotopic, if
and only if it is possible
to go from $(D,C)$ to $(D',C')$ by a finite sequence
$$(D_0,C_0)=(D,C),(D_1,C_1),\cdots,(D_s,C_s)=(D',C')$$
such that, for each $i$ with $1\leq i<s$, {\bf either}
$$(D_i,C_i)\sim_1 (D_{i+1},C_{i+1})\eqno(*)$$
("move 0"), {\bf or} $(D_i,C_i)$ and $(D_{i+1},C_{i+1})$ are
related by one of the (various flavors of) the five moves in
Figures  \ref{fig:reid1_3} and \ref{fig:yetter}) (i.e. either by
one of the three Reidemeister moves in Fig. \ref{fig:reid1_3}, or
one of the two Yetter moves in Fig. \ref{fig:yetter}). When this
is the case, we shall write
$$(D,C)\sim_2 (D',C')\;.$$
\end{Th}
\smallskip

(Let us be a bit pedantic, and say that when (*) holds,
$(D_i,C_i)$ is obtained from $(D_{i+1},C_{i+1})$ via Move 0: this
will be a useful point of view later on---cf. Fig.
\ref{fig:isotopy_move0})
\smallskip

There is yet a third equivalence relation, \underline{regular isotopy}
which plays a role, via the Kauffman amplitudes shortly to be defined, and which we
shall denote by $\sim_3$. It is important to note that only the equivalence-relation
$\sim_1$ is used for the 1-morphisms (stills) in the 2-category of movies.

Let us denote by $\mathcal{D}^{\tau}_0$ the set of all
$\tau$-regular link-diagrams in $\RR^{\;2}$. This is contained in
a larger set $\mathcal{D}^{\tau}$ which constitutes a braided
monoidal category with duals, the category of `tangles', as
explained in [FY]. For the purposes of the present paper, we do
not need the details of this construction; we only need the
quotient $\Sigma^{\tau}$  of $\mathcal{D}^{\tau}$ by the
equivalence relation extending $\tau$-equivalence, or rather, we
only need the purely combinatorial construction of this category,
which we next present in some detail.

\begin{Def}\label{FEvent}By an \underline{elementary event} will be meant an element of the
set of strings $$\mathcal{E}=\{Cup,Cap,NE,NW\}\;.$$ By a
\underline{framed event} will be meant one of two kinds of object:
{\bf EITHER} an ordered triple (always to be enclosed in {\bf
square} brackets) of the form $$[m,E,n]\eqno (1a)$$ where $m$ and
$n$ are natural numbers, and  $$E\in \mathcal{E} \;,$$ {\bf OR} a
symbol of the form $$1_n \eqno(1b)$$ where $n$ is a natural
number..
\end{Def}

\begin{figure}[h]
\centering
\includegraphics[scale=0.6]{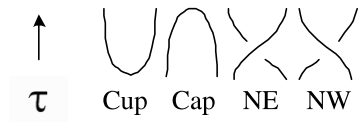}
\caption{\label{fig:events}The 4 Elementary Events}
\end{figure}

\begin{figure}[h]
\begin{minipage}[b]{6.0cm}
\centering
\includegraphics[scale=0.5]{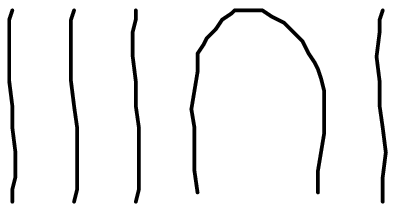}
\caption{\label{fig:3cap1} Framed Event [3,Cap,1]}
\end{minipage}
\hfill
\begin{minipage}[b]{5.5cm}
\centering
\includegraphics[scale=0.5]{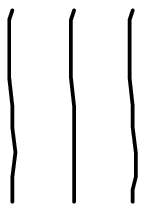}
\caption{\label{fig:1_3} Framed Event $1_3$  }
\end{minipage}
\end{figure}

\noindent {\bf Note:} The 4 elementary events are to be thought of
as associated with the 4 $\tau$-regular link diagrams given in
Figure \ref{fig:events}. Also (with $\tau$ always drawn
vertically, from the bottom of the page to its top) $[m,E,n]$ is
to be thought of, as associated with the $\tau $-regular
link-diagram, obtained by adjoining to the diagram for $E$, $m$
vertical strings to the left, and $n$ vertical strings to the
right---e.g., Figure \ref{fig:3cap1} represents the framed event
$[3,Cap,1]$. Finally,$1_n$ is associated with the regular
link-diagram, consisting of $n$ vertical lines---e.g. Figure
\ref{fig:1_3} represents the framed event $1_3$. In particular,
$1_0$ is associated with the 'empty' link-diagram, containing no
strings at all.

\medskip

\begin{Def}\label{source}By the \underline{source} and \underline{target} of an elementary
event $E$ will be meant the integers specified by the following
table. The \underline{source} and \underline{target} of a framed event
are then defined by the formulas
$$source([m,E,n]):=m+n+source(E)\;;target([m,E,n]):=m+n+target(E)$$
and $$source(1_n):=n,\;\;target(1_n):=n$$
\end{Def}
\begin{tabular}{ccccc}
Event & Source & Target\\
\hline
NW & 2 & 2\\
NE & 2 & 2\\
Cup & 0 & 2\\
Cap & 2 & 0
\end{tabular}
\medskip

\noindent {\bf Note:} Thus, with the conventions in use in the
present paper, the source of an event, or framed event, is the
number of lines leading into its base, and the target is the
number of lines going out its top---e.g. the framed event
$[3,Cap,1]$ in Fig. \ref{fig:3cap1} has source 6 and target 4.

\begin{Def}\label{Def:still} By a \underline{still} $S$ will be meant
a
finite non-empty sequence of framed  events,
$$(F_1,\cdots,F_s)$$
such that
for each of these framed events $F_i$ (except the last one) the target of $F_i$ is
equal to the source of the succeeding framed event $F_{i+1}$.We shall usually
use a multiplicative notation to denote such a concatenation of framed events,
writing $$S=F_1F_2\cdots F_s$$
 We denote by $\Sigma$
the set of such stills. By the \underline{source} of a non-empty still $S$ will
be meant the source of its first framed event, while its \underline{target} is defined
to be the target of its last framed event. By the \underline{{\bf empty still}}
will be meant the still $1_0$; note that, by the preceding
definitions,the associated link-diagram to the empty still is empty, and
the source and target of the empty still are both 0.
\end{Def}

\noindent{\bf NOTE: } Such a still is called a \underline{`word' } by Carter,Rieger
and Saito---cf.{\cite{CS}}, p.35.

\begin{figure}[h]
\centering
\includegraphics[scale=0.9]{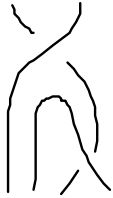}
\caption{\label{fig:still_example}Still [2,NW,0][1,Cap,1][0,NE,0]}
\end{figure}

\noindent {\bf EXAMPLE}:

 Consider the still $S\in\Sigma$ given by
$$S=[2,NW,0][1,Cap,1][0,NE,0]$$
which is illustrated by Figure \ref{fig:still_example} (which is
to be read from bottom to top, while $S$ is read left-to-right.
(\underline{Caution:} some other authors use other conventions.)

Note this sequence of 3 framed events
satisfies the compatibility condition of Def.\ref{Def:still}, since
$$target([2,NW,0])=4=source[1,Cap,1]$$
 and $$target([1,Cap,1])=2=source([0,NE,0]\;.$$
Note also that $S$ has source 4 and target 2.

\bigskip

The set $\Sigma$, gives rise in the following way to a category
 (which, by a slight abuse of notation, will also be denoted by $\Sigma$):

The \underline{\bf {objects}} of the category $\Sigma$ are the natural
numbers; the \underline{\bf {morphisms}} are the stills. If $m$ and $n$ are natural
numbers,$Hom_{\Sigma}(m,n)$ is the set of all stills S such that
$$source(S)=m \mbox{ and }target(S)=n$$
in which case we shall also write
$$S:m\rightarrow n\;.$$
\begin{Def}\label{def:stillComp}
If $S,S'$ are the stills defined, respectively, by the
sequences\footnote{By the conventions used in this paper, these
compositions are read from {\bf bottom} to {\bf top}, e.g. $F_1$
occurs at the bottom of the still $S$, and $F_s$ at the top.}
$$S=F_1\cdots F_s,\;\;S'=F'_1\cdots F'_t$$
of framed events, then the composite still $S\circ S'$ is defined in the
category $\Sigma$, if and only if
$$source(S')=target(S)\;,$$
in which case we set
$$S\circ S'=F_1\cdots F_sF'_1\cdots F'_t$$
\end{Def}
\medskip

 The following operation on the
category $\Sigma$ will prove useful below:
\begin{Def}\label{def:circ}
Let $m$ and $n$ be natural numbers.\\
a)If $F$ is a framed event, we define $m\bullet F \bullet n$ to be the framed event
$$
\left\{
\begin{array}{ll}
[m+m',E,n+n'] & \mbox{if }F=[m',E,n']\\
1_{m+n+p} & \mbox{if }F=1_p
\end{array}
\right.
$$
b)Let
$$S=F_1\cdots F_p$$
be a still (where $1\leq p$ and $F_1,\cdots,F_p$ are framed events.)\\
Then we define $m\bullet S \bullet n$ to be the still
$$m\bullet S \bullet n:=(m\bullet F_1 \bullet n)\cdots(m\bullet F_p \bullet n)$$
\end{Def}
\medskip

In other words (by a slight generalization of the construction in Def.\ref{FEvent}), the
link-diagram for the
still $m\bullet S \bullet n$ is obtained from that for the still $S$, by
appending $m$ vertical strings to the left, and n to the right.
\smallskip

Let us next note the following two (commuting) symmetries $t$ and $f$ on the
collection $\Sigma$ of stills:
First, $f$ is the symmetry which "flips" over-crossings into under-crossings; in more detail:
\begin{Def}\label{still:f}\underline{\bf THE SYMMETRY $f$:}\\
(i)We define the action of $f$ on the elementary events in $\mathcal{E}$ by:
$$f\cdot Cap= Cap, f\cdot Cup= Cup, f\cdot NE= NW,f\cdot NW= NE$$
(ii)We define the action of $f$ on framed events by:
$$f\cdot[m,E,n]=[m,f\cdot E,n]\mbox{ for $E\in \mathcal{E}$, and }f\cdot 1_n=1_n$$
(iii)Finally, if $S$ is the still given by
$$S=F_1\cdots F_s\;,$$
where the $F_i$ are framed events, then we define
$$f\cdot S:=(f\cdot F_1)\cdots(f\cdot F_s)$$
\end{Def}
Next, $t$ is the symmetry which reverses the time-orientation $\tau$; in more detail:
\begin{Def}\label{still:t}\underline{\bf THE SYMMETRY $t$:}\\
(i)We define the action of $t$ on the elementary events in $\mathcal{E}$ by:
$$t\cdot Cap= Cup, t\cdot Cup= Cap, t\cdot NE= NW, t\cdot NW= NE$$
(ii)We define the action of $t$ on framed events by:
$$t\cdot[m,E,n]=[m,t\cdot E,n]\mbox{ for $E\in \mathcal{E}$, and }t\cdot 1_n=1_n$$
(iii)Finally, if $S$ is the still given by
$$S=F_1\cdots F_s\;,$$
where the $F_i$ are framed events, then we define
$$t\cdot S:=(t\cdot F_s)\cdots(t\cdot F_1)$$
\end{Def}
\smallskip

We next review one well-known further aspect of $\Sigma$ which
will be used in our later constructions---namely, the construction
(due to L.Kauffman) which assigns to every still $S\in \Sigma$ a
`Kauffman amplitude' $<S>$:\footnote{It is more usual to call this
the `Kauffman bracket'. Our present terminology is intended to
distinguish this from the other {\bf Kauffman bracket}, which has
a similar definition and notation, but whose domain of definition
is the set of non-time-oriented link-diagrams (of non-oriented
links or tangles). }

\bigskip
\noindent Let V be free on $\{e_0,e_1\}$ over the ring
$\ZZ[q,q^{-1}]$ of Laurent polynomials in one indeterminate $q$
over the ring $\ZZ$ of integers. Thus $V^{\otimes n}$ is free of
rank $2^n$ over $\ZZ [q,q^{-1}]$, with free basis
$$\{e(i_1,\cdots,i_n):=e_{i_1}\otimes e_{i_2}\otimes\cdots e_{i_n} |\mbox{
 all i's=0 or 1}\}$$
Let $S:m\rightarrow n$ be a still with source $m$ and target $n$.
(With the  conventions described above, this means $m$ strings
lead into the bottom of the time-oriented diagram describing $S$,
and $n$ lead up out the  top.) Then, we construct as follows a
linear transformation
$$<S>:V^{\otimes m}\rightarrow V^{\otimes n}\;.$$

In the first place, to the four elementary events in $\mathcal{E}$,
 we assign the following four $\ZZ [q,q^{-1}]$-linear maps as amplitudes:

$<Cup>: V^{\otimes 0}=\ZZ[q,q^{-1}]\rightarrow V^{\otimes 2}$
maps 1 to $-qe(0,1)+q^{-1}e(1,0)$ ;

\smallskip
$<Cap>: V^{\otimes 2}\rightarrow \ZZ[q,q^{-1}]$ maps both $e(0,0)$ and $e(1,1)$ to 0,
maps $e(0,1)$ to $q$ and maps $e(1,0)$ to $-q^{-1}$;

\smallskip
$<NE>:V^{\otimes 2}\rightarrow V^{\otimes 2}$  maps $e(i,i)$ to $qe(i,i)$ (for i=0,1) and maps
$e(0,1)$ to $q^{-1}e(1,0)$ and $e(1,0)$ to $q^{-1}e(0,1)+(q-q^{-3})e(1,0)$;

\smallskip
$<NW>:V^{\otimes 2}\rightarrow V^{\otimes 2}$  maps $e(i,i)$ to $q^{-1}e(i,i)$ (for i=0,1)
and maps $e(1,0)$ to $qe(0,1)$ and $e(0,1)$ to $(q^{-1}-q^3)e(0,1)+qe(1,0)$.

\bigskip
This construction extends in the obvious way from elementary
events to framed events; namely, given an elementary event $E$
with source $s$ and target $t$, and given the related framed event
$F=[m,E,n]$ (with source $m+s+n$ and target $m+t+n$), the Kauffman
amplitude $$<F>:V^{\otimes m+s+n}\rightarrow V^{\otimes m+t+n}$$
is defined to be the tensor product $Id_{v^m}\otimes <E>\otimes
Id_{V^n}$. Also, $<1_n>$ is defined to be the identity map
$Id_{V^n}$.
\smallskip

Finally, consider the general still $S$. We may write $$S=F_1F_2\cdots F_s$$
where each $F_i$ is a framed event, and where, for
$1\leq i <s$, $$source(F_{i+1})=target(F_i)\;.$$
Then we define the Kauffman amplitude of $S$ to be
$$<S>=<F_s>\circ <F_{s-1}>\circ\cdots\circ <F_1>:V^{\otimes source(S)}\rightarrow
V^{\otimes target(S)}\;.$$

\medskip
\noindent {\bf REMARK:} For this construction, cf. for instance, Kauffman
 (\cite{K}, Part I, $\S 9^o$.) We have
here made one trivial modification in this construction---namely, Kauffman's original amplitudes\\
$<Cup>,<Cap>$ are respectively -i,i times those we shall use in the algorithm
presented in this paper. With this modification
 we can work, not in $\ZZ[i,q,q^{-1}]$, but rather in the smaller ground-ring $\ZZ [q,q^{-1}]$,
over which we shall later have a rather large number (over two
thousand)
 of simultaneous linear equations to solve.
This modification is `trivial' in one important sense: if $S$ is
a still which arises from a link-diagram, so that
$$source(S)=0=target(S)$$
then if S has $m$ minima (cups), it also has $m$
maxima (caps), so Kauffman's original amplitude,
and the modification to be used here, differ by a factor of
$$(-i)^m.i^m=1\;,$$
i.e. agree at least for such stills. Also, the present
modification spoils, neither the $sl_q(2)$-equivariance of the
amplitude, nor the fact that the Kauffman amplitude respects the
relation $\sim_3$ of regular isotopy.
\bigskip

\subsection{Elementary Transitions}\label{SS:ET}
The set
$$\mathcal{E}=\{Cup,Cap,NE,NW\}\;.$$
of `elementary events', discussed above, has in the
higher-dimensional theory of knotted surfaces, an analog
$\mathcal{ET}$ (made up of two parts, one of cardinality 66 and
the other consisting of 32 infinite sequences), which is now to be
discussed. The elements of $\mathcal{ET}$ will here be called {\bf
elementary transitions}; they are called {\bf fundamental
elementary string interactions} in [CRS]. Each elementary
transition $\mathcal{E}$ is an ordered pair of
stills\footnote{Note that, for the purposes of the present paper,
an elementary transition is a purely {\bf combinatorial} object---
it is no more and no less than an ordered pair of stills which
occurs somewhere in our list. For the {\bf topological}
significance of these, as isotopy classes of specific 2-tangles in
$\RR ^3$, cf. [CRS], [CS2] and [LL].}, of which the first is
called the {\bf source} of $\mathcal{E}$, and the second the {\bf
target} of $\mathcal{E}$. An elementary transition $\mathcal{E}$
is uniquely determined by its source $S$ and target $T$, in which
case we shall write
\begin{equation}\label{ET:source}
 \mathcal{E}=[S\Longrightarrow T]\;.
 \end{equation}
\medskip

The reader is asked to bear with the (perhaps somewhat pedantic)
detailed enumeration and labelling in this section of these
elementary transitions. This is necessary, because we shall need
carefully to assign individual `amplitudes' to each of these, as
explained below. More conjecturally, the detailed combinatorial
data thus obtained may be also helpful in future investigations,
in suggesting the still higher-dimensional analogues.
\medskip

The tangle category $\Sigma/(\sim_2$) sketched in \S\ref{SS:Stills}, may be thought of as
obtained via the following generators and relations construction:\\
i) The generators are the framed events, obtained by applying to the 4 elements of
$\mathcal{E}$, the framing construction of Def. \ref{def:circ} (together with
the identity morphisms $1_n$), and then taking the $\sim_2$-classes
of the results.\\
ii)The defining relations are furnished by the three Reidemeister moves,
the two Yetter moves, and the `zero' move. (It
is now necessary, carefully to draw these 6 moves in all possible
"flavors", and this will be done in the present section.)
\medskip

The analogous generators-and-relations presentation, for the
theory of knotted surfaces in $\RR ^4$, contains:\\ for
generators, objects("flickers") obtained by suitably framing the
elementary transitions; and for relations, the collection of
"movie-moves" explained below in \S\ref{SS:moviemoves}. It is an
observation, not original with the authors of the present paper
(rather, it is an idea which seems well-known among experts in
knotted surfaces) that there are  \underline{{\bf
generators}}---i.e., elementary transitions---for the
higher-dimensional theory, which correspond precisely to the
\underline{{\bf relations}} for the lower-dimensional theory (as
well as some `new' elementary transitions apparently not of this
sort, namely those of type ET6 and ET7.)---and that this is a
possible guide for what to expect, as the existing theory gets
extended in the future into yet higher dimensions.\\
Enough speculation---let us next carefully examine the intricate collection of
``elementary transitions" in some detail:\\
\smallskip

\begin{figure}
  \centering
  \includegraphics[scale=0.55]{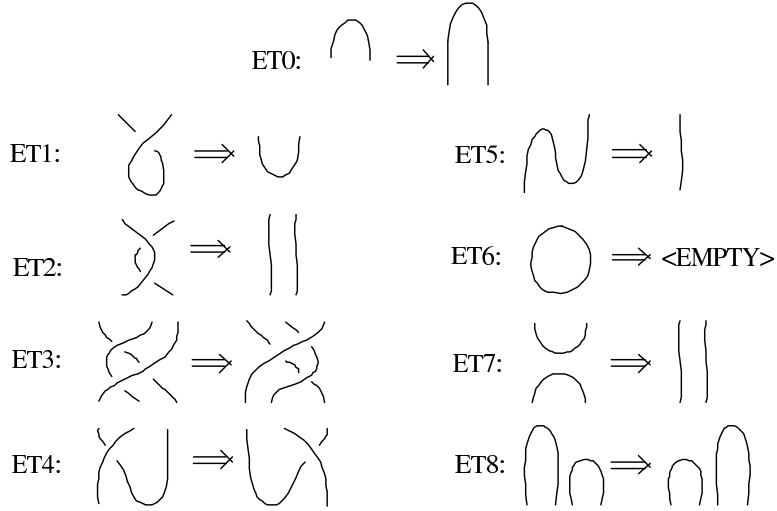}
  \caption{Elementary Transitions}\label{fig:alltypes}
\end{figure}

The  elementary transitions fall into 9 types, samples of which
are pictured in Fig. \ref{fig:alltypes}.

\smallskip

The second, third and fourth
 of these types, correspond (in accordance with
the principle stated above) to the three Reidemeister moves. Let
us next examine more closely the type ET1: This comes in eight
"flavors", i.e. there are 8 elementary transitions corresponding
to Reid I, the first Reidemeister move (cf. Fig.
\ref{fig:reid1_3}); these 8 are labelled as follows, and
illustrated in Figure \ref{fig:ET1}. (In \S\ref{AMPS} below, when
we come to assign "amplitudes" to elementary transitions, these 8
are assigned distinct amplitudes--- and indeed, this needs to be
done in many distinct ways. Thus it is essential to be able to
deal with these 8 different flavors on an individual basis---and
similarly for the other types treated below.)
\begin{eqnarray*}
ET1I&=[[0,Cup,0][0,NE,0]\Longrightarrow [0,Cup,0]],\\
ET1R&= [[0,Cup,0]\Longrightarrow[0,Cup,0][0,NE,0]],\\
ET1t&=[[0,NW,0][0,Cap,0]\Longrightarrow [0,Cap,0]],\\
ET1tR&=[[0,Cap,0]\Longrightarrow [0,NW,0][0,Cap,0]],\\
ET1f&=[[0,Cup,0][0,NW,0]\Longrightarrow [0,Cup,0]],\\
ET1fR&= [[0,Cup,0]\Longrightarrow[0,Cup,0][0,NW,0]],\\
ET1ft&=[[0,NE,0][0,Cap,0]\Longrightarrow [0,Cap,0]],\\
ET1ftR&=[[0,Cap,0]\Longrightarrow [0,NE,0][0,Cap,0]].
\end{eqnarray*}

\begin{figure}
 \centering
\includegraphics[scale=0.6]{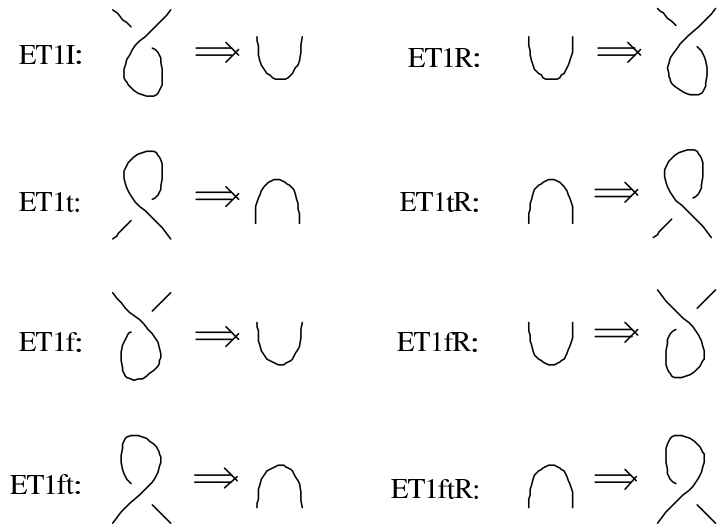}
\caption{\label{fig:ET1} Elementary Transitions of Type ET1}
\end{figure}
\begin{figure}
  \centering
  \includegraphics[scale=0.6]{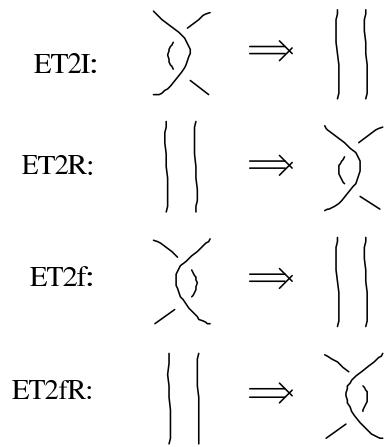}
  \caption{Elementary Transitions of Type ET2}\label{fig:ET2}
\end{figure}

(Note this ordering derives from the ordering $R<t<f$.)\\
Next, there are 4 elementary transitions corresponding to the second Reidemeister move, Reid II:
\begin{eqnarray*}
ET2I&=[[0,NE,0][0,NW,0]\Longrightarrow 1_2],\\
ET2R&=[1_2\Longrightarrow [0,NE,0][0,NW,0]],\\
ET2f&=[[0,NW,0][0,NE,0]\Longrightarrow 1_2],\\
ET2fR&=[1_2\Longrightarrow [0,NW,0][0,NE,0]].
\end{eqnarray*}
These are illustrated in Figure \ref{fig:ET2}.

\begin{figure}
  \centering
  \includegraphics[scale=0.5]{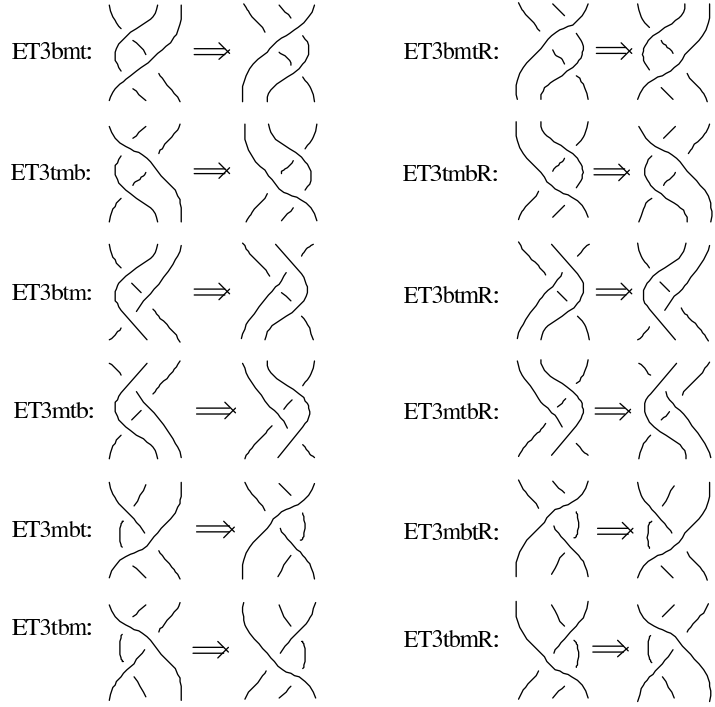}
  \caption{Elementary Transitions of Type ET3}\label{fig:ET3}
\end{figure}

For the third Reidemeister move, Reid III, the number of associated elementary
transitions increases to 12 (This is, in some ways, one of the most difficult of the
9 types. It is the type connected with the Yang-Baxter equation.)
These are listed as follows (and pictured in Figure \ref{fig:ET3}):
\begin{eqnarray*}
ET3bmt&=[[0,NE,1][1,NE,0][0,NE,1]=>[1,NE,0][0,NE,1][1,NE,0]],\\
ET3bmtR&=[[1,NE,0][0,NE,1][1,NE,0]=>[0,NE,1][1,NE,0][0,NE,1]],\\
ET3btm&=[[0,NW,1][1,NE,0][0,NE,1]=>[1,NE,0][0,NE,1][1,NW,0]],\\
ET3btmR&=[[1,NE,0][0,NE,1][1,NW,0]=>[0,NW,1][1,NE,0][0,NE,1]],\\
ET3mbt&=[[0,NE,1][1,NE,0][0,NW,1]=>[1,NW,0][0,NE,1][1,NE,0]],\\
ET3mbtR&=[[1,NW,0][0,NE,1][1,NE,0]=>[0,NE,1][1,NE,0][0,NW,1]],\\
ET3mtb&=[[0,NW,1][1,NW,0][0,NE,1]=>[1,NE,0][0,NW,1][1,NW,0]],\\
ET3mtbR&=[[1,NE,0][0,NW,1][1,NW,0]=>[0,NW,1][1,NW,0][0,NE,1]],\\
ET3tbm&=[[0,NE,1][1,NW,0][0,NW,1]=>[1,NW,0][0,NW,1][1,NE,0]],\\
ET3tbmR&=[[1,NW,0][0,NW,1][1,NE,0]=>[0,NE,1][1,NW,0][0,NW,1]],\\
ET3tmb&=[[0,NW,1][1,NW,0][0,NW,1]=>[1,NW,0][0,NW,1][1,NW,0]],\\
ET3tmbR&=[[1,NW,0][0,NW,1][1,NW,0]=>[0,NW,1][1,NW,0][0,NW,1]].
\end{eqnarray*}

\begin{figure}
  \centering
  \includegraphics[scale=0.5]{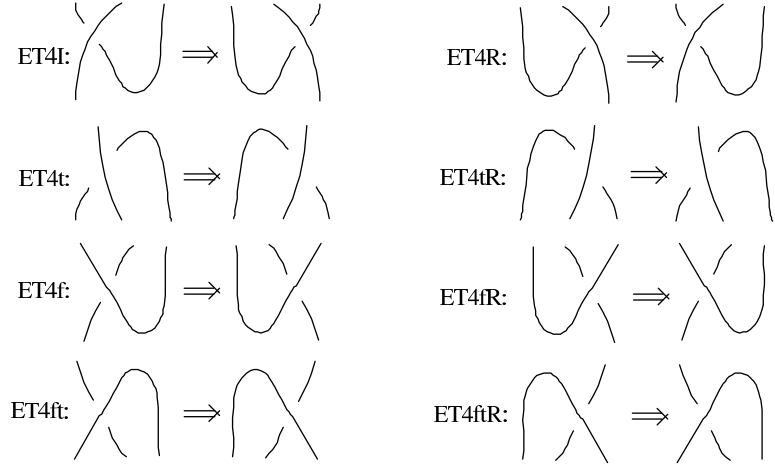}
  \caption{Elementary Transitions of Type ET4}\label{fig:ET4}
\end{figure}

There are 8 elementary transitions, associated to the Yetter $"\Psi"$ move,
tabulated as follows (and pictured in Figure \ref{fig:ET4}):
\begin{eqnarray*}
ET4I&=[[1,Cup,0][0,NE,1]=>[0,Cup,1][1,NW,0]],\\
ET4R&=[[0,Cup,1][1,NW,0]=>[1,Cup,0][0,NE,1]],\\
ET4t&=[[0,NW,1][1,Cap,0]=>[1,NE,0][0,Cap,1]],\\
ET4tR&=[[1,NE,0][0,Cap,1]=>[0,NW,1][1,Cap,0]],\\
ET4f&=[[1,Cup,0][0,NW,1]=>[0,Cup,1][1,NE,0]],\\
ET4fR&=[[0,Cup,1][1,NE,0]=>[1,Cup,0][0,NW,1]],\\
ET4ft&=[[0,NE,1][1,Cap,0]=>[1,NW,0][0,Cap,1]],\\
ET4ftR&=[[1,NW,0][0,Cap,1]=>[0,NE,1][1,Cap,0]],\\
\end{eqnarray*}

\begin{figure}
  \centering
  \includegraphics[scale=0.5]{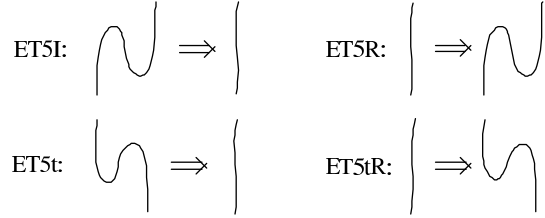}
  \caption{Elementary Transitions of Type ET5}\label{fig:ET5}
\end{figure}

There are 4 elementary transitions, associated to the second
Yetter "switchback" move, tabulated as follows (and pictured in
Figure \ref{fig:ET5}):
\begin{eqnarray*}
ET5I&=[[1,Cup,0][0,Cap,1]=>1_1],\\
ET5R&=[1_1=>[1,Cup,0][0,Cap,1]],\\
ET5t&=[[0,Cup,1][1,Cap,0]=>1_1],\\
ET5tR&=[1_1=>[0,Cup,1][1,Cap,0]]
\end{eqnarray*}

\begin{figure}
  \centering
  \includegraphics[scale=0.5]{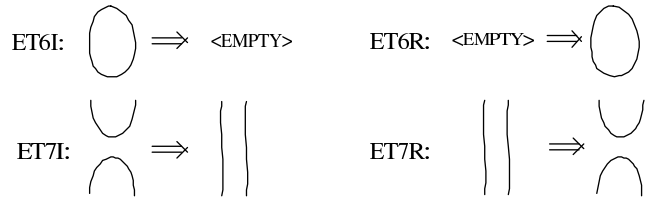}
  \caption{Elementary Transitions of Type ET6 and ET7}\label{fig:ET6_7}
\end{figure}

The next two types of elementary transition, seem to have no
direct analog in the relations (`moves') of the lower-dimensional
theory. The first is sometimes called the "birth-death"
transition; it comes in precisely two flavors, tabulated as
follows, and pictured in Figure \ref{fig:ET6_7}:
\begin{eqnarray*}
ET6I&=[[0,Cup,0][0,Cap,0]=>1_0],\\
ET6R&=[1_0=>[0,Cup,0][0,Cap,0]]
\end{eqnarray*}
The next is the "surgery" or "syzygy" transition; it comes in the two flavors
tabulated below, and also pictured in Figure \ref{fig:ET6_7}:
\begin{eqnarray*}
ET7I&=[[0,Cap,0][0,Cup,0]=>1_2],\\
ET7R&=[1_2=>[0,Cap,0][0,Cup,0]]
\end{eqnarray*}

(Perhaps $ET6$ may be thought of as a higher-dimensional analogue of $Cup,Cap$, and
$ET7$ of $NE,NW$.)
\smallskip

There still remain two types of elementary transitions to discuss,
namely $ET0$ and $ET8$. These are of a more formal nature, and may
be regarded as the elementary transitions analogous to "move 0" in
our earlier discussion associated with Theorem \ref{prop:Yetter}.
(We need to be careful to get these two types ET0 and ET8
straight---especially for the rather difficult proofs of
Props.\ref{prop:semiRespects }
 and \ref{prop:semi17to20 } below---
though in fact these two types play no role in the linear
equations used later to determine amplitude assignments. Like
"move 0", they sometimes seem like more trouble than they are
worth...)
\smallskip

\begin{figure}
  \centering
  \includegraphics[scale=0.5]{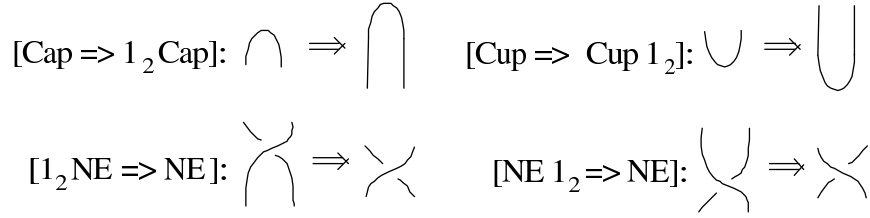}
  \caption{Four Elementary Transitions of Type ET0}\label{fig:ET0}
\end{figure}

The elementary transitions of type ET0, relate to the way the
framed events $1_n:n->n$ play the role (in the tangle category
approach to knotted circles in $\RR^3$) of identity
morphisms.Thus, if $\pi:m->n$ is an elementary event with source
$m$ and target $n$, then associated to $\pi$ we have the 4
elementary transitions
$$ET0\pi l:=[\pi=>1_m \pi], ET0\pi lR:= [1_m \pi=>\pi],ET0\pi r:=[\pi=>
\pi 1_n] $$ and $$ ET0\pi rR=[\pi 1_n=>\pi]\;.  $$ Since $\pi$ can
be any one of 4 elementary events, we thus obtain 16 elementary
transitions of type ET0, of which four are pictured in Figure
\ref{fig:ET0}. Each of these modifies a still, by inserting or
removing
 an appropriate collection of vertical lines.
\smallskip

\begin{figure}
\centering
\includegraphics[scale=0.6]{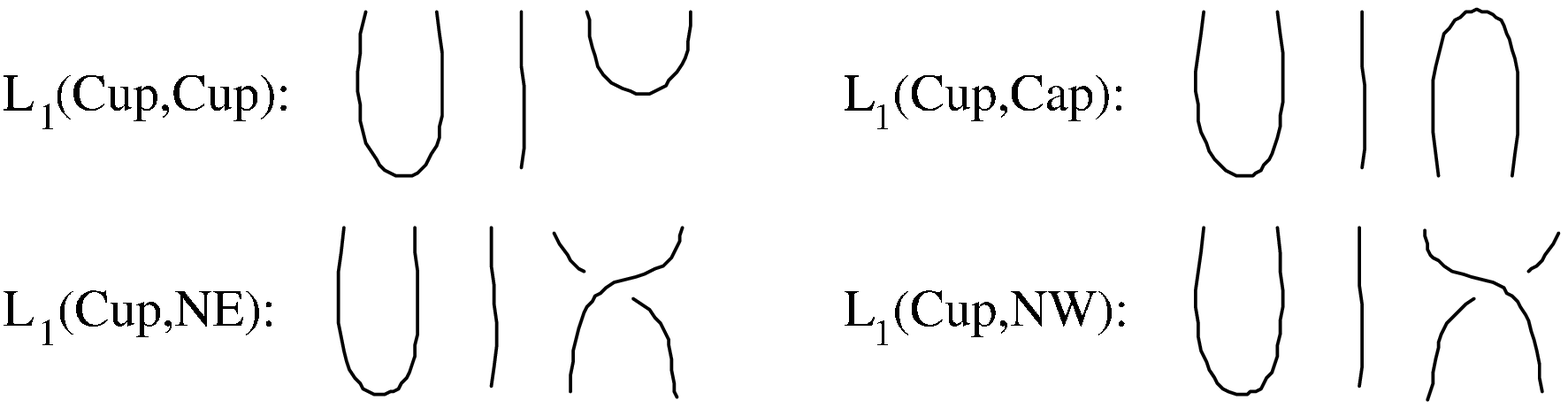}
\caption{\label{fig:L1Cup}$L_1(Cup,*)$}
\end{figure}

Finally, there are an infinite
number of elementary transitions of type ET8, which we next define.\\
Let $E,E'$ be elementary states---i.e. elements (possibly equal) of
$$\mathcal{E}=\{Cup,Cap,NE,NW\}\;\;,$$
and let $n$ be a natural number. Associated to such data, there
 is a law in the category of tangles,
which is roughly speaking a commutative law---it asserts the
isotopy of two "products" $L_n(E,E')$ and $U_n(E,E')$ (each of two
elementary events, separated by $n$ vertical strings) which are
next to be defined---and so in the
category of movies we get corresponding elementary transitions of type ET8. In more detail:\\
In the first place, we denote by $L_n(E,E')$,where $n$ is any natural number,
the still which is defined---roughly speaking---by placing $E$
{\bf below} and {\bf to the left} of $E'$, with
$n$ vertical strands between.). More precisely, here is a case-by-case
definition of $L_n(E,E')$:\\
(cf.Figure \ref{fig:L1Cup} ,where the first four of these are pictured
with $n=1$.)\\
\smallskip

\noindent $L_n(Cup,Cup):= [0,Cup,n][n+2,Cup,0],\\
L_n(Cup,Cap):=[0,Cup,n+2][n+2,Cap,0],$\\
$L_n(Cup,NE):=[0,Cup,n+2][n+2,NE,0],\\
L_n(Cup,NW):=[0,Cup,n+2][n+2,NW,0],$\\
$L_n(Cap,Cup):=[0,Cap,n][n,Cup,0],\\ L_n(Cap,Cap):=[0,Cap,n+2][n,Cap,0],$\\
$L_n(Cap,NE):=[0,Cap,n+2][n,NE,0], \\
L_n(Cap,NW):=[0,Cap,n+2][n,NW,0],$\\
$L_n(NE,Cup):=[0,NE,n][n+2,Cup,0],\\
L_n(NE,Cap)=[0,NE,n+2][n+2,Cap,0],$\\
$L_n(NE,NE):=[0,NE,n+2][n+2,NE,0],\\
L_n(NE,NW):=[0,NE,n+2][n+2,NW,0],$\\
$L_n(NW,Cup):=[0,NW,n][n+2,Cup,0],\\
L_n(NW,Cap)=[0,NW,n+2][n+2,Cap,0],$\\
$L_n(NW,NE):=[0,NW,n+2][n+2,NE,0],\\
L_n(NW,NW):=[0,NW,n+2][n+2,NW,0].$\\
These 16 separate formulas are all subsumed in:\\
$$L_n(E,E'):=[0,E,n+source(E')][n+target(E),E',0]\mbox{ for } E,E'\in  \mathcal{E}\;.$$
\smallskip

\begin{figure}
\centering
\includegraphics[scale=0.6]{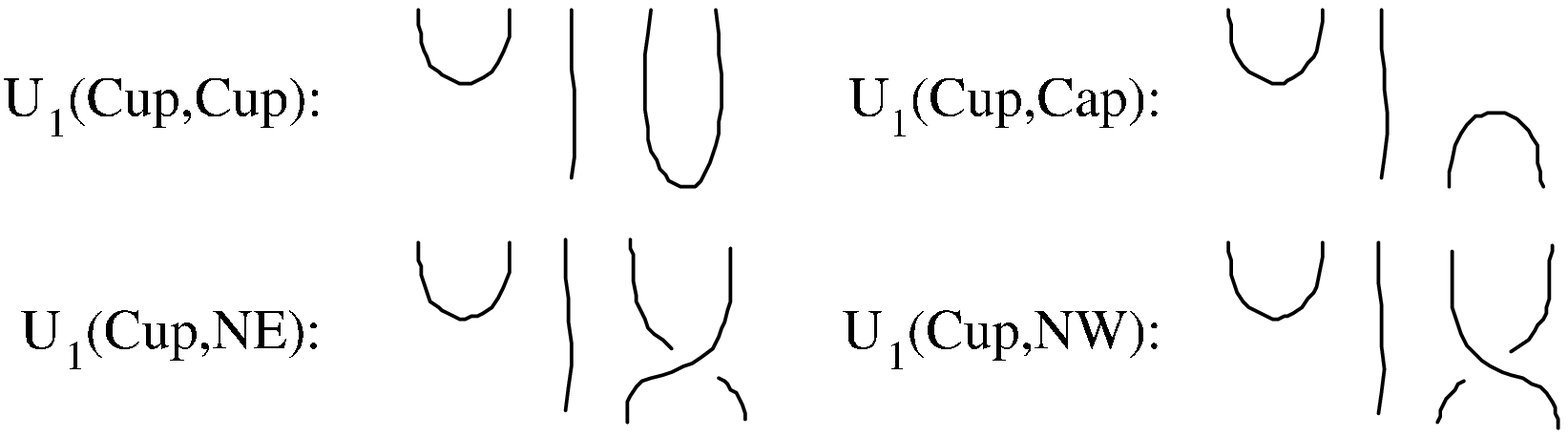}
\caption{\label{fig:U1Cup}$U_1(Cup,*)$}
\end{figure}

Similarly, $U_n(E,E')$is defined, roughly speaking, to be the
still which contains $E$ {\bf above} and to the {\bf left} of
$E'$, with $n$ vertical strands between. (`L' for lower,`U' for
upper...) More precisely, here is a case-by-case definition:
(cf. also Figure \ref{fig:U1Cup}.)\\
$U_n(Cup,Cup):=[n,Cup,0][0,Cup,n+2],\\
U_n(Cup,Cap):=[n,Cap,0][0,Cup,n],$\\
$U_n(Cup,NE):=[n,NE,0][0,Cup,n+2],\\
U_n(Cup,NW):=[n,NW,0][0,Cup,n+2],$\\
$U_n(Cap,Cup):=[n+2,Cup,0][0,Cap,n+2],\\
U_n(Cap,Cap):=[n+2,Cap,0][0,Cap,n],$\\
$U_n(Cap,NE):=[n+2,NE,0][0,Cap,n+2],\\
U_n(Cap,NW):=[n+2,NW,0][0,Cap,n+2],$\\
$U_n(NE,Cup):=[n+2,Cup,0][0,NE,n+2],\\
U_n(NE,Cap):=[n+2,Cap,0][0,NE,n],$\\
$U_n(NE,NE):=[n+2,NE,0][0,NE,n+2],\\
U_n(NE,NW):=[n+2,NW,0][0,NE,n+2],$\\
$U_n(NW,Cup):=[n+2,Cup,0][0,NW,n+2],\\
U_n(NW,Cap):=[n+2,Cap,0][0,NW,n],$\\
$U_n(NW,NE):=[n+2,NE,0][0,NW,n+2],\\
U_n(NW,NW):=[n+2,NW,0][0,NW,n+2].$\\
These 16 separate formulas are all subsumed in:\\
$$U_n(E,E'):=[n+source(E),E',0][0,E,n+target(E')]\mbox{ for } E,E'\in  \mathcal{E}\;.$$
\smallskip

\begin{figure}
\centering
\includegraphics[scale=0.5]{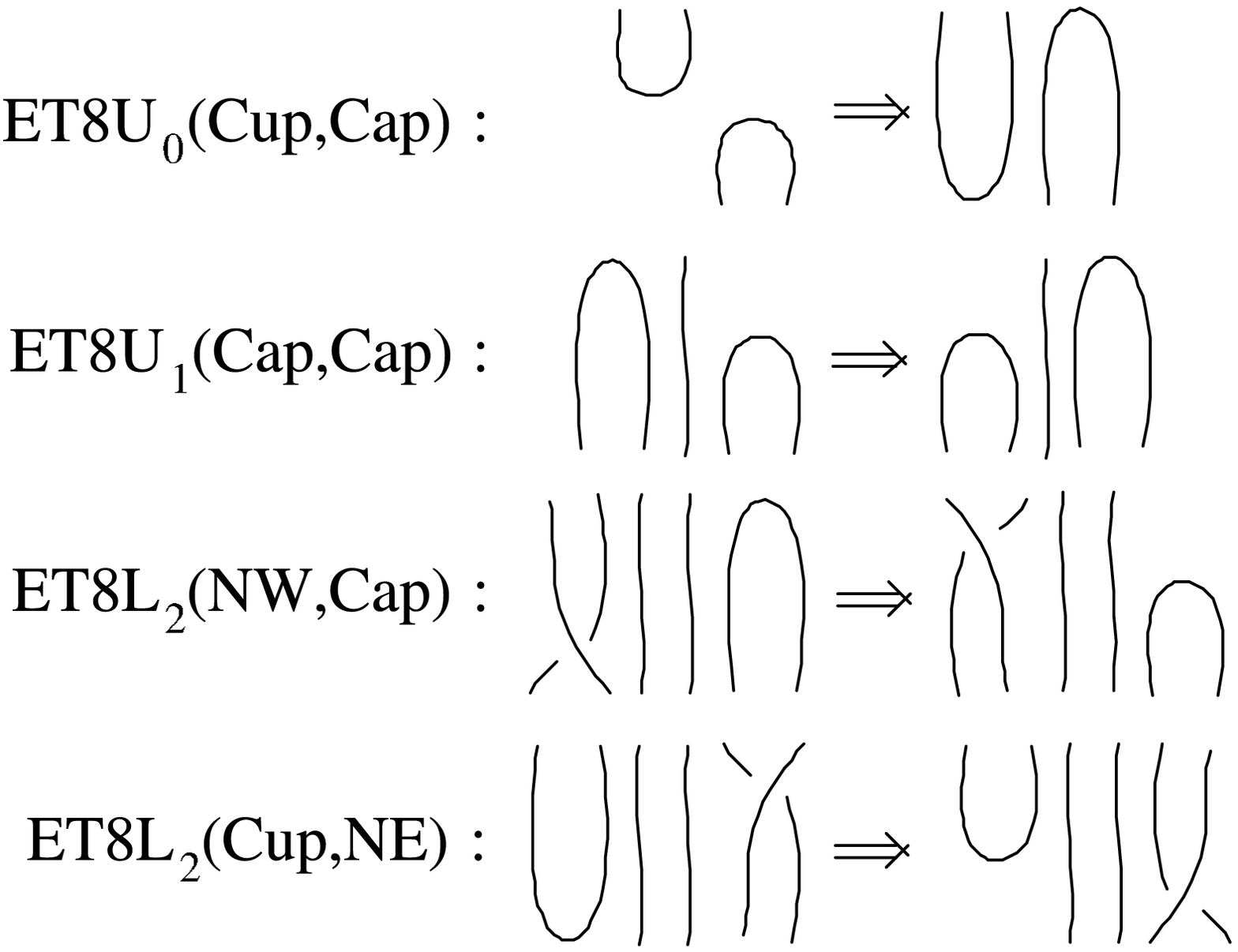}
\caption{\label{fig:et8} Some Elementary Transitions of type
$ET8$}
\end{figure}

Finally, with every ordered pair $E,E'$ of elementary events, and every
natural number $n$, we associate two
elementary transitions of type ET8, as follows:
$$ET8L_n(E,E'):=[L_n(E,E')=>U_n(E,E')],\;ET8U_n(E,E'):=[U_n(E,E')=>L_n(E,E')]$$
The authors refuse to draw here pictures for all 32 families of these, but perhaps
the 6 pictures in Fig. \ref{fig:et8} will help convey the general idea.
\smallskip

We shall refer to elementary transitions of type $ET8$ as
`commutation transitions'. In [CS](cf.p.34) these transitions (and
the corresponding isomorphisms in the lower-dimensional category)
are referred to as `exchanging the levels of distant critical
points'---perhaps in analogy to the axiom in relativistic quantum
field theory, according to which the commutators of amplitudes of
events `distant' (in the sense they are separated by a space-like
interval), are given by suitable Dirac deltas. Of course, in the
present case, what are `commuted' are \underline{{\bf adjacent}}
framed events in the $\tau$ direction,{\bf distant} in the sense
they may be separated by one or more vertical lines.
\smallskip

\begin{tabular}{l}
\hline {\bf This completes our construction and labelling of the
set $\M{ET}$.}\\
 \hline
\end{tabular}
\smallskip

We conclude this subsection by making a few elementary
observations, on some of the combinatorial structure on this set
$\M{ET}$, which will be useful
for the constructions in the present paper:\\

Here is a table, showing how the elementary transitions are divided into types and flavors:
\begin{center}
\begin{tabular}{l|c|c|c|c|c|c|c|c|c}
Type & 0 & 1 & 2 & 3 & 4 & 5 & 6 & 7 & 8\\
\hline
Number of flavors&16 & 8 & 4 & 12 & 8 & 4 & 2 & 2 & $\infty$
\end{tabular}
\end{center}
And, just to be totally explicit, an elementary transition listed
above as of type $ETi$ ($0\leq i\leq 8$), will also be said to be
{\bf of type} $i$. Also, if $E$ is an elementary transition if
type $i$, we shall sometimes write $t(E)=i$.
\smallskip

Let us next note  the following three basic symmetries of
elementary transitions (whose definition is based on the earlier
symmetries on stills given by Defs.\ref{still:f} and
\ref{still:t}):

\begin{Def}\label{Def:R:ET}
 \underline{\bf THE SYMMETRY R:}\\
Let $\mathcal{E}=[S=>T]$ be an elementary transition; then by
$R\cdot \mathcal{E}$ will be meant the elementary transition
$[T=>S]$.
\end{Def}
\begin{Def}\label{Def:t:ET}
\underline{\bf THE SYMMETRY t:}\\
Let $\mathcal{E}=[S=>T]$ be an elementary transition; then by
$t\cdot \mathcal{E}$ will be meant the elementary transition
$[(t\cdot S)=>(t\cdot T)]$
\end{Def}
\begin{Def}\label{Def:f:ET}
 \underline{\bf THE SYMMETRY f:}\\
Let $\mathcal{E}=[S=>T]$ be an elementary transition; then by
$f\cdot \mathcal{E}$ will be meant the elementary transition
$[(f\cdot S)=>(f\cdot T)]$
\end{Def}
Note that these three symmetries on elementary transitions are
well-defined, i.e. examination of the above collection of pictures
shows that, if $[S=>T]$ is an elementary transition, so are
$[T=>S]$, $[t\cdot T=> t\cdot S]$ and $[f\cdot S=>f\cdot T]$.
Also, these three symmetries commute pairwise, with each having
square equal to the identity, and so give rise to an action on
$\mathcal{ET}$ of the finite Abelian group
$$
\mathcal{SYM}=\{I,R,t,tR,f,fR,ft,ftR\} \cong \ZZ/(2)\times \ZZ/(2)\times \ZZ/(2)
$$
of order 8, generated by $R,t$ and $f$.
\smallskip

\noindent {\bf REMARK:} It might seem appropriate to add a fourth
symmetry---call it $m$--- consisting in reflecting in a vertical
midline in each still (without changing the assigned over- and
under-crossings, or the movie-ordering of stills.) However,
observe that $mt=tm$ consists in rotating by $180^{\circ}$ in each
still, thus resulting in a movie easily shown to be isotopic to
the original movie. {\bf For this reason, $m$ will not be needed,
and will not occur, in our further discussions of symmetries.}
(With two exceptions: in the proof of Prop.\ref{prop:semi17to20 }
below, and also in the proof of Prop.\ref{prop:Assoc31}, while the
symmetry $m$ is not absolutely needed in the proof itself, it
clarifies and simplifies the proof.)
\medskip

Utilizing these definitions of symmetries, let us now clarify the
rationale
 behind our notation above for the elementary
transitions---it is important to get this straight, since it
underlies our later notation for the 102 variables (called
`amplitude-parameters'), in the system of linear equations we
shall need to solve.

Consider first the 8 flavors of $ET1$ (cf.Fig.\ref{fig:ET1}).The group
$$\mathcal{SYM}=\{I,R,t,tR,f,fR,ft,ftR\}$$
acts freely on this set of 8 flavors. We choose (quite
arbitrarily) one of these 8 to be the `origin' ET1I, and then for
any $\pi$ in $\mathcal{SYM}$, we denote $\pi \cdot ET1I$ by
$ET1\pi$.

The 4 flavors of $ET2$ (cf.Fig.\ref{fig:ET2}) are labelled similarly, with this minor
variation: Note first, the symmetry $t$ in $\mathcal{SYM}$ fixes each of the 4 flavors of $ET2$.
We choose (again quite arbitrarily) one of these four to be $ET2I$, and then
the other 3, labelled respectively $ET2R,ET2f,ET2fR$ are obtained from our choice of
$ET2I$ by acting on it, respectively, by $R,f,fR$.

The labelling for the flavors of $ET4,5,6,7$ (as supplied in Figures
\ref{fig:ET4}, \ref{fig:ET5}
and \ref{fig:ET6_7}) should now be self-explanatory.

Our labelling for the 12 flavors of $ET3$ (cf. Fig. \ref{fig:ET3}) involves
the following different ideas:
\smallskip

Consider any one of the 12 pictures in Fig. \ref{fig:ET3}. It is
of the form $S=>T$, where each of the stills $S$ and $T$ involves
3 strings and 3 intersections. Of the 3 strings which make up $S$
(in the picture chosen), one lies above the other two, one lies
above another string and below the other one, and one lies beneath
the other two--- label these t, m,b respectively (for top, middle,
bottom).Let us (arbitrarily) choose the top (rather than the
bottom) of $S$ to label with these 3 letters---and note that we
have labelled the picture, by this permutation of b,t,m  {\bf OR}
by this permutation followed by R.(see Fig. \ref{fig:ET3} for the
two transitions corresponding to the permutation btm) For a given
permutation $\pi$ of $\{b,t,m\}$ correspond {\bf TWO} elementary
transitions of type $ET3$, one labelled $ET3\pi$ and the other
labelled $ET\pi R$ ,where the former has the appearance
$$ET3\pi=[[0,*,1][1,*,0][0,*,1]=>[1,*,0][0,*,1][1,*,0]]$$
(with each * replaced by a suitable value of $NE$ or $NW$) and the latter
has the appearance
$$ET\pi R=[[1,*,0][0,*,1][1,*,0]=>[0,*,1][1,*,0][0,*,1]]\;.$$
(As a check on these computations, note that each of these 12
elementary transitions $S=>T$ of  type ET3 is a palindrome, in the
sense that if we read $S$ in our present notation from right to
left, we obtain the notation for $T$.)
\smallskip

\begin{Prop}
Let $E$ be an elementary transition, with source $S$ and target $T$; then
$$source(S)=source(T)\mbox{ and } target(S)=target(T)\;.$$
\end{Prop}
\noindent {\bf PROOF:} Direct examination of all pictures in this
section, together with Def.\ref{Def:still}.
\begin{Def}\label{Def:in:out}
Let $E:S=>T$ be an elementary transformation; then by $in(E)$ will be meant the common value
of $source(S)$ and $source(T)$, while by $out(E)$ will be meant the common value of
$target(S)$ and $target(T)$.
\end{Def}

\noindent {\bf EXAMPLE:} Examination of Fig.\ref{fig:ET1} shows
that
$$in(ET1I)= 0,\mbox{ and }out(ET1I)=2            $$

Because of the great importance of [LL] in the theory of
2-knots, and also because of our need below to refer to [LL]
in the discussion of movie-moves,
it seems worthwhile to conclude the present section,
by presenting the following table, which shows how
 to translate between the notation for elementary transitions explained above,
and the notation for elementary transitions explained in [LL](p.43 and 44) and
[BL](p.46):\\
\begin{center}
\begin{tabular}{|l|c|c|c|c|c|c|c|c|c|}
ET notation & 1tR & 2R & 3tmb & 3tbm & 3mbt & 4ftR & 5I & 6R & 7I\\
\hline
(LL) notation & W & $\Pi$ & $S_0$ & $S_1$ & $S_2$ & H & T & I & E
\end{tabular}
\end{center}
Also, in [LL] and [BL], the symbols
$$\overline{\alpha}, \alpha^{\ast}, \alpha^{\dag}$$
(where $\alpha$ is an elementary transition) are used for the symmetries which in the
notation of the present paper, translate (respectively) to
$$f\cdot \alpha, R\cdot \alpha,Rt\cdot \alpha \;.$$
Hence, in particular,
$$\overline{S_0}=ET3bmt,\overline{S_1}=ET3btm,\overline{S_2}=ET3mtb$$
The composition of stills in [LL] and [BL] occurs (alas!) in the opposite order to that
used in the present paper.

Finally, in loc.cit., the notation $N_{Y_{m,n},Z_{i,j}}$, where
$Y,Z$ denote elementary events (i.e. elements of $\{Cup,Cap,NE,NW
\}$) and $m,n,i,j$ are natural numbers, denotes
>>>>>>>>>>>>>>>>>>>>>>>>>>>>
TO BE CONTINUED
>>>>>>>>>>>>>>>>>>>>>>>>>>>>>>>>>>>>>>>>>>>>>>>>>
 the elementary transition of type 8 (illustrated
in Figure \ref{fig:et8}):
$$m\bullet L_N(Y,Z)\bullet j => m\bullet U_N(Y,Z)\bullet j $$
i.e.\\
$s[m,Y,N+source(Z)+j][m+target(Y)+N,Z,j]s=>8$\\
$f[m+source(Y)+N,Z,j][m,Y,N+target(Z)+j]f$\\
where N is determined by
\begin{equation}\label{eqn:N}
N=n-j-target(Z)=i-m-source(Y)\;.
\end{equation}
Note that (\ref{eqn:N}) has the  following three consequences:
$$n-j\geq target(Z),\;i-m\geq source(Y)$$
(the second of which is mentioned in [BL] and [LL]), and
$$n-j-target(Z)=i-m-source(Y)\;.$$
---which are the necessary and sufficient conditions for the symbol $N_{Y_{m,n},Z_{i,j}}$
to make sense. Mention THIS paper is $sl_q(2)$-based.
\subsection{Flickers and Movies}\label{S:movies}
We next shall define the combinatorial concept (due in its final
form to Carter,Rieger and Saito,[CRS]), of a \underline{{\bf
`movie'}}. Very roughly speaking, this is a sequence of stills,
each obtained from the preceding one by a permissible kind of
transition, next to be explained, called a
\underline{{\bf`flicker'}} in the present paper. More precisely:
\medskip

Recall first that, in the lower-dimensional theory sketched in $\S
\ref{SS:Stills}$, a \emph{`framed event'} $[m,E,n]$ was obtained
by enhancing an \emph{`elementary event'} $E$ by $m$ vertical
strings on the left, and $n$ on the right
(cf.Fig.\ref{fig:events}). A \emph{`flicker'} is the
next-higher-dimensional analog of this: it is obtained from an
elementary transition $\mathcal{E}:U\Longrightarrow V$ upon
enhancing $\mathcal{E}$ by a `bottom still' $B$ below, $m$
vertical strings to the left of $\mathcal{E}$ and $n$ to the
right, and then a `top still' T above. (cf. Figure
\ref{fig:flicker} for a schematic illustration). In more detail
(cf. Def. 5.):
\begin{figure}
  \centering
  \includegraphics[scale=0.8]{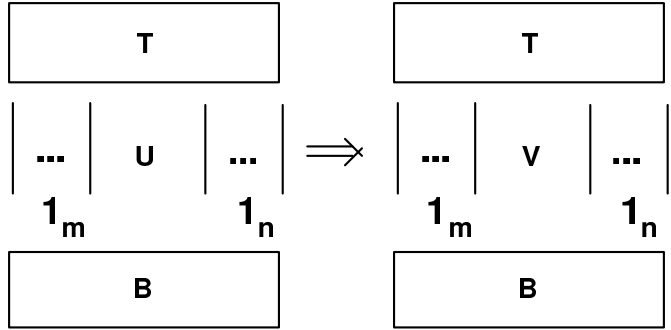}
  \caption{Flicker $\mathcal{F} = (U \Longrightarrow V, B, m, n, T)$ }\label{fig:flicker}
\end{figure}

\begin{Def}\label{defFlicker}
By a \underline{\bf {"flicker"}}. will be meant one of two types of object:\\
{\bf TYPE ONE:} an ordered 5-tuple
$$\mathcal{F}=(\mathcal{E},B,m,n,T) \eqno(2)$$
where $\mathcal{E}$ is an elementary transition, $B$ and $T$ are stills (either
or both of which may be empty), $m$ and $n$ are
natural numbers(which may be 0), such that the two following relations are satisfied:
$$target(B)=m+in(\mathcal{E})+n \eqno(3a)$$
$$source(T)=m+out(\mathcal{E})+n \eqno(3b)$$
(with $in(\mathcal{E}), out(\mathcal{E})$ as defined in Def.\ref{Def:in:out}.)
This being the case, we then define the \underline{\emph {source}} and
\underline{\emph {target}} of the  flicker $\mathcal{F}$, to be the stills
$$source(\mathcal{F})=T\circ(m\bullet source(\mathcal{E})\bullet n)\circ B \eqno(4a)$$
$$target(\mathcal{F})=T\circ(m\bullet target(\mathcal{E})\bullet n)\circ B \eqno(4b)$$
(cf. Defs.\ref{def:stillComp} and \ref{def:circ}), and we then write
$$\mathcal{F}:source(\mathcal{F})\Longrightarrow target(\mathcal{F})$$
If $\mathcal{E}$ is of type $i$ (so $0\leq i\leq 8$) then we shall
also say that the flicker $\mathcal{F}$ is \underline{\emph{of
type i}}.
\smallskip

\noindent {\bf TYPE TWO:} a symbol of the form $$1_S \eqno(5)$$
where $S$ is a still; both the source and target of this flicker
are defined to be $S$.
\end{Def}
\medskip

\begin{figure}
  \centering
  \includegraphics[scale=0.5]{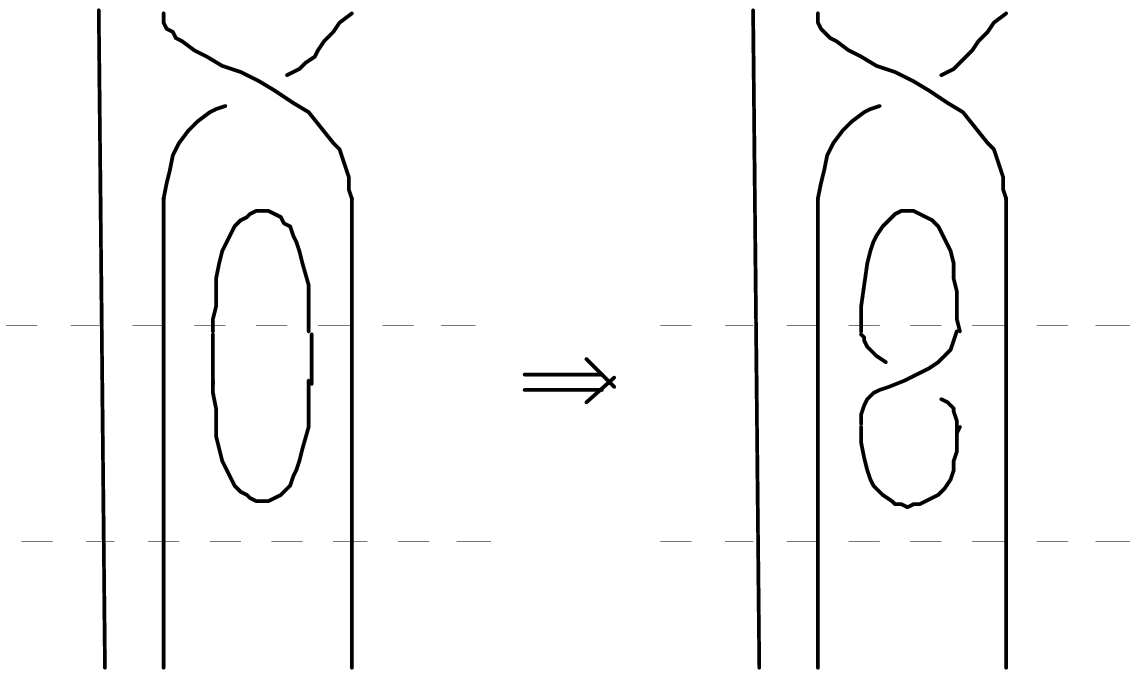}
  \caption{$(\mathcal{E},1_3,2,1,[2,Cap,1][1,NW,0])$}\label{figFlickerEx}
\end{figure}

\noindent {\bf EXAMPLE:} Let $\mathcal{E}=ET1R$ (cf.Fig.
\ref{fig:ET1}), then the flicker
$$(\mathcal{E},1_3,2,1,[2,Cap,1][1,NW,0])$$
is represented by Figure \ref{figFlickerEx}.

\medskip
\begin{Prop+Def}
Let $\mathcal{F}$ denote a flicker, with source $S$ and target $T$; then
$$source(S)=source(T) \mbox{ and } target(S)=target(T)\;.$$
By $in(\mathcal{F})$ will be meant the common value
of $source(S)$ and $source(T)$, while by $out(\mathcal{F})$ will be meant the common value of
$target(S)$ and $target(T)$.
\end{Prop+Def}
\noindent {\bf PROOF:} We must consider two cases:\\
If $$\mathcal{F}=(\mathcal{E},B,m,n,T) $$ is of type (2) in Def.\ref{defFlicker}, then
by (4a) and (4b), $source(\mathcal{F})$ and $target(\mathcal{F})$ both have sources equal to
$source(B)$, and targets equal to $target(T)$

\noindent On the other hand, if $\mathcal{F}$ has the form $1_U$,
then $S=T=U$.

 Q.E.D.
\medskip

\begin{Def}\label{def:circFlic}
Let $\mathcal{F}$ be a flicker, and let $p,q$ be natural numbers; then
by $p\bullet \mathcal{F}\bullet q$
will be meant the flicker defined as follows:\\
Case 1) If $\mathcal{F}=(\mathcal{E},B,m,n,T) $ is of type (2) in Def.\ref{defFlicker}, then
(cf. Def.\ref{def:circ})
$$p\bullet \mathcal{F}\bullet q:=(\mathcal{E},p\bullet B\bullet q,p+m,q+n,p\bullet T\bullet q)$$
Case 2) If $\mathcal{F}=1_U$ then
$$p\bullet \mathcal{F}\bullet q:=1_{p\bullet U\bullet q}$$
\end{Def}
We next extend to flickers, the action of the group
$\mathcal{SYM}$ on elementary transitions, given in
Defs.\ref{Def:R:ET}, \ref{Def:t:ET} and \ref{Def:f:ET}. It
suffices to define the action on flickers, of the commuting
generators $R,t$ and $f$:
\begin{Def}\label{def:symFlics}
Let $\mathcal{F}$ be a flicker; then we
define $R\cdot \mathcal{F}$, $t\cdot \mathcal{F}$ and $f\cdot \mathcal{F}$ as follows:\\
Case 1) If $\mathcal{F}=(\mathcal{E},B,m,n,T) $ is of type (2) in Def.\ref{defFlicker}, then
we set\\
$R\cdot \mathcal{F}:=(R\cdot \mathcal{E},B,m,n,T)$,\\
$t\cdot \mathcal{F}:=(t\cdot \mathcal{E},t\cdot T,m,n,t\cdot B)$,\\
$f\cdot \mathcal{F}:=(f\cdot \mathcal{E},f\cdot B,m,n,f\cdot T) $.\\
Case 2) If $\mathcal{F}=1_U$ (for some still $U$) then we set
$$R\cdot \mathcal{F}:=\mathcal{F};\; \; t\cdot \mathcal{F}:=
1_{t\cdot U}\mbox{ and }f\cdot \mathcal{F}:=1_{f\cdot U}\;.$$
\end{Def}

\medskip

The  following definition should be compared with Def.\ref{Def:still}.
\begin{Def}\label{DefMovie}
A \underline{\emph {Carter-Rieger-Saito movie}} (or simply
\underline{\emph {movie}}, when there is no danger of confusion
with other senses of `movie' \footnote{Some different senses of
the word `movie' may be found, for instance, in [CS1], which
involves a sense proposed earlier by Carter and Saito involving
only one rather than two time-orientations, in ([CS2], \S 1.4])
and in [Ka], Chapter 8.}) is defined to be a finite non-empty
sequence of flickers,
$$M=(F_1,...,F_s)\eqno (6)$$
such that the flickers $F_i$ satisfy
$$target(F_i)=source(F_{i+1})\mbox{ for }1\leq i<s\;.$$
By the \underline{\emph {source}} of a movie M, will be meant the source
of its first flicker, i.e. $source(F_1)$ if (6) holds, while
its \underline{\emph {target}}
is defined to be the target of its last flicker, i.e. $target(F_s)$.
If the stills $S,T$ are the source (~resp. target) of the movie $M$, we write
$M: S\Longrightarrow T$.
We  shall often use a multiplicative notation to express (6), writing instead
$$M=F_1F_2...F_s \;.$$
Finally, we define $in(M)$ to be the common value of all $in(F_i)$, and
$out(M)$ to be the common value of all $out(F_i)$.
\end{Def}
\medskip

As elucidated in [BL] and [LL], movies take part in an intricate structure,
a small portion of which which has already been defined above. For the
purposes of the present paper, we shall only need the further structure on movies given by
the following four definitions (which will play a role in the discussion below (\S
\ref{SS:moviemoves}) of movie-moves).
\begin{Def}\label{moviesComposite}
Let $$M=F_1\cdots F_s,\; M'=F'_1\cdots F'_t$$ be two movies. We
shall say that $M$ \underline{{\bf is composable with }}$M'$ if
the target of $M$ equals the source of $M'$, ie if
$$target(F_s)=source(F'_1)$$---in which case
 the \underline{{\bf composite}} movie $M\circ M'$  is defined to be the movie
consisting of the sequence
$$(F_1,\cdots,F_s,F'_1,\cdots,F_t')$$
of $s+t$ flickers.
\end{Def}
{\bf NOTE:} In [BL] this is called the {\bf vertical composite};
movies modulo movie-moves have a structure of 2-category, which
means there is also defined another operation on movies (modulo
movie-moves), the {\bf horizontal composite}, which will not be
utilized in this paper.We here choose the ordering-convention for
vertical composition, which coincides with that in [BL] and [LL].
\begin{Def}\label{bulletMovie}
Let $M=(F_1,\cdots,F_s)$ be a movie, and let $m,n$ be natural numbers.\\
Then (cf.Def.\ref{def:circFlic}) we define $m\bullet M \bullet n$ to be the movie
$$(m\bullet F_1 \bullet n,\cdots,m\bullet F_s \bullet n)$$
\end{Def}
\begin{Def}
Let
$$M=(F_1,\cdots,F_s)$$
be a movie, with
$$in(M)=m, out(M)=n$$
(so each still $F_i$ of $M$ has $m$ strings leading into the bottom, and
$n$ strings leading out the top.) Let $S,T$ be stills, such that
$$target(S)=m, source(T)=n\;;$$
then we define $S\circ M\circ T$ to be the movie
$$((S\circ F_1\circ T),\cdots,(S\circ F_s\circ T))  $$
(cf. Def.\ref{def:stillComp})
\end{Def}
Finally, we define the action of the group $\mathcal{SYM}$ on
movies. It suffices to define the action on movies, of the
commuting generators $R,t$ and $f$:
\begin{Def}\label{def:symMovies}
Let
$$M=(F_1,\cdots,F_s)$$
 be a movie; then we define :\\
$R\cdot M:=(R\cdot F_s,\cdots,R\cdot F_1)$\\
$t\cdot M:=(t\cdot F_1,\cdots,t\cdot F_s)$,\\
$f\cdot M:=(f\cdot F_1,\cdots,f\cdot F_s) $.\\
\end{Def}
\subsection{Notation for Stills, Flickers, and Movies }\label{SS:notation}
In order for the reader to utilize the program 2KnotsLib (in
particular, for computing $\M{U}$-regular amplitude-invariants of
movies), we now need to explain the notation used, both for some
of the mathematical arguments in this paper, and also used to
input stills, flickers and movies into that program.
\smallskip

\noindent \underline{{\bf NOTE:}} Unfortunately, in trying to
ensure the two programs C++2KnotsLib and Java2KnotsLib were
written fairly independently (and so could be used to check one
another's results), slightly different methods of inputting data
were developed for these two programs. The notation explained in
this sub-section, is that originally developed for the library
C++2KnotsLib and its driver programs. As indicated above, this
notation will also be used in many of the proofs in the present
paper.
\begin{description}
\item[a)]The other library,
Java2KnotsLib, cheerfully accepts all the notation explained in
this sub-section,
with one exception:\\
what the C++ version denotes by $1_n$, the Java version only
accepts as $[n_1,null,n_2]$ with $n_1+n_2=n$ (case insensitive,
i.e. $Null$ is OK).
\item[b)]The Java version \emph {also} accepts
some other slight variations in this section's notation, as
explained in the ReadMe for the Java version. The C++ version
\emph{only} accepts the notation next to be explained.
\end{description}
\smallskip

The notation for stills used in this paper, will be that explained
in the preceding $\S$ \ref{SS:Stills}
---so, for example, the notation for the still in
Fig.~\ref{fig:still_example} is $[2,NW,0][1,Cap,1][0,NE,0]$. Note
the framed event $1_n$ is written without square brackets. This is
practically the notation used in [CS]---slightly modified to
facilitate inputting via keyboard where the [CS] notation is,
rather, more suitable for coding into \TeX--- e.g. $[2,Cup,0]$
here replaces the Carter-Saito $\cup _{2,0}$. We here translate
their $X,\overline{X}$ by $NW,NE$---this choice of translation
involves an arbitrariness caused by the symmetry $f$. We again
remind the reader, that our time-ordering convention for stills
runs from bottom to top of the page---this is (alas) the opposite
of the convention used in [CS] and in [CRS].

\begin{figure}
\centering
\includegraphics[scale=0.6]{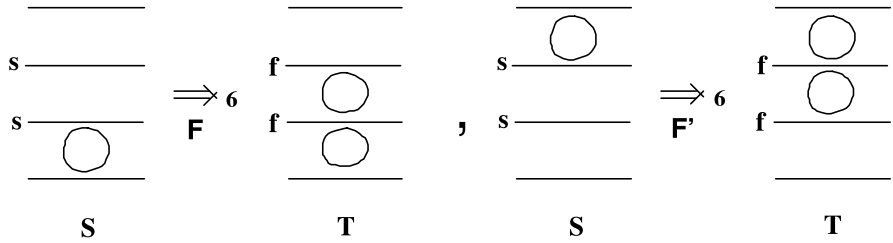}
\caption{\label{fig:S_to_T}Two flickers for same still pair (S,T)}
\end{figure}

We next need a notation for flickers. One unambiguous notation is
furnished by (2) and (5) of Def. \ref{defFlicker}, but we also
need a notation which describes a flicker in terms of its source
and target.
For this purpose, a little care is needed:\\
\noindent {\bf CAUTION:} Let $F$ be a flicker, with source $S$ and
target $T$. As was noted above, if $F$ is an elementary
transition, $F$ is uniquely determined by $S$ and $T$. {\bf This
is not true for flickers in general.} For example, (cf.
[CS],p.35), if
$$S=[0,Cup,0][0,Cap,0]\;, T=[0,Cup,0][0,Cap,0][0,Cup,0][0,Cap,0]\eqno(*)$$
we have $$F,F':S=>T$$ for the  \emph{distinct} flickers
$$F=(ET6R,1_0,0,0,[0,Cup,0][0,Cap,0]),F'=(ET6R,[0,Cup,0][0,Cap,0],0,0,1_0)$$
with notation as in (2) in Def.\ref{defFlicker}. (These are
pictured in Fig. \ref{fig:S_to_T}.) Thus, we need a notation for
flickers, which
 determines explicitly which `path' from $S$ to $T$ is intended.

The notation used here, which we call the `\underline{{\bf
sf-notation}}', solves this problem in the way now to be
described. (Note however that two other, rather different,
solutions may be found in the literature: one in in [CRS], [CS]
and the other in [BL], [LL]).
\smallskip

Consider the flicker represented (in the notation of Def.
\ref{defFlicker}) as
$$\mathcal{F}=(\mathcal{E},B,m,n,T)\;.$$
Let
$$B=F_1\cdots F_u\;,\;T=F'_1\cdots F'_v$$
express the stills $B,T$ as sequences of framed events (here read,
the reader is reminded, from the bottom of the page to the top.)
We wish to express. in our notation, the way this flicker is
derived from its source
$$source(\mathcal{E})=(F_1\cdots F_u)\circ (m\bullet source(\mathcal{E})
 \bullet n)\circ (F'_1\cdots F'_v)$$
and its target
$$target(\mathcal{E})=(F_1\cdots F_u)\circ (m\bullet target(\mathcal{E})
 \bullet n)\circ (F'_1\cdots F'_v)$$
Then here is the `sf-notation' for this flicker F:
\begin{eqnarray}
F &=[(F_1\cdots F_u) s(m\bullet source(\mathcal{E})
 \bullet n) s (F'_1\cdots F'_v)]=> \nonumber\\
  &[(F_1\cdots F_u)f(m\bullet target(\mathcal{E})
 \bullet n)f(F'_1\cdots F'_v)] \nonumber
\end{eqnarray}
Thus, in the source of $F$, s---for `start'---frames the source
\footnote{This source $m\bullet source(\mathcal{E})\bullet n$ may
be `empty', as in the two flickers next described---in which case
it is denoted by $1_0$. Similarly for the target $m\bullet
target(\mathcal{E})\bullet n$.} $m\bullet
source(\mathcal{E})\bullet n$ of the framed elementary transition
$m\bullet \mathcal{E}\bullet n$, while in the target of F, f---for
finish---frames the target $m\bullet target(\mathcal{E})\bullet n$
of $m\bullet \mathcal{E}\bullet n$ . The addition of these two s
and f's to our notation, precisely resolves the ambiguity
explained at the start of this section.\\
{\bf NOTE:} If  the flicker F is of type $i$, there is the option
of expressing this fact by writing $i$ immediately after  the
symbol $=>$. The program works perfectly well if this $i$ is
omitted (it knows how to compute $i$ from the source and target
decorated with s and f), but including the $i$ allows the program
to give a better diagnostic message if the input contains some
error.
\medskip

\begin{figure}
\centering
\includegraphics[scale=0.7]{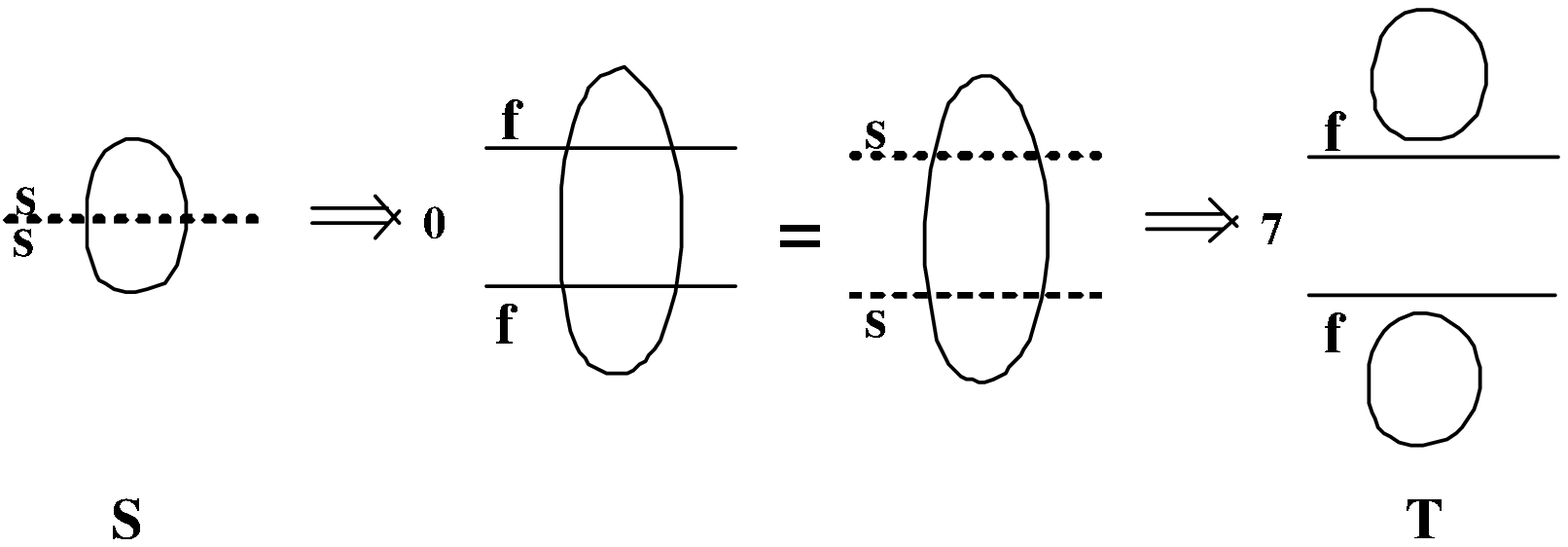}
\caption{\label{fig:MovieA}Movie A }
\end{figure}

For instance, returning to the situation pictured in Fig.
\ref{fig:S_to_T}, the two flickers $F,F':S=>T$ described above,
may now be written (unambiguously) as
\begin{verbatim}
F=[
s1_0s[0,Cup,0][0,Cap,0]=>f[0,Cup,0][0,Cap,0]f[0,Cup,0][0,Cap,0] ]
\end{verbatim}
and
\begin{verbatim}
F'=[
[0,Cup,0][0,Cap,0]s1_0s=>[0,Cup,0][0,Cap,0]f[0,Cup,0][0,Cap,0]f ]
\end{verbatim}
respectively. (There is in each case the option of writing $=>$ as
$=>6$.) For another example, we observe there is yet a third path
(this time involving two steps) from $S$ to $T$, (pictured in Fig.
\ref{fig:MovieA}) consisting of the flicker
\begin{verbatim}
[ [0,Cup,0]ss[0,Cap,0] => [0,Cup,0]f1_2f[0,Cap,0] ]
\end{verbatim}
of type ET0, followed by the flicker
\begin{verbatim}
[ [0,Cup,0]s1_2s[0,Cap,0] =>
[0,Cup,0]f[0,Cap,0][0,Cup,0]f[0,Cap,0]
\end{verbatim}
of type ET7. We may combine these two flickers, to obtain the
following movie
(with source $S$ and target $T$) in sf-notation:\\
\begin{verbatim}
MOVIE A: [0,Cup,0]ss[0,Cap,0] =>0 [0,Cup,0]sf1_2fs[0,Cap,0] =>7
[0,Cup,0]f[0,Cap,0][0,Cup,0]f[0,Cap,0]
\end{verbatim}

We shall see many more examples below, of movies in
sf-notation---for now, let us make just five observations on this
first example in this paper, Movie A, of a movie with more than
one flicker:
\begin{itemize}
\item Just as in this example, so also for the sf-notation for any movie
containing more than one flicker--- the first still contains two
s's, the last still contains 2 f's, and all intermediate stills
contain both. {\bf We make the convention that, for a movie
consisting of a single flicker $1_S:S=>S$,
 we omit both s and f.}
\item In a single still it is  permitted that two s's or two f.s shall be adjacent
---$ss$ is then interpreted as $s1_0s$ (so that the first flicker of
Fig.\ref{fig:S_to_T} could also be input as
\begin{verbatim}
F=[ ss[0,Cup,0][0,Cap,0]=>f[0,Cup,0][0,Cap,0]f[0,Cup,0][0,Cap,0]]
\end{verbatim}
\item  It is permitted
that an s be adjacent to an f---in such a situation, their
ordering is irrelevant.(For example, it would make no difference
if, in the second still in Movie A above, one had instead written
$$[0,Cup,0]fs1_2fs[0,Cap,0]$$
or...)
\item We have used a double arrow $=>$ to indicate a flicker,
motivated by the fact that, in the Baez-Langford theory, flickers
are 2-morphisms. Following this double arrow, we write a number 0
to 8, indicating the {\bf type} of the flicker---in our program,
this is {\bf optional}---these integers may be omitted on some or
all double arrows when inputting a movie.If the option to supply
some or all of these numbers is chosen, it helps the program check
the input for errors.
\item The square brackets [] framing the description of a flicker
are optional---for instance, they are omitted for all the movies
in \S 4.
\end{itemize}
\medskip

Many more examples of this notation will occur in \S
\ref{S:examples}, giving further illustrations of the points just
made. When using the program(s) explained in this paper, the
notation accepted by the program(s) is that just explained---the
further necessary details (how the program input begins and ends a
movie, how comments are written, etc.) will be found carefully
explained for each program in a ReadMe file supplied on the URL
from which the program is to be downloaded.

\subsection{Movie-isotopy and the Carter-Rieger-Saito Movie-move Theorem}\label{SS:movieIsotopy}
\begin{Def}
A movie $M$ will be called \underline{{\bf compact}} if its source and target
are both the empty still, i.e. if
$$source(M)=target(M)=1_0$$
\end{Def}

In [CS] and [CRS] it is explained how every compact movie defines
an isotopy-class of knotted surfaces (in $\RR^4$), and it is
proved that every isotopy-class of knotted surfaces arises this
way. The question immediately suggests itself: Given two compact
movies $M,M'$, when do they give rise to the same isotopy-class of
knotted surfaces? Let us say that two compact movies $M,M'$ are
\underline{{\bf isotopic}} when this is the case.
\smallskip

Note that \emph{a priori}, this relation between movies is (at
least, in terms of its definition) topological rather than purely
combinatorial; the work of Carter, Rieger and Saito converts this
to an equivalent purely combinatorial notion, in a way next to be
described.
\smallskip

 This notion of isotopy between movies has an immediate
generalization to non-compact movies, as follows. Recall that the
composite movie $A\circ B$ is defined whenever $A,B$ are movies
with
$$source(A)=target(B)\;.$$
\begin{Def}\label{def:isotopic}
Let $M,M'$ be movies, not necessarily compact. $M$ and $M'$ are defined to be
\underline{{\bf isotopic}} when the two following statements hold:\\
a) $source(M)=source(M')$ and $target(M)=target(M')$,\\
and\\
b) For all movies $A,B$ of the type
$$A:[1_0]=>source(M), B:target(M)=>[1_0]$$
the movies $A\circ M \circ B$ and $A\circ M' \circ B$ (both necessarily
compact) are isotopic.
\end{Def}

We now turn to the beautiful theorem of Carter, Rieger and Sato
(extending earlier work of Roseman, and of Carter and Saito---cf.
[R],[R2], and [CS1],[CS2])
 which answers (though perhaps non-effectively) the question just raised above.
 The statement of this theorem, involves
 an infinite collection $\mathcal{MM}$ of \underline{{\bf ``movie-moves"}}, first
constructed completely \footnote{This enumeration of
$\mathcal{MM}$ by Carter, Rieger and Saito, presents one possible
difficulty, at least in terms of the explicitness needed for
programming purposes. This difficulty was fixed up in the later
work of Baez and Langford([BL],[LL])---as explained in the
discussion at the beginning of \S \ref{SS:moviemoves}.} by
Carter,Rieger and Saito ([CRS]), which play for knotted surfaces
the same role played in the lower-dimensional theory by: the
Reidemeister and Yetter moves (or more precisely, by the framings
of these.)

A more detailed study of $\mathcal{MM}$ will be the subject of
\S\ref{SS:moviemoves} below.
A few preliminary observations for now:\\
The elements of $\mathcal{MM}$ are ordered pairs of movies (not necessarily
compact). We shall sometimes use the notation $(M==M')$ to denote an ordered pair $(M,M')$
in $\mathcal{MM}$---as this notation is intended to suggest, $M$ is always isotopic to $M'$
when $(M==M')$ is a movie-move. (However, the converse is certainly not true.)

 We are now ready to state (without proof---for proof, the reader is referred to [CS]
and [CRS]) what may be regarded as the fundamental theorem
concerning the matters here under study. (Note this is a
higher-dimensional analogue of the Theorem \ref{prop:Yetter}
asserted above.)
\begin{Th}\label{th:main}
(\emph{Carter-Rieger-Saito}) Let $M$ and $M'$
be two movies; then, $M$ and $M'$ are isotopic, if
and only if it is possible
to go from $M$ to $M'$ by a finite sequence\
$$M=M_0,M_1,\cdots,M_s=M'$$
 of movies, such that, for each $i$ with $1\leq i<s$,
$$(M_i,M_{i+1}) \in \mathcal{MM}\;.$$
\end{Th}

\subsection{The Set $\mathcal{MM}$ of Movie-Moves}\label{SS:moviemoves}
The construction of a set $\mathcal{MM}$ so that Th. \ref{th:main}
shall hold, seems to require much greater effort than the
analogous construction, by explicit enumeration, of $\mathcal{ET}$
in \S \ref{SS:ET} above.
 The first completely explicit
enumeration of the set of Carter-Rieger-Saito movie-moves--- in
particular, sufficiently explicit to enable putting these matters
in a program--- seems (to the present authors) to be that
furnished by Baez and Langford in ([BL], p.47-50) and also in
([LL],p.44-50). The earlier construction given in [CRS], and a bit
later in [CS2], while essentially correct topologically, seems to
present the following combinatorial difficulty:
\medskip

After giving a list of movie-moves on pages 75--78 in [CS], which
supply a proper subset of $\mathcal{MM}$, Carter and Saito then
explain, (on p.74), rules for enlarging this subset to the full
collection $\mathcal{MM}$:
 {\bf ``Furthermore, we include
the following four variations to the list".} The fourth of these
`variations', involves the following instructions:
\smallskip

{\bf ``Change $X$ to $\overline{X}$ and vice versa in the relations consistently
 whenever possible.....}\footnote{Note: In the notation used in the present paper,
 $X ,\overline{X}$ are replaced by $NW,NE$ respectively.}
{\bf Thus a given sentence may also be valid with such a replacement, and there is a
move on sentences when these (and similar) replacements are valid. Add such
variations to the list...}"
\smallskip

As the present authors interpret these instructions for
constructing `variations' on the original smaller list of
movie-moves, they involve an easier combinatorial part,
(recognizing which changes are `consistent'), and
a much harder topological part (recognizing which variations are `valid'):\\
As applied to a movie-move $M==M'$, the easier part involves changing some
 (not necessarily all) NE's to NW's
and vice-versa, in all possible ways such that the resulting pair $(\overline M,\overline M')$
is again acceptable ---``consistent"--- in some
sense not explicitly stated in [CRS]. Results seem
compatible with those in [BL] if we interpret ``consistent'' as
meaning, that $\overline M$ and $\overline M'$ are again movies, such that
$$source(\overline M)=source(\overline M')\mbox{ and }target(\overline M)=
target(\overline M')\;.$$ The harder topological part then seems
to require recognizing, for which of these newly constructed
consistent pairs, $\overline M$ and $\overline M'$ are isotopic
(which seems to be the meaning here of `valid')
---while not undecidable, this still involves a hefty sequence of
rather nontrivial topological exercises
 in order to construct $\mathcal{MM}$ following these instructions. Let us
note also that the set of movie-moves $(M,M')$ to which these
instructions
must be applied, is infinite.\\
In other words, this construction of the set $\M{MM}$ seems not to
be purely \emph{combinatorial}---while this is perhaps not so
relevant to the original purposes of [CRS], it is an important
obstacle to the purposes of the present paper.

 These auxiliary topological
exercises are not required (or have already been carried out)
 in the purely combinatorial construction of $\mathcal{MM}$ presented by Baez and Langford ( in
[BL] and [LL]), which the present authors found
 invaluable for the insertion of the relevant
movie-moves into the program for computing amplitudes made
available in footnote 2, as well as for the purely mathematical
considerations in \S\ref{SS:balanced}. (We are indebted to Carter
for suggesting to us, that the work of Baez and Langford would be
useful for the requirements of our program in this connection.)
\medskip

\noindent \underline{CAUTION:} Prof. Langford has informed us
about the following erratum in the construction in [BL] and [LL]
just discussed: namely, there is a typo in the movie-move numbered
14 in [BL] and 21 in [LL],which can be corrected as follows:
replace the fourth flicker in the left-hand movie by
$$B_{i;m,n}J_{A;m+1,n}X_{m,n+2}$$
both in [LL] and in [BL]; and, in [LL], replace the third flicker
in the right-hand movie by
$$\cap_{m,n+2}W_{i;m+1,n}\;.$$
\medskip

Here are some useful definitions, suggested by the preceding discussion:
\begin{Def}\label{Def:valid}
Let $M,M'$ be two movies.\\
a)The ordered pair $(M,M')$ will be called \underline{{\bf grammatical}}
provided $M$ and $M'$ have the same source and the same target.\\
b)The ordered pair $(M,M')$ will be called \underline{{\bf valid}}
provided that $M$ and $M'$ are isotopic.
\end{Def}
\medskip

\noindent {\bf REMARKS:} `valid' implies `grammatical', but not
conversely. Note that the assertion: `$(M,M')$ is grammatical' is
purely combinatorial, and its truth-value is easily decided, while
the assertion `$(M,M')$ is valid' seems to lie much deeper. The
present authors do not know if `validity' (in this sense) is
effectively decidable. All movie-moves are valid, hence are
grammatical.
\begin{Prop+Def}\label{Prop:in}
If the ordered pair of movies $\mathcal{M}=(M,M')$ is grammatical (hence \emph {a fortiori}
if it is valid) then
$$in(M)=in(M'), out(M)=out(M')$$
(in which case let us define $in(\mathcal{M})$ to be the common value of $in(M),in(N)$,
and similarly $out(\mathcal{M})$ to be the common value of $out(M),out(N)$)
\end{Prop+Def}
\noindent {\bf PROOF:} Immediate.
\medskip

Let us now begin to examine more carefully the collection of
movie-moves constructed in [BL], pp.47-50. The construction in
[LL] is essentially the same, but here the formulation in [BL]
will be preferred, because it has the desirable feature of
utilizing the same numbering for the movie-moves, that is used in
[CS] and in [CRS]. {\bf In the following, when we say a movie-move
is of type $i$, it is always this numbering to which reference is
made.} In [BL], the notation they use is that explained at the end
of \S \ref{SS:ET} in the present paper; in our discussion of the
work of Baez and Langford, we shall translate their notation into
the ``$sf$-notation'' explained above in \S\ref{SS:notation}.

\begin{tabular}{c}
\hline
AT THIS POINT, AN IMPORTANT COMPLICATION ARISES:       \\
\end{tabular}

To be quite precise concerning the set \MMM, we must note that
while Baez-Langford ([BL]) have 30 basic types of movie-moves,
Carter-Rieger-Saito ([CRS]) have a 31st type of movie-move, (the
last in the list on pp.75--~79 in [Carter-Saito]), which is
omitted from the list in Baez-Langford. On p.47 of [BL], Baez and
Langford say concerning this, ``We omit their 31st movie-move,
since it follows from the definition of a 2-category.". Thus, two
possible candidates for $\MMM$ seem to present themselves: the one
consisting of 31 types, as defined by [CRS] with the difficulty
explained above, the other as defined more explicitly in [BL] and
[LL] and consisting only of the first 30 types in [CRS].
\medskip

 For our present paper, the set $\MMM$ must be
defined so that Th.\ref{th:main} holds. We shall prove below that
Th.\ref{th:main} would NOT hold (at least, as here stated), if we
took $\MMM$ to consist only of movie-moves of types 1--30. Thus,
for the purposes of the present paper, we must consider (as do
Carter-Rieger-Saito in their original proof of Th.\ref{th:main})
$\mathcal{MM}$ as also including all these movie-moves of type 31.
We shall thus adopt the notation
\begin{equation}\label{eqn:MMdecomp}
\mathcal{MM}=\mathcal{MM}_{BL}\cup \mathcal{MM}_{31}
\end{equation}
 where
$\mathcal{MM}_{BL}$ denotes the set of movie-moves of types 1--30,
as explicitly defined in [BL]--- while $\mathcal{MM}_{31}$ denotes
the set of movie-moves of type 31, as explicitly defined  in \S
\ref{SS:31} below, and further discussed in \S \ref{SS:31a}. More
generally, we denote the set of all movie-moves of type $i$ by
$\MMM_i$, for $1\leq i\leq 31$.

 We must note explicitly that, as
 a consequence of [BL] and [LL], Th.\ref{th:main} indeed holds with
the definition of $\MMM$  just indicated.
\subsection{Movie-move 31}\label{SS:31}
We now turn to the explicit definition of movie-moves of type 31,
 in terms of the $sf$-notation explained in \S \ref{SS:notation}. First, here
is how the movie-moves of type 31 are defined on p.76 of [CS]:
\smallskip

``$(WZ_1Z_2V,WZ_1'Z_2V,WZ_1'Z_2'V)\leftrightarrow
(WZ_1Z_2V,WZ_1Z_2'V,WZ_1'Z_2'V)$ where the changes $Z_i$ to $Z_i'$
 for i=1,2 are FESI's."
 \medskip

We next translate this into our $sf$-notation. (It must be
admitted this looks a bit clumsier in our present more explicit
$sf$-notation, which however seems better-suited for some later
computations in \S 3.3 below.) \footnote{By ``FESI" [CS] mean,
what we are here calling ``elementary transition". To obtain
Def.\ref{MM31} from the [CS] definition, we replace:
$$F_i:Z_i\mbox{ to }Z_i' \mbox{ by }E_i:U_i=>V_i, W\mbox{ by }W,V\mbox{ by }X\;.$$ }

\begin{Def}\label{MM31}
Let
$$E_1:U_1=>V_1\mbox{ and }E_2:U_2=>V_2$$
be two elementary transitions (not necessarily different); and let
$$m_1,n_1,m_2,n_2$$
be natural numbers such that
\begin{equation}\label{compos31}
m_1+out(E_1)+n_1=m_2+in(E_2)+n_2\;;
\end{equation}
 then we define (cf. Fig. \ref{figMovieMove31})
 \begin{equation}\label{eqn:MM31}
\mathcal{M}(E_1,E_2;m_1,n_1,m_2,n_2)=(M_1,M_2)
\end{equation}

\begin{figure}
  \centering
  \includegraphics[scale=0.5]{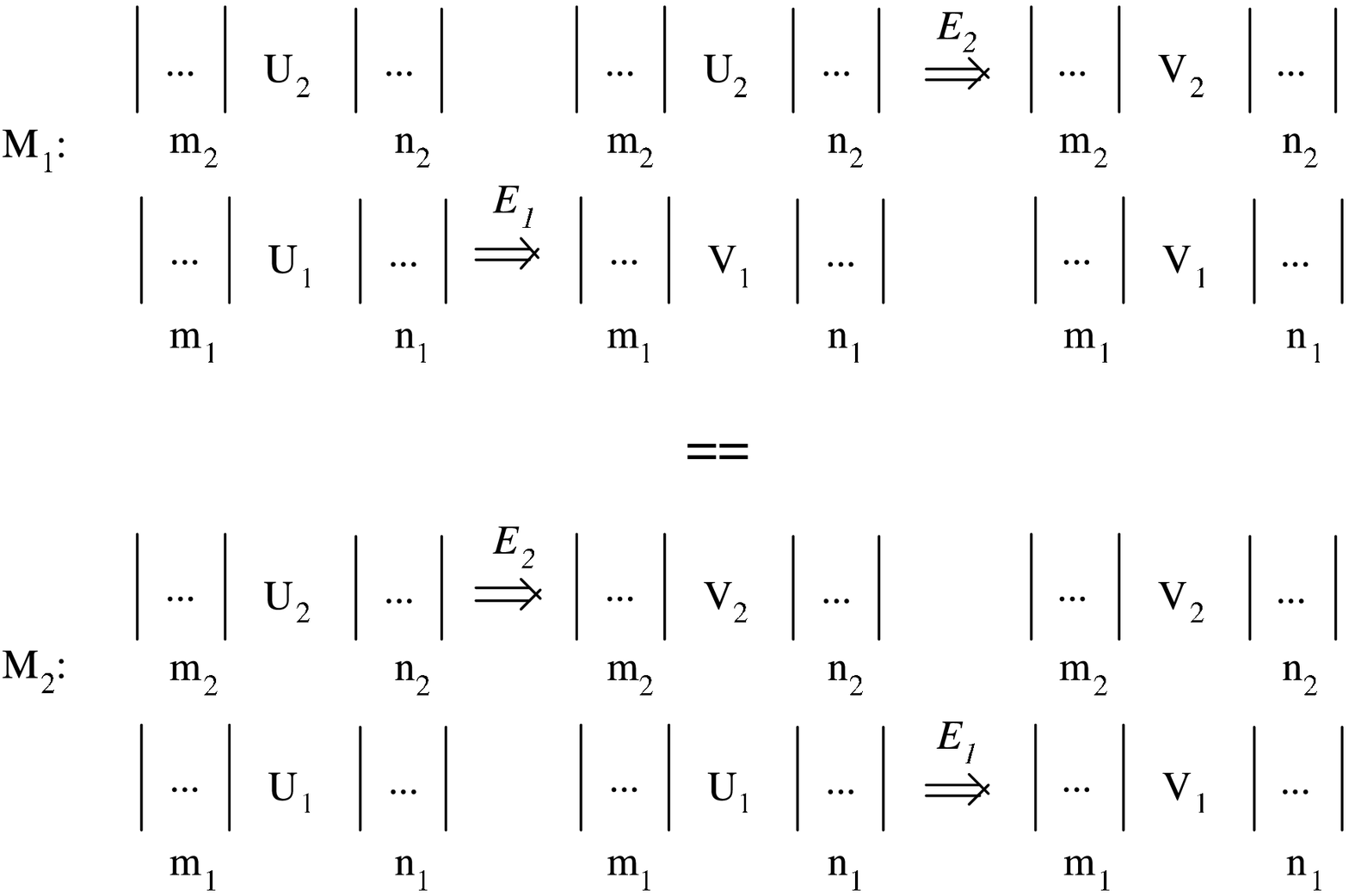}
  \caption{Movie Move 31}\label{figMovieMove31}
\end{figure}

to be the ordered pair of movies given (in $sf$-notation), by:
\smallskip

\noindent $M_1= [$\\
$(m_2\bullet U_2\bullet n_2)\circ s(m_1\bullet U_1\bullet n_1)s$\\
$=> t(E_1)$\\
$s(m_2\bullet U_2\bullet n_2)s\circ f(m_1\bullet V_1\bullet n_1)f$\\
$=> t(E_2)$\\
$f(m_2\bullet V_2\bullet n_2)f\circ (m_1\bullet V_1\bullet n_1)$\\
 $]$

 and by\\
$M_2= [$\\
$s(m_2\bullet U_2\bullet n_2)s\circ (m_1\bullet U_1\bullet n_1)$\\
$=> t(E_2)$\\
$f(m_2\bullet V_2\bullet n_2)f\circ s(m_1\bullet U_1\bullet n_1)s$\\
$=> t(E_1)$\\
$(m_2\bullet V_2\bullet n_2)\circ f(m_1\bullet V_1\bullet n_1)f$\\
$]$\\
(Recall that if $E$ is an elementary transition of type $i$, we
have set $t(E)=i$.)
 We then define a \underline{{\bf movie-move of type $31$}} to
be an ordered pair of movies
\begin{equation}\label{eqn:MM31+}
(W\circ (M_1)\circ X, W\circ (M_2)\circ X)
\end{equation}
where $W,X$ are stills and either $(M_1,M_2)$ or $ (M_2,M_1)$ is
an ordered pair of the form (\ref{eqn:MM31}) with each of which
$W,X$ are composable. We denote the set of all movie-moves of type
31 by $\mathcal{MM}_{31}$.
\end{Def}

\subsection {Some of the Structure on the Set $\mathcal{MM}$ of
Movie-Moves}\label{SS:MmStruct} In the first place, we have
$$\MMM=\bigcup_{1\leq i\leq 31}\MMM_i=(\MMM_{BL})\bigcup \MMM_{31}$$
where, for $1\leq i\leq 31$, $\MMM_i$ is the set (infinite for
each such $i$) of all movie-moves of type $i$, and where
$$\MMM_{BL}=\bigcup_{1\leq i\leq 30}\MMM_i\;.$$

Let us next note that every movie-move is valid, and hence
grammatical (in the sense of Def. \ref{Def:valid}).

We next describe some structure on the set $\MMM^+$ of all
grammatical ordered pairs of movies---i.e., the set
 $$\{(M,M') | M\mbox{ and }M'\mbox{ are movies, }
 source(M)=source(M'),target(M)=target(M')\}$$
(It will then be proved that this structure is  inherited via
 $$\MMM_{BL} \subseteq \MMM^+$$
 by the set $\MMM_{BL}$, at least in the sense given by
 Prop.\ref{MMclosed} below.)
\begin{Def}\label{def:compMM}
Let $(M,M')$ be a grammatical ordered pair of movies, and let $B,C$ be two stills.
Then we shall say $B,C$ are \underline{{\bf composable}} with $(M,M')$, provided that
\begin{equation}\label{eqn:stillMovie}
target(C)=in(M)(=in(M')),\;source(B)=out(M)(=out(M'))
\end{equation}
(cf.Prop.\ref{Prop:in}) in which case we define $B\circ (M,M')\circ C$ to  be
the  pair of movies given by
\begin{equation}\label{eqn:stillMMove}
B\circ (M,M')\circ C :=(B\circ M\circ C ,B\circ M'\circ C)
\end{equation}
\end{Def}
\smallskip

\begin{Def}\label{def:circPair}
Let $(M,M')$ be a grammatical pair of movies, and let $\pi\in \mathcal{SYM}$;
then we define $\pi\circ (M,M')$ to be the pair $(\pi\circ M,\pi\circ M')$
\end{Def}
\begin{Def}\label{def:bulletPair}
Let $(M,M')$ be a grammatical pair of movies, and let $m,n$ be
natural numbers; then we define $m\bullet (M,M')\bullet n$ to be
the pair
$$(m\bullet M\bullet n, m\bullet M'\bullet n)$$
\end{Def}
\begin{Def}\label{def:converse}
Let $(M,M')$ be a grammatical pair of movies; then we define $(M,M')^C$
--- the \underline{{\bf converse}} of $(M,M')$ --- to be $(M',M)$
\end{Def}
We omit the easy proof of the following straightforward but useful proposition:
\begin{Prop}\label{propdef:MM}Let $M,M'$ be a grammatical pair of movies.\\
a) Suppose the stills B,C are composable with $(M,M)$; then also
$$B\circ (M,M')\circ C$$ is
grammatical. If also (M,M') is valid, so is $B\circ (M,M')\circ C$.\\
b)If $\pi\in\mathcal{SYM}$ then $\pi\circ (M,M')$ is grammatical.
If also $(M,M')$ is valid,
so is $\pi\circ (M,M')$.\\
c)Let $m,n$ be natural numbers; then
$$m\bullet (M,M')\bullet n=(m\bullet M\bullet n,m\bullet M'\bullet n)$$
is grammatical, and is valid if ($M,M'$) is.\\
d) $(M,M')^C$ is grammatical, and is valid if $(M,M')$ is.
\end{Prop}
\noindent {\bf REMARK:} In a sense, this proposition underlies the
construction of the list $\mathcal{MM}$ of movie-moves, both in
[BL] and in [CRS]. The sense in which this is so will become
clearer in a little while.
\medskip

\begin{figure}
\centering
\includegraphics[scale=0.5]{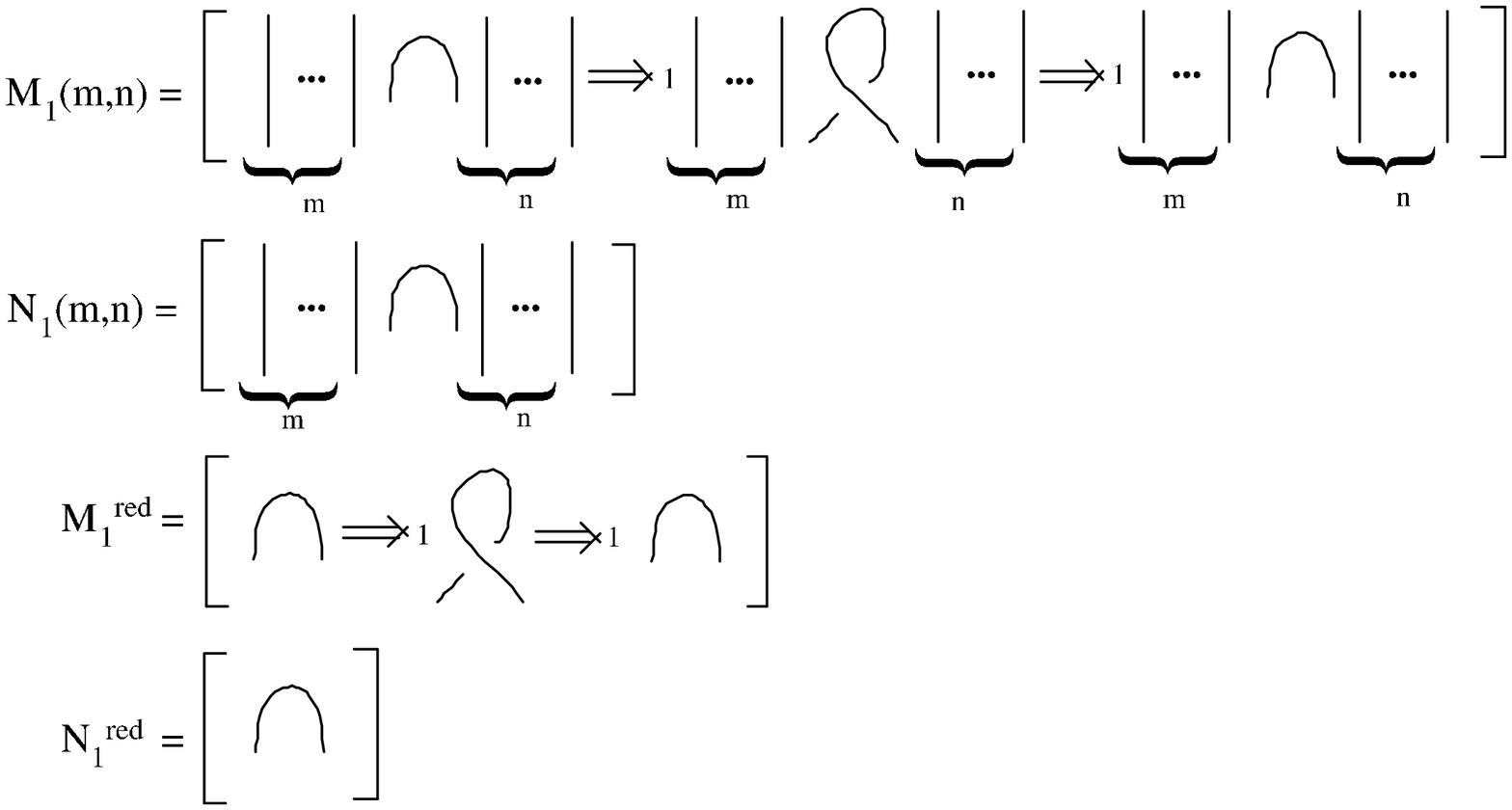}
\caption{\label{fig:mm1} Movie Move I. $M_1(m,n) == N_1(m,n)$ and
$M^{red}_1 == N^{red}_1$ }
\end{figure}

We have proceeded about as far as we can in this subsection,
without actually examining the detailed explicit construction of
$\MMM_{BL}$ in [BL]. To begin this examination, let us consider
the case $\MMM _1$:\\
$\bullet$ {\bf EXAMPLE A):} Thus, we begin by examining (in the
light of the preceding definitions) the first movie-move in [BL].
(This
 is illustrated by Fig.52 on p.86 of [CS], where
this movie-move has the bemusing name, ``an elliptic confluence of
branch points". See also Fig. \ref{fig:mm1} in the present paper.)
For movie-move 1, [BL] (p.47) lists $W_{m,n}W^*_{m,n}=1$. In the
notation of the present paper, this is the movie-move
$(M_1(m,n)==N_1(m,n))$, where $
M_1(m,n)=[ \\
s[m,Cap,n]s\\
=>1\\
 sf[m,NW,n][m,Cap,n]sf\\
 =>1\\
  f[m,Cap,n]f\\
]$\\
 and
$N_(m,n)=[[m,Cap,n]].$\\
$m$ and $n$ represent arbitrary natural numbers, so we have here
an infinity of movie-moves of type 1 (though not {\bf all}
movie-moves of type 1 are in this list.)
\medskip

\underline{{\bf All}} movie-moves of type 1, can be derived (by
processes next to be explained) from the single `reduced'
movie-move
$$(M^{red}_1==N^{red}_1)=(M(0,0),N(0,0))$$
obtained by setting $m=n=0$ in the movie-move above.
\medskip

\begin{figure}
\centering
\includegraphics[scale=0.5]{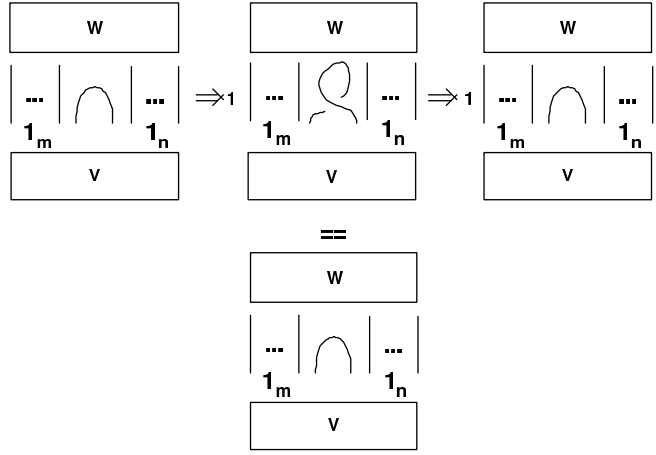}
\caption{\label{fig:gen_mm1} General Movie-Move of Type 1
(Carter-Rieger-Saito)}
\end{figure}

To see how this works, note first that the most general movie-move
of type 1, in the Carter-Rieger-Saito notation, reads
\begin{equation}\label{eqn:MM1}
(W\cap _{m,n}V,W\cap _{m,n}X_{m,n}V,W\cap _{m,n}V)\leftrightarrow
(W\cap _{m,n}V)
\end{equation}
(See also Fig. \ref{fig:gen_mm1} in the present paper.) Here $V$
and $W$ denote arbitrary stills, subject to the composability
conditions
$$source(W)=m+n,\;target(V)=m+n+2$$

In the notation explained in \S \ref{SS:notation}, combined with
the above Def.\ref{def:compMM} (and switching from the [CRS]
convention which reads stills from top to bottom, to the opposite
convention utilized in the present paper) this movie-move (read
from left to right) may be written:
$$V\circ (m\bullet M^{red}\bullet n)\circ W ==
V\circ (m\bullet N^{red}\bullet n )\circ W \eqno(A)$$
where we define $[M^{red}==N^{red}]$, as above,  to be the movie-move defined by:\\
$M^{red}=[s[0,Cap,0]s\;=>1\;\;  sf[0,Cap,0][0,NW,0]sf\;=>1\;\;  f[0,Cap,0]f ],$\\
$N^{red}=[[0,Cap,0]].$\\
---i.e. the infinitely many movie-moves given by (\ref{eqn:MM1}), are all derived via $(A)$
from the {\bf single} 'reduced' movie-move
$$\M{M}^{red}_1=(M^{red}==N^{red})\;.$$ Finally, the most general movie-move of
type 1, is obtained via this framing-process from one of the 8
`reduced' movie-moves
$$\pi\circ (M^{red},N^{red})\;\;\;\;\;(\pi\in\mathcal{SYM})$$
i.e. has the form
$$V\circ (m\bullet (source(\pi\circ \M{M}^{red}_1)\bullet n)\circ W ==
V\circ (m\bullet(target(\pi\circ \M{M}^{red}_1)\bullet n)\circ
W\;\; \eqno(B)$$
\bigskip

Precisely the same pattern holds throughout the entire
construction of $\mathcal{MM}_{BL}$ in [BL], that was just
exhibited for movie-move 1. To explain this pattern in more
detail, two further definitions will be helpful:
\begin{Def}\label{def:strictDerived}
Given a grammatical pair of movies $(M,N)$, a pair
$(\overline{M},\overline{N})$ will be called \underline{{\bf
strictly derived from}} $(M,N)$ when either
$(\overline{M},\overline{N})$ or its converse  has the form $(A)$,
i.e., when there exist natural numbers $m$ and $n$, and stills $V$
and $W$, such that
$$\mbox{either }(\overline{M},\overline{N})\mbox{ or }(\overline{N},\overline{M})=
(V\circ (m\bullet M\bullet n)\circ W,\;V\circ (m\bullet N\bullet
n)\circ W ) \;.$$
\end{Def}
(Of course, this equation presupposes that $V$ and $W$ satisfy the
composability conditions of Def. \ref{def:compMM})
\begin{Def}\label{def:extendedDerived}
Given a grammatical pair of movies
$$\M{M}=(M,N)\;,$$ a pair
$(\overline{M},\overline{N})$ will be called \underline{{\bf
derived from}} $(M,N)$ \underline{{\bf in the extended sense,}}
when there exists $\pi\in\mathcal{SYM}$ such that
$(\overline{M},\overline{N})$ is strictly derived from $\pi \circ
\M{M}$.
\end{Def}
\medskip

{\bf The above discussion  of movie-moves of type 1 asserts, in
terms of the two preceding definitions, that all movie-moves of
type 1 are derived from a \emph{single} ``reduced'' one
$(M^{red}(1),N^{red}(1))$ by the derivation process (in the
extended sense).}

Having found this formulation. we next observe that it essentially
holds more generally ---with a suitable definition of `reduced'
movie-moves--- for the entire list of movie-moves in [BL],p.47-50
 (all of which have type $<31$). To see this, it must first be noted that, for $1\leq
i\leq 30$, the movie-moves of type i listed in ([BL],p.47---50)
consist of one or more movie-moves involving two arbitrary natural
numbers $m$ and $n$: i.e., the movie-moves there listed, are the
elements of a set (possibly infinite, depending on $i$), this set
being explicitly described in [BL], which has the form
$$\{(M_{\alpha}^{m,n}(i)==N_{\alpha}^{m,n}(i))| \alpha\in \mathcal{P}_i \}$$
We then define (still for $1\leq i\leq 30$ ) the \underline
{BL-{\bf reduced movie-moves of type} $i$} to be the elements of
the set\footnote{We shall also define later the set
$\mathcal{MM}^{red}(31)$ in a somewhat more complicated way.}
\begin{equation}\label{eqn:BLreduced}
\mathcal{MM}^{red}(i):=\{(M_{\alpha}^{0,0}(i)
==N_{\alpha}^{0,0}(i))| \alpha\in \mathcal{P}_i \}
\end{equation}
obtained by setting $m=n=0$ in the former set.\\

Example A illustrates this construction. Here are two more
examples:\\
$\bullet$ {\bf EXAMPLE B):} \underline{i=5}\\
For movie-moves of type 5, ([BL]) lists 6 families involving two
natural numbers $m$ and $n$. These families are indexed by the set
$\mathcal{P}_5$ consisting of the 6 elementary transitions which
in the [BL] notation are\\
$S_0$,$S_1$,$S_2$,$S_0^*$,$S_1^*$,$S_2^*$; and, in the notation
of the present paper, are:\\
ET3tmb,ET3tbm,ET3mbt,ET3tmbR,ET3tbmR,ET3mbtR. \\
If $\M{E} =[S=>T]$ is one of these 6 elementary transitions in
$\mathcal{P}_5$, the corresponding family of movie-moves of type
5, is
$$((M_{\M{E}}^{m,n}(5)==N_{\M{E}}^{m,n}(5))$$
where
$$ M_{\mathcal{E}}^{m,n}(5)=[[m,S,n]=>3\;\; [m,T,n]=>3\;\;
[m,S,n]]$$ and
$$N_{\mathcal{E}}^{m,n}(5)=[[m,S,n]]$$
---while the 6 BL-reduced movie-moves in
$\mathcal{MM}^{red}(5)$ are obtained by setting $m=n=0$
in the preceding.\\
$\bullet$ {\bf EXAMPLE C):} \underline{i=17}\\
$\mathcal{P}_{17}$ is infinite (as will become become clear below
in the proof of Prop. \ref{prop:semi17to20 })--- and so there are
infinitely many BL-reduced movie-moves of type 17.
\medskip

To sum up, the list in [BL],p.47-50, gives rise (by the procedure
just described) to a collection (unfortunately still infinite) of
``BL-reduced" movie-moves, from which every other movie-move may
be obtained by the process explained in Def.
\ref{def:extendedDerived}
---more precisely:
\begin{Prop}\label{prop:red}
Let $1\leq i\leq 30$; then the collection
\begin{equation}\label{red}
\mathcal{MM}^{red}(i):=\{(M_{\alpha}^{red}(i)==N_{\alpha}^{red}(i))|
\alpha\in \mathcal{P}_i \}
\end{equation}
of movie-moves, constructed above, has the property that an
ordered pair $\M{M}$ of movies, is a  movie-move of type $i$, if
and only if  $\M{M}$ is derived in the extended sense (i.e, in the
sense of the preceding Def.\ref{def:extendedDerived}) from one of
the movie-moves in the list (\ref{red}).
\end{Prop}
\smallskip
As an immediate Corollary, we have (as promised above):
\begin{Prop}\label{MMclosed}
If $(M,N)$ is any movie-move of type $\leq 30$, every pair derived
from it in the extended sense is again a movie-move.
\end{Prop}

 We define the set of all
BL-reduced movie-moves by
\begin{equation}\label{MMred}
\mathcal{MM}^{red}_{BL}:=\bigcup_{1\leq i\leq 30}
\mathcal{MM}^{red}_i
\end{equation}
\medskip

\noindent \underline{{\bf NOTE:}} Let $1\leq i \leq 30$. The set
of movie-moves of type $i$, as defined in [BL], is a {\bf proper}
subset of the set denoted in this paper by $\MMM_i$---the latter
set consists of all movie-moves derived from the former set via
Def.\ref{def:compMM}. (This more generous definition of $\MMM$ and
$\MMM_i$ is needed, in this paper, so that Th.\ref{th:main} may
hold.) The relation between these two sets is very simple, and may
become clearer, upon comparing Figures \ref{fig:mm1} and
\ref{fig:gen_mm1}.

\section{$\mathcal{U}$-Regular Amplitude-Invariants for Movies,
\mbox{2-Knots} and 2-Tangles}\label{AMPS} Throughout this section,
it will be convenient to replace the ground-ring $A=\ZZ
[q,q^{-1}]$ by its quotient-field $$\KK :=\QQ (q)$$ (mainly
because it is easier to handle systems of linear equations over
$\KK$ than over $A$.) We denote by $V_1$, the two-dimensional
vector-space $V\otimes_A \KK$ over $\KK$. If $S$ is a still, we
shall (by a slight abuse of notation) denote the $\KK$-linear
transformation
$${<}S{>}\otimes_A \KK : V_1^{source(S)}\rightarrow V_1^{target(S)}$$
simply by ${<}S{>}$, and will continue to refer to it as the
Kauffman bracket of $S$.
\subsection{$\mathcal{A}$-Amplitudes for Elementary Transitions, Flickers
and Movies}\label{S:amp}
\begin{Def}\label{Def:ampassign}
By an \underline{{\bf amplitude-assignment}} will be meant a map
$\mathcal {A}$ which assigns to every elementary transition $E\in
\mathcal{ET}$, a $\KK$-linear map
$$<E>_{\mathcal{A}}:V_1^{\otimes in(E)}\longrightarrow V_1^{\otimes out(E)}\;.$$
which will be called the $\mathcal{A}$-{\bf amplitude} of $E$.
\end{Def}

Such an assignment of `amplitudes' to elementary transitions,
extends immediately to assignments of `amplitudes' to flickers and
to movies, as follows:
\begin{Def}\label{Def:ampFlicker}
Let $\mathcal{A}$ be an amplitude-assignment, and let $F$ be a
flicker. Then the $\mathcal{A}$-{\bf amplitude} of $F$ is defined
to be the $\KK$-linear map
$$<F>_{\mathcal{A}}:V_1^{\otimes in(F)}\longrightarrow V_1^{\otimes out(F)}$$
constructed as follows:
\medskip

\noindent By Def.\ref{defFlicker}, $F$ has one of the two forms
(2) or (5).

Suppose first (2) holds, so $F=(\mathcal{E},B,m,n,T) $ (as
illustrated by Fig. \ref{fig:flicker}). Then we set
$$<F>_{\mathcal{A}}:=<T>\circ(Id_{V_1^{\otimes m}}\otimes <\mathcal{E}>_{\mathcal{A}}
\otimes Id_{V_1^{\otimes n}})\circ <B> :  V^{\otimes
in(F)}\longrightarrow V_1^{\otimes out(F)}$$ (where $<B>,<T>$
denote the Kauffman brackets of the stills $B,T$.)

\noindent On the other hand, if (5) holds, so $F=1_S:S=>S$, where
$S$ is a still with source $m$ and target $n$, then we define
$$<F>_{\mathcal{A}}:V^m\longrightarrow V^n$$
to be simply the Kauffman bracket $<S>$ (which we note is, in the
case that (5) holds, of course independent of the choice of
$\mathcal{A}$).
\end{Def}
\medskip

\begin{Def}\label{def:ampMovie}
Let $\mathcal{A}$ be an amplitude-assignment, and let
$$M=F_1\cdots F_s$$ be a movie. Let us label the $s+1$ stills of
the movie $M$ as follows:
$$S_1=source(F_1),S_i=source(F_i)=target(F_{i-1})
\mbox{ for $1<i\leq s$, }S_{s+1}=target(F_s)\;.$$ Then we
construct as follows, a $\KK$-linear map
$$<M>_{\mathcal{A}}: V_1^{\otimes in(M)}
\longrightarrow V_1^{\otimes out(M)}\;.$$
---which we shall call the  $\mathcal{A}$-{\bf amplitude of M}---\\
namely, this map is defined to be the alternating sum
\begin{equation}\label{eqn:ampMovie}
<M>_{\mathcal{A}}:= <S_1>-<F_1>_{\mathcal{A}}+<S_2>
-<F_2>_{\mathcal{A}}+\cdots -<F_s>_{\mathcal{A}}+<S_{s+1}>\;.
\end{equation}
\end{Def}

\subsection{$\MU$-balanced and Semi-normal Amplitude-Assignments}\label{SS:balanced}
Define
$$\underline{31}:=\{i\in \ZZ\; |\; 1\leq i\leq 31\}\;,$$
Throughout the following discussion, we shall let $\mathcal{U}$
denote any proper subset of $\underline{31}$; and
$$\mathcal{U}^c:=\{i\in \underline{31}\; |\; i\notin \mathcal{U}\}$$
will denote the non-empty complement of $\mathcal{U}$ in
$\underline{31}$.
\medskip

The reader is asked to review the definition (in the Introduction)
 of $\MU$-regular isotopy of movies (Def.\ref{def:Uregular}). Note
that this equivalence relation on movies, is stronger than (or
equal to) isotopy.

 As a technique for manufacturing $\MU$-regular isotopy invariants of
movies, we shall utilize the $\MA$-amplitude construction of
Def.\ref{def:ampMovie}.
\medskip

Thus, we now turn to the question: \underline
{when does the construction $<M>_{\mathcal{A}}$}\\
\underline{furnished by Def.\ref{def:ampMovie}, yield
$\MU$-regular isotopy-invariants of movies?}
\begin{Def}\label{def:balanced}
Let $\mathcal{A}$ be an amplitude-assignment.
We define $\mathcal{A}$
to be \underline{{\bf $\MU$-balanced}} if it has this property:\\
Whenever $M$ and $M'$ are $\MU$-isotopic movies---whence, by
Prop.\ref{Prop:in},
$$in(M)=in(M'),out(M)=out(M')$$
---we then have
\begin{equation}\label{eqn:ampBal}
<M>_{\mathcal{A}}=<M'>_{\mathcal{A}}:V_1^{\otimes
in(M)}\longrightarrow V_1^{\otimes out(M)}.
\end{equation}
\end{Def}
\begin{Def}\label{def:respect}
Let $\mathcal{A}$ be an amplitude-assignment, and let
$$\mathcal{M}=(M_1,M_2)$$
be an ordered pair of movies.Then we
say that \underline{{\bf $\mathcal{A}$} respects $\mathcal{M}$}
if the following holds:
$$<M_1>_{\mathcal{A}}=<M_2>_{\mathcal{A}}$$
\end{Def}

We note the following immediate consequence of our
definitions:\\
\underline{An amplitude-assignment $\mathcal{A}$ is
$\MU$-balanced, iff. it respects all movie-moves}\\
\underline{ whose type is not in \MU.}\\
As it stands, to see if $\mathcal{A}$ is $\MU$-balanced, this
criterion requires us to test, for every $i\in \MU^c$, all
movie-moves in the infinite set $\mathcal{MM}_i$. However, the
considerations in the remainder of this sub-section, and (to deal
with movie-move 31) in the next, allow us to cut down, (at the
cost of restricting to a smaller collection of {\bf `semi-normal'}
amplitude-assignments, explained in Def. 48 below) to a finite
number of tests, thus making feasible (as will be seen) a computer
approach to these matters.
\begin{Prop}\label{red:comp}
Let $\mathcal{A}$ be an amplitude-assignment, and let
$$\mathcal{M}=(M==N)\in \mathcal{MM}$$
be a movie-move. Suppose $\mathcal{A}$ respects $\mathcal{M}$.\\
a) Let $B$ and $C$ be stills, such that
$$target(C)=in(M), \;source(B)=out(M)$$
(i.e., such that $B,C$ is composable with $\mathcal{M}$). Then
$\mathcal{A}$ respects
$$B\circ \mathcal{M} \circ C:=((B\circ M \circ C)==(B\circ N\circ C))$$.\\
b) Let $m,n$ be natural numbers; then
$\mathcal{A}$ respects
$$m\bullet \mathcal{M}\bullet n:=((m\bullet M\bullet n)==(m\bullet N\bullet n))\;.$$.
\end{Prop}
{\bf PROOF:}By hypothesis,
$$<M>_{\mathcal{A}}=<M'>_{\mathcal{A}} :V_1^{\otimes in(M)}\longrightarrow V_1^{\otimes out(M)}$$
We have
$$<B\circ M\circ C>_{\mathcal{A}}=<B>\circ <M>_{\mathcal{A}}\circ <C>=
<B>\circ <N>_{\mathcal{A}}\circ <C>=<B\circ N\circ
C>_{\mathcal{A}}$$ which proves a). Similarly,
$$ <m\bullet M\bullet n>_{\mathcal{A}}=1_{V_1^{\otimes m}}
\otimes <M>_{\mathcal{A}}\otimes 1_{V_1^{\otimes n}}
=1_{V_1^{\otimes m}}\otimes <N>_{\mathcal{A}}\otimes
1_{V_1^{\otimes n}}=<m\bullet N\bullet n>_{\mathcal{A}}$$ which
proves b).
\\Q.E.D.
\medskip

As an immediate corollary, we have:
\begin{Prop}\label{prop:respDerived}
If an amplitude-assignment $\mathcal{A}$ respects every element in a collection
$\mathcal{C}$ of ordered pairs of movies, then it respects every ordered pair of
movies strictly derived from $\mathcal{C}$ (in the sense of Def. \ref{def:strictDerived}).
\end{Prop}

Next note, that if $E:S=>T$ is an elementary transition of type either
0 or 8, then $S$ and $T$ have  the same Kauffman amplitude. For this reason, the
following definition makes sense:
\begin{Def}\label{def:seminormalAA}
An amplitude-assignment $\mathcal{A}$ will be called \underline{{\bf semi-normal}}
when, for every elementary transition $E:S=>T$ of type either 0 or 8,
$$<E>_{\mathcal{A}}=<S>=<T>\;.$$
\end{Def}
As noted at the end of \S\ref{SS:moviemoves}, we shall use the
same numbering for movie-moves, as that employed in [CS],[CRS] and
(except for the inclusion in the present paper of movie-move 31 in
$\mathcal{MM}$) in [BL].
\begin{Prop}\label{prop:semiRespects }
If $\mathcal{A}$ is a semi-normal amplitude-assignment, then
$\mathcal{A}$ respects all movie-moves of types 15,16 or 22.
\end{Prop}
{\bf PROOF:} Let $(M==M')$ be a movie-move of type 15, 16 or 22.
We begin by showing that every
flicker in $M$ and in $M'$ is of type 8:\\
Indeed, it is easy to see that if this is so for $(M==M')$, and if
$$\pi\in\mathcal{SYM}$$ then also the movie-move $\pi\circ (M==M')$
 has the same property. Now, all the movie-moves of types
15,16 or 22, result by acting by $\mathcal{SYM}$ on the three
(sets of) movie-moves listed under these numbers on p.48 and 49 of
[BL], and inspection shows each of the latter has the required
property.
\smallskip

Thus every flicker F in $M$ and in $M'$ is of type 8, and since $\mathcal{A}$
is semi-normal, we then have
$$<source(F)>=<target(F)>=<F>_{\mathcal{A}}\;.$$
 Let
$$M=F_1\cdots F_s,\;\;M'=F'_1\cdots F'_t\;;$$
then all of the stills of $M$ (i.e. all the sources and targets of the
flickers $F_1,\cdots,F_s$) have the same Kauffman amplitude, say $A$.
 Also, since $(M,M')$ is a movie-move, it is grammatical,
so $F_1=F'_1$. Hence all the stills of $M'$ also have the Kauffman amplitude $A$, and
so
$$<M>_{\mathcal{A}}=A-A+\cdots +A=A=<M'>_{\mathcal{A}}\;,$$
Q.E.D.

\begin{Prop}\label{prop:semi17to20 }
If $\mathcal{A}$ is a semi-normal amplitude-assignment, then $\mathcal{A}$ respects
all movie-moves of types 17, 18, 19 and 20.
\end{Prop}
{\bf PROOF:}Throughout the following argument, we assume $i\in\{17,18,19,20\}$, and
then set
$$I(i)=
\left\{
\begin{array}{ll}
3 & \mbox{ when $i=17$,} \\
5& \mbox{ when $i=18$,} \\
1& \mbox{ when $i=19$,} \\
4& \mbox{ when $i=20$.}
\end{array}
\right.
$$
We  set $\nu(i)$ equal to the number of elementary transitions of type $I(i)$,
as given by the following table:
\begin{center}
\begin{tabular}{|l|c|c|c|c|c|c|c|c|c|}
i & 17 & 18 & 19 & 20\\
\hline
$\nu (i)$ & 12 & 4 & 8 & 8
\end{tabular}
\end{center}
Finally, we assume the $\nu (i)$ elementary transitions of type $I(i)$, enumerated
and with notation as in \S \ref{SS:ET}, as follows:
$$ET(I(i),j)=\{[F_{ij1}\cdots F_{ijs(j)}  \;\;\;=>I(i) \; \;\;\;G_{ij1}\cdots
G_{ijt(j)}]| 1\leq j\leq \nu(i)\}$$
\smallskip

\begin{figure}
\centering
\includegraphics[scale=0.5]{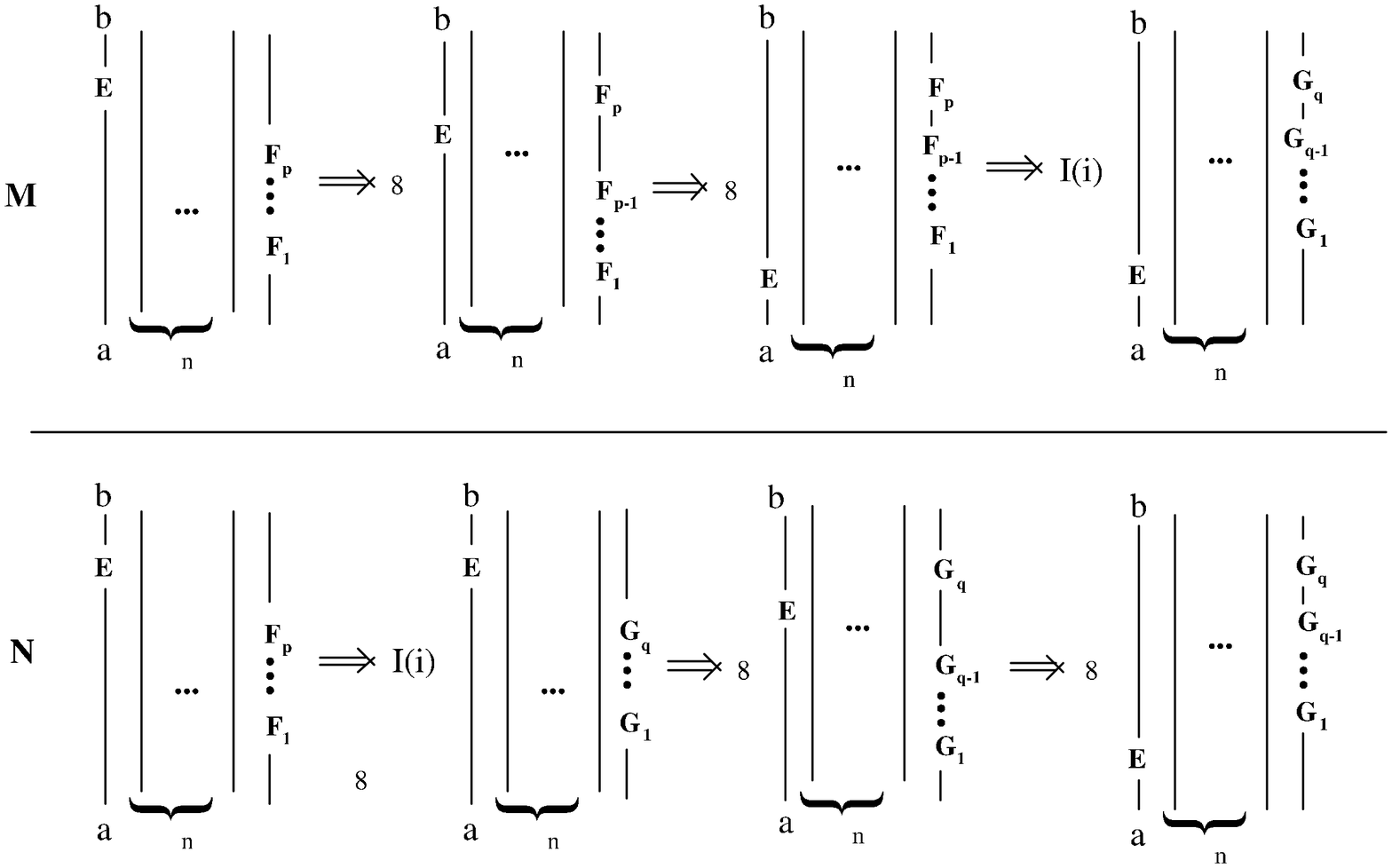}
\caption{\label{fig:L_i_movies} $\mathcal{L}_{i}$}
\end{figure}

The reader is invited at this point to
examine the listings given for movie-move 17--20 on p.49 of [BL];
such an examination
shows that every movie-move of type i, is derived in the extended sense
(as explained above in Def. \ref{def:extendedDerived}) from a movie-move in the two
lists $\mathcal{L}_{i},\;\mathcal{L}'_{i}$ next to be constructed,
whose union
$$ \mathcal{L}_{i}\bigcup \mathcal{L}'_{i}=\mathcal{MM}^{red}_i$$
makes up the set of BL-reduced movie-moves of type i.\\
\underline{{\bf STEP ONE:} $\mathcal{A}$ {\bf RESPECTS ALL MOVIE-MOVES IN}
 $\mathcal{L}_{i}$:}\\
Let us first construct $\mathcal{L}_{i}$:\\
Each movie-move $(M==N)$ in the list $\mathcal{L}_{i}$ is constructed from the following
data:
\begin{Def}\label{def:datum}
An \underline{i-{\bf datum }}$\mathcal{D}_i(E,n,\mathcal{T})$ is
defined to consist of the following three items:
\begin{itemize}
\item An elementary event $E\in \mathcal{E}=\{Cup,Cap,NE,NW\}$, with source $a$ and target $b$:
$$source(E)=a, target(E)=b$$
\item A natural number $n$
\item An elementary transition $\mathcal{T}=[F_1\cdots F_p => G_1\cdots G_q]$ of type $I(i)$
\end{itemize}
\end{Def}
We denote by $\mathcal{M}\mathcal{D}_i(E,n,\mathcal{T})$ the movie-move
$(M,N)$ of type i, constructed
from this data in the manner next to be described. This movie-move is
pictured in Figure \ref{fig:L_i_movies}.
(Then $\mathcal{L}_{i}$ is defined to be the set of all
such $\mathcal{M}\mathcal{D}_i(E,n,\mathcal{T})$)
Note that the common source $S$ of $M$ and of $N$, consists of $E$ above and to the left of
the source $F_1\cdots F_s$ of $\mathcal{T}$, with $n$ vertical strings in between.
 Similarly, the common target $T$ of $M$ and $N$,
consists of $E$ below and to the left of the target $G_1\cdots G_t $ of $\mathcal{T}$, again
with $n$ vertical strings in between.

$M$ and $M'$ then furnish two paths $S=>T$, constructed as follows:\\
In the first place, in accordance with Fig. \ref{fig:L_i_movies},
the left-hand movie $M$ consists of $p+1$ flickers of type 8, followed by one
flicker of type $I(i)$, and is given in $sf$-notation by:
\begin{verbatim}M=[
[a+n,F1,0]...[a+n,F(p-1),0]s[a+n,Fp,0][0,E,n+target(Fp)]s =>8
[a+n,F1,0]...s[a+n,F(p-1),0]f[0,E,n+target(F(p-1))]s[b+n,Fp,0]f => 8
......................................................
......................................................... =>8
f[0,E,n+source(F1)]s[b+n,F1,0]f...[b+n,F_p,0]s =>I(i)
[0,E,n+source(G1)]f[b+n,G1,0]...[b+n,Gq,0]f
]
\end{verbatim}

Similarly, the right-hand movie $N$ consists of one flicker of
type $I(i)$, followed by $q+1$ flickers of type 8, and is, in $sf$-notation:
\begin{verbatim}N=[
s[a+n,F1,0]...[a+n,F(p-1),0][a+n,Fp,0]s[0,E,n+target(Fp)] =>I(i)
f[a+n,G1,0]...[a+n,G(q-1),0]s[a+n,Gq,0]f[0,E,n+target(Gq)]s =>8
........................................................ =>8
........................................................ =>8
f[0,E,n+source(G1)][b+n,G1,0]f...[b+n,Gq,0]
]
\end{verbatim}

\medskip

Having thus constructed $\mathcal{L}_{i}$, let us next verify that the semi-normal
$\mathcal{A}$ respects the preceding movie-move $(M==N)$
in $\mathcal{L}_{i}$. We begin by computing $<M>_{\mathcal{A}}$.

Here we must distinguish two cases:\\
\underline{Case 1:} i=17, 18 or 20, \underline{Case 2:} i=19\\
\underline{Ad Case 1:}
Here I(i)=3, 4 or 5; since the Kauffman amplitude is unchanged by elementary transitions of
these three types, and also by type 8,
it follows that all 5 stills of $M$ have the same Kauffman amplitude,
say
$$A:V^{\otimes (source(E)+n+in(\mathcal{T}))}\rightarrow V^{\otimes
 (target(E)+n+out((\mathcal{T}))}$$
Since $\mathcal{A}$ is, by hypothesis, semi-normal, $A$ is then also the
$\mathcal{A}$-amplitude of the $p+1$ flickers in $M$ of type 8 (i.e., the first $p+1$ flickers
of $M$). As for the last flicker $\mathcal{F}$ in $M$, of type $I(i)$,
its $\mathcal{A}$-amplitude is (by Def. \ref{Def:ampFlicker})
\begin{eqnarray}
<\mathcal{F}>_{\mathcal{A}}&=&(Id_{V^{\otimes (source(E)+n)}} \otimes
<\mathcal{T}>_{\mathcal{A}})\circ (<E>\otimes Id_{V^{\otimes (n+target(Fp)}}) \nonumber \\
&=&<E>\otimes V^{\otimes n} <\mathcal{T}>_{\mathcal{A}} \nonumber
\end{eqnarray}
Hence, by Def. \ref{def:ampMovie},
\begin{eqnarray}
<M>_{\mathcal{A}}&=&A-A+\cdots +A-<\mathcal{F}>_{\mathcal{A}}+A \nonumber \\
&= & (p+2)A-pA-(<E>\otimes V^{\otimes n} \otimes <\mathcal{T}>_{\mathcal{A}})\nonumber \\
&= & 2A-(<E>\otimes V^{\otimes n} \otimes <\mathcal{T}>_{\mathcal{A}}) \;. \nonumber
\end{eqnarray}
Since $M$ and $N$ have the same first still, $A$ is also the common Kauffman
amplitude of all stills of $N$.
A precisely similar computation, now shows that the right-hand movie $N$ has the same
$\mathcal{A}$-amplitude as does $M$:
\begin{eqnarray}
<N>_{\mathcal{A}}&=&(q+2)A-qA-(<E>\otimes Id_{V^{\otimes (n+target(Fp))}})\circ
(Id_{V^{\otimes (source(E)+n)}} \otimes <\mathcal{T}>_{\mathcal{A}}) \nonumber \\
&=&2A-(<E>\otimes V^{\otimes n} \otimes <\mathcal{T}>_{\mathcal{A}}) \;. \nonumber
\end{eqnarray}
Hence, as asserted, $\mathcal{A}$ respects every movie-move in $\mathcal{L}_{i}$
---provided $i=$ 17, 18 or 20.\\
\underline{Ad Case 2:} Here $i=19$, so $I(i)=1$, i.e. $\mathcal{T}$ is of type $ET1$.
There is a slight extra complication in this case, because Reidemeister I moves do not
preserve the Kauffman amplitude. Thus, instead of all stills in $M$ and in $N$
having the same Kauffman amplitude (as in Case 1) there are now among these stills
{\bf two} Kauffman amplitudes. Namely, in $M$ the $q+1$ stills involving $F$'s all have the
same Kauffman amplitude $<S>$ as does the common source $S$ of $M$ and $N$,since they
arise from $S$ (the first still of $M$) by a sequence of ET8's
(which preserve the Kauffman bracket), leaving the one still $<T>$ in $M$ with a possibly
different Kauffman amplitude. Similarly, in the right-hand movie $N$,the $q+1$
stills involving $G's$ all have Kauffman amplitude $<T>$, and
there is one remaining still with amplitude $<S>$. (Note: Since here
$\mathcal{T}$ is of type 1, $(p,q)$ is $(2,0)$ or $(0,2)$. Using this fact seems
unhelpful---it seems to require
unnecessarily subdividing the argument into lots of cases. It seems neater to
stick with $p$ and $q$...)

Next let us consider the $\mathcal{A}$-amplitudes of the flickers of $M$ and $N$.
Since $\mathcal{A}$ is semi-normal, the flickers of type 8 divide as follows:
In $M$ there are $q$ of $\mathcal{A}$-amplitude $<S>$, and in $N$ there are $p$
of $\mathcal{A}$-amplitude $<T>$. Finally, $M$ contains one flicker $\mathcal{F}$
(the last) of type ET1, and $N$ contains one flicker $\mathcal{F}'$ (the first)
of type ET1. Exactly as in Case 1, one verifies that $\mathcal{F}$ and $\mathcal{F'}$
have the same $\mathcal{A}$-amplitude:
$$<\mathcal{F}>_\mathcal{A}=<\mathcal{F}'>_\mathcal{A}=
<E>\otimes V^{\otimes n} <\mathcal{T}>_{\mathcal{A}}$$

Putting these facts together, we finally obtain:
\begin{eqnarray}
<\mathcal{M}>_\mathcal{A}&=&(p+1)<S>+<T>-p<S>-<\mathcal{F}>_\mathcal{A}\nonumber \\
&=& <S>+<T>-<E>\otimes V^{\otimes n} <\mathcal{T}>_{\mathcal{A}}\nonumber \\
&=&<\mathcal{N}>_\mathcal{A}\nonumber
\end{eqnarray}
i.e. in both cases, $\mathcal{A}$ respects all movie-moves $(M,N)$ in $\mathcal{L}_i$, as
was to be proved.
\smallskip

\underline{{\bf STEP TWO:} $\mathcal{A}$ {\bf RESPECTS ALL MOVIE-MOVES IN}
 $\mathcal{L}'_{i}$:}\\
Given the $i$-datum $\mathcal{D}_i(E,n,\mathcal{T})$ described above, a second movie-move
 $$\mathcal{M}'\mathcal{D}_i(E,n,\mathcal{T})=(M1,N1)$$
 (these then making up the collection $\mathcal{L}'_{i}$)
is constructed as follows:\footnote{As a check, note that if $m$ is reflection in
the vertical midline, i.e. is the symmetry which interchanges left and right, then
$(M1,N1)~=~(m\circ N~,~m\circ M)$ If we apply this reflection to Fig. \ref{fig:L_i_movies}, we obtain
the diagram for $\mathcal{L}'_i$.}
\begin{verbatim}
M1=[
[0,F1,n+a]...[0,F(p-1),n+a]s[0,Fp,n+a][n+target(Fp),E,0]s=>8
[0,F1,n+a]...s[0,F(p-1),n+a]f[n+target(F(p-1)),E,0]s[0,Fp,n+b]f=>8
......................................................
.......................................................=>8
f[n+source(F1),E,0]s[0,F1,n+b]f...[0,F(p-1),n+b][0,Fp,n+b]s=>I(i)
[n+source(G1),E,0]f[0,G1,n+b]...[0,G(q-1),n+b][0,Gq,n+b]f
]
and
N1=[
s[0,F1,n+a]...[0,F(p-1),n+a][0,Fp,n+a]s[n+target(Fp),E,0]=>I(i)
f[0,G1,n+a]...[0,Q(q-1),n+a]s[0,Gq,n+a]f[n+target(Gq),E,0]s=>8
[0,G1,n+a]...s[0,Q(q-1),n+a]f[n+target(Gq),E,0]s[0,Gq,n+b]f=>8
......................................................
......................................................=>8
[n+source(G1),E,0]f[0,G1,n+b]...[0,G(q-1),n+b][0,Gq,n+b]f
]
\end{verbatim}
(The instructions in [BL].p.49 tell us to obtain the latter movie-moves from the former
``relations"---i.e., movie-moves---
 ``where ... $N^*_{Y_{m,n},\chi_{j,k}}$ is replaced by $N_{\chi_{j,k},Y_{m,n}}$",
Translation into the present notation gives precisely the preceding M1,N1.
It is readily verified, by  amplitude-computations similar to the preceding ones, that
(with notation as above)
$$<M1>_{\mathcal{A}}=<N1>_{\mathcal{A}}=$$
 $$=
\left\{
\begin{array}{ll}
2A-<\mathcal{T}>_{\mathcal{A}}\otimes Id_{(V^{\otimes n})}\otimes <E> & \mbox{ when
$i$=17,18 or 20}\\
<S>+<T>-<\mathcal{T}>_{\mathcal{A}}\otimes Id_{(V^{\otimes n})}\otimes <E> & \mbox{ when
$i$=19}
\end{array}
\right.
$$
Thus, as asserted, $\mathcal{A}$ respects every movie $(M1,N1)$ in $\mathcal{L}'_{i}$.
\smallskip

\underline{{\bf STEP THREE:} $\mathcal{A}$ {\bf RESPECTS }
 $\mathcal{SYM}\circ \mathcal{L}_{i}$ {\bf AND} $\mathcal{SYM}\circ \mathcal{L}'_{i}$}
\smallskip

We are finally ready to tackle the general problem, of proving that every semi-normal
amplitude assignment $\mathcal{A}$, respects every movie-move of
type 17, 18, 19 or 20. For this
purpose, we use the construction of movie-moves of these four types, given in [BL], p.49,
augmented by the addenda on the bottom of loc.cit., p.50. Translated into the terminology of
the present paper, this characterizes the movie-moves of one of these types $i$,
 as the set of all pairs of
movies, derived in the extended sense (cf. Def.\ref{def:extendedDerived})
from the movie-moves in the set
$$ \mathcal{L}_{i} \bigcup  \mathcal{L}'_{i}\;.$$
---or equivalently, as the set of all pairs of movies, strictly derived
(cf. Def.~\ref{def:strictDerived}) from the movie-moves in the set
\begin{equation}\label{eqn:LUL}
(\mathcal{SYM}\circ \mathcal{L}_{i})\; \bigcup\;(\mathcal{SYM}\circ \mathcal{L}'_{i})
\end{equation}
Hence, by Prop \ref{prop:respDerived}, to prove that $\mathcal{A}$ respects all movie-moves
of such types $i$, it suffices to show that it respects all elements of (\ref{eqn:LUL}).

For this purpose, let us first note that the set (\ref{eqn:LUL}) is made up of the lists
$\mathcal{L}_{i}$ and $\mathcal{L}'_{i}$ just described above, together with
two more lists $\mathcal{L}''_{i}$
and $\mathcal{L}'''_{i}$ whose construction is next to be explained:
$$
(\mathcal{SYM}\circ \mathcal{L}_{i})\; \bigcup\;(\mathcal{SYM}\circ \mathcal{L}'_{i})=
\mathcal{L}_{i}\bigcup\mathcal{L}'_{i}\bigcup\mathcal{L}''_{i}\bigcup\mathcal{L}'''_{i}$$

\begin{figure}
\centering
\includegraphics[scale=0.5]{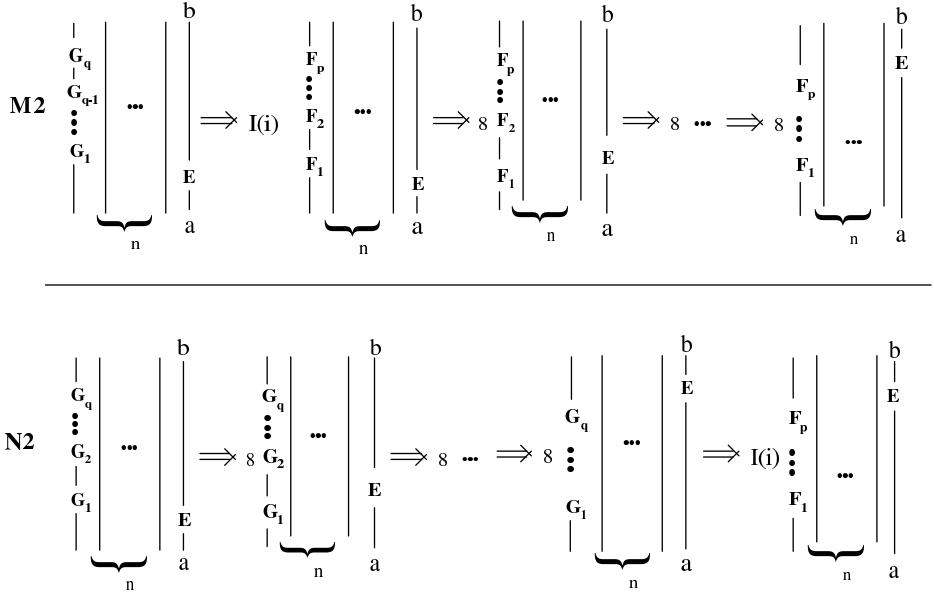}
\caption{\label{fig:L_i_2prime_movies} $\mathcal{L}_{i}''$}
\end{figure}

Figure \ref{fig:L_i_2prime_movies}  gives a sketchy picture of the movie-move
$$\mathcal{M}''\mathcal{D}_i(E,n,\mathcal{T})=(M2==N2)$$. This
movie-move and the related movie-move\footnote{
If $m$ is reflection in the vertical midline,
i.e. is the symmetry which interchanges left and right, then
$(M3,N3)~=~(m\circ N2~,~m\circ M2)$ If we apply this reflection to
Fig. \ref{fig:L_i_2prime_movies}, we obtain the diagram for
$\mathcal{L}'''_i$.}
$$\;\mathcal{M}'''\mathcal{D}_i(E,n,\mathcal{T})=(M3==N3)$$
make up (respectively)
$$\mathcal{L}''_{i}=t\circ\mathcal{L}_i=R\circ\mathcal{L}_i$$
and
 $$\mathcal{L}'''_{i}=t\circ\mathcal{L}'_i=R\circ\mathcal{L}'_i\;.$$

Examination of the figure for $\mathcal{L}''_{i}$ (and a similar
figure for $\mathcal{L}'''_{i}$)  shows that we have
$$\mathcal{M}''\mathcal{D}_i(E,n,\mathcal{T})=
(R\circ \mathcal{M}\mathcal{D}_i(E,n,\mathcal{R\circ T}))^C$$
and
$$\mathcal{M}'''\mathcal{D}_i(E,n,\mathcal{T})=
(R\circ \mathcal{M}'\mathcal{D}_i(E,n,\mathcal{R\circ T}))^C$$
which implies
\begin{equation}\label{eqn:C}
\mathcal{L}''_i=(\mathcal{L}_i)^C \mbox{ and }\mathcal{L}'''_i=(\mathcal{L}'_i)^C
\end{equation}
Now observe that if $\mathcal{A}$ respects a movie-move $(M==N)$, it of course respects
the converse movie-move $(N=M)$. Since we have proved that $\mathcal{A}$ (assumed
semi-normal) respects all movie-moves in $\mathcal{L}_i$ and in $\mathcal{L}'_i$,
it follows from eqn. \ref{eqn:C} that $\mathcal{A}$ also respects all movie-moves
 in $\mathcal{L}''_{i}$ and in $\mathcal{L}'''_{i}$. As noted above, this shows
$\mathcal{A}$ respects all movie-moves of type $i$, \\
Q.E.D.
\bigskip

Let $A=\{15,16,17,18,19,20,22\}$, and let $B$ be the set of
integers between 1 and 30 which are not in $A$. Thus,
\begin{equation}\label{eqn:listB}
B:=\{i\; |\; 1\leq i\leq 14\}\bigcup\{\; 21\; \}\bigcup \{i\; |\; 23\leq i\leq 30 \}
\end{equation}
 Let us continue to assume that $\mathcal{A}$ is a semi-normal
amplitude assignment. It follows from Props.\ref{prop:semiRespects
} and \ref{prop:semi17to20 }, that $\mathcal{A}$ is respects all
movie-moves  whose type lies in $A$. The following proposition
shows that if $i\in B$, a finite number of tests suffice to check
whether $\MA$ respects the (infinitely many) movie-moves of type
$i$; more precisely:
\begin{Prop}\label{finiteTests}
Let $\mathcal{C}$ denote the set of all BL-reduced movie-moves whose type lies in $B$,
and let $\mathcal{A}$ be a semi-normal amplitude-assignment. Then:\\
a) $ \mathcal{C}$, hence also $\mathcal{SYM} \circ \mathcal{C}$, is a finite set.\\
b) If $i\in B$, then $\mathcal{A}$ respects all movie-moves of
type i, if and only if it respects each of the finitely many
movie-moves in
$$(\mathcal{SYM} \circ \mathcal{C})\cap \MMM_i=\mathcal{SYM} \circ (\mathcal{C}\cap \MMM_i) $$.
\end{Prop}
\smallskip

\noindent \underline{{\bf PROOF:}}\\
\underline{Ad a)}  Here is a table (easily obtained by examining
the lists of movie-moves on p.47--50 of [BL]) giving for all $i$
in $B$, the cardinality---which we shall denote by
$\mathcal{N}(i)$--- of the set $\mathcal{MM}^{red}_i$ (of all
BL-reduced movie-moves of type $i$):
\smallskip

\begin{tabular}{|c|l|l|l|l|l|l|l|l|l|l|l|l|l|l|l|l|}
\hline
i in B  & 1 & 2&3&4&5&6&7&8&9&10&11&12&13&14&21&23 \\ \hline
$\mathcal{N}(i)$ & 1 & 1&1&1&6&6&2&1&1&1&2&2&1&3&1&2 \\ \hline
\end{tabular}

\begin{tabular}{|c|l|l|l|l|l|l|l|}
\hline
i in B & 24 & 25 & 26 & 27 & 28 & 29 & 30 \\ \hline
$\mathcal{N}(i)$&1&2&1&1&1&1&1 \\ \hline
\end{tabular}
\smallskip

Since $\mathcal{N}(i)$ is finite for every $i\in B$, it follows that $\mathcal{C}$
is a finite set, of cardinality
$$\Sigma_{i\in B}\mathcal{N}(i)=40$$
Hence also $\mathcal{SYM} \circ \mathcal{C}$ is finite, with cardinality
bounded above by
$$\# (\mathcal{SYM})\times 40=320.$$
\underline{Ad b)}  This follows immediately from Prop.\ref{prop:respDerived}.
\medskip

Up to this point, we have left open the question of what
conditions are needed, for an amplitude-assignment to respect all
movie-moves of type 31. This requires special considerations, to
which we now turn.
\subsection{Reduced Movie-moves of Type 31}\label{SS:31a}
In the following proposition, it will be convenient to use $<<m>>$
as an abbreviation for $Id_{V_1^{\otimes m}}$ (where $m$ is any
natural number.)
\begin{Prop}\label{prop:31A}
Let $\mathcal{A}$ be an amplitude-assignment; let
$$E_1:U_1=>V_1\mbox{ and }E_2:U_2=>V_2$$
be two elementary transitions and let
$$m_1,n_1,m_2,n_2$$
be natural numbers such that (\ref{compos31}) holds---i.e., such
that\\
$m_1+out(E_1)+n_1=m_2+in(E_2+n_2.$

 Construct the movie-move
$$\mathcal{M}(E_1,E_2;m_1,n_1,m_2,n_2)=(M_1==M_2)$$
(as in Definition 31.) Then:\\
(a) $\mathcal{A}$ respects $(M_1==M_2)$, if and only if the
following equation holds:
\begin{eqnarray}
& & (<<m_2>>\otimes <U_2>\otimes
 <<n_2>>)\circ (<<m_1>>\otimes <V_1>\otimes <<n_1>>)-\hspace{1cm} \label{eqn:resp31} \\
& & (<<m_2>>\otimes
<V_2>\otimes <<n_2>>)\circ (<<m_1>>\otimes <U_1>\otimes <<n_1>>)= \nonumber \\
& & (<<m_2>>\otimes
<U_2>\otimes <<n_2>>)\circ (<<m_1>>\otimes <E_1>_{\mathcal{A}}\otimes <<n_1>>)+ \nonumber \\
& &  (<<m_2>>\otimes
<E_2>_{\mathcal{A}}\otimes <<n_2>>)\circ (<<m_1>>\otimes <V_1>\otimes <<n_1>>)- \nonumber \\
& & (<<m_2>>\otimes
<E_2>_{\mathcal{A}}\otimes <<n_2>>)\circ (<<m_1>>\otimes <U_1>\otimes <<n_1>>)- \nonumber \\
& & (<<m_2>>\otimes <V_2>\otimes <<n_2>>)\circ (<<m_1>>\otimes
<E_1>_{\mathcal{A}}\otimes <<n_1>>) \nonumber
\end{eqnarray}
(b) If $\mathcal{A}$ is semi-normal, and at least one of $E_1,E_2$
is of type $ET8$, then (\ref{eqn:resp31}) holds (and hence
$\mathcal{A}$ respects the movie-move $(M_1==M_2))$.
\end{Prop}
\smallskip

\noindent \underline{{\bf PROOF:}}\\
The left movie $M_1$ consists (in movie-time ordering, and using
the $sf$-notation) of the three stills $S_1,S_2,S_3$ given by:
\begin{eqnarray*}
S_1&=&s(m_1\bullet U_1\bullet n_1)s\circ (m_2\bullet U_2\bullet n_2) \\
S_2 &=& f(m_1\bullet V_1\bullet n_1)f\circ s(m_2\bullet U_2\bullet n_2)s \\
S_3 &=& (m_1\bullet V_1\bullet n_1)\circ f(m_2\bullet V_2\bullet
n_2)f
\end{eqnarray*}
joined by the two flickers $F_1,F_2$ given by
\begin{eqnarray*}
F_1 &=& (m_1\bullet E_1\bullet n_1)\circ (m_2\bullet U_2\bullet n_2)\\
F_2 &=& (m_1\bullet V_1\bullet n_1)\circ (m_2\bullet E_2\bullet
n_2)
\end{eqnarray*}
We then have
$$<M_1>_{\mathcal{A}}=<S_1>+<S_2>+<S_3>-<F_1>_{\mathcal{A}}-<F_2>_{\mathcal{A}}
$$
Similarly, the right movie $M_2$ consists of the three stills:
\begin{eqnarray*}
S'_1&=&  S_1\\
S'_2 &=& s(m_1\bullet U_1\bullet n_1)s\circ f(m_2\bullet V_2\bullet n_2)f \\
S'_3 &=& S_3
\end{eqnarray*}
and the two flickers
\begin{eqnarray*}
F'_1 &=& (m_1\bullet U_1\bullet n_1)\circ (m_2\bullet E_2\bullet n_2)\\
F'_2 &=& (m_1\bullet E_1\bullet n_1)\circ (m_2\bullet V_2\bullet
n_2)
\end{eqnarray*}
so that
$$<M_2>_{\mathcal{A}}=<S_1>+<S'_2>+<S_3>-<F'_1>_{\mathcal{A}}-<F'_2>_{\mathcal{A}}$$
Since the terms $<S_1>$ and $<S_3>$ are common to these two
movie-amplitudes, we have
$$<M_1>_{\mathcal{A}}-<M_2>_{\mathcal{A}}=<S_2>-<S'_2>
-<F_1>_{\mathcal{A}}-<F_2>_{\mathcal{A}}+<F'_1>_{\mathcal{A}}+<F'_2>_{\mathcal{A}}$$
Replacing these amplitudes by their explicit values \footnote{The
reader is reminded that (with our present conventions)
still-composition occurs from the bottom to the top ('past' to
`future') of link-diagrams, and that the Kauffman bracket is then
\emph{contravariant}, i.e. $<S_1\circ S_2>=<S_2>\circ <S_1>$.}now
yields eqn. (\ref{eqn:resp31}),
 as the
necessary and sufficient condition that
$$<M_1>_{\mathcal{A}}=<M_2>_{\mathcal{A}}\;,$$ i.e., that
$\mathcal{A}$ respect $(M_1==M_2)$---which proves (a).
\smallskip

 We now turn to (b). Assume $\mathcal{A}$ is semi-normal, and consider the case
that $E_1$ is of type $ET8$ (--- we omit the precisely similar
argument for the case that it is $E_2$ that is of type $ET8$.)
This assumption implies that
$$<U_1>=<V_1>=<E_1>_{\mathcal{A}}:V^{\otimes in(E_1)}\rightarrow V^{\otimes out(E_1)}$$
This implies the equality of the 3 $\KK$-linear maps
$<<m_1>>\otimes <U_1>\otimes <<n_1>>,<<m_1>>\otimes <V_1>\otimes
<<n_1>>$\\
and $<<m_1>>\otimes <E_1>_{\mathcal{A}}\otimes <<n_1>>$ ---let us
denote their common value by
$$\mathcal{S}:V^{\otimes (m_1+in(E_1)+n_1)}\rightarrow V^{\otimes (m_1+out(E_1)+n_1)}\;.$$
Then (\ref{eqn:resp31}) is, in this case, the result of applying
$\circ \mathcal{S}$ on the right of the obviously correct equation
\begin{eqnarray*}
<<m_2>>\otimes <U_2>\otimes <<n_2>>-<<m_2>>\otimes <V_2>\otimes <<n_2>>=\\
<<m_2>>\otimes <U_2>\otimes <<n_2>>+<<m_2>>\otimes <E_2>\otimes <<n_2>>-\\
<<m_2>>\otimes <E_2>\otimes <<n_2>>-<<m_2>>\otimes <V_2>\otimes
<<n_2>>\;,
\end{eqnarray*}
Q.E.D.
\smallskip

We now turn to the definition of {\bf reduced} movie-moves of type
31---these will then be defined to make up the finite set
$\mathcal{MM}^{red}_{31}$. This is analogous to (but rather more
complicated than) the above concept of $BL$-reduced (for
movie-moves of types $\leq 30$), and is chosen so that the
analogue of Prop.\ref{finiteTests} holds (this analogue being
Prop.\ref{prop:Assoc31} below.)
\medskip

\begin{figure}
  \centering
  \includegraphics[scale=0.5]{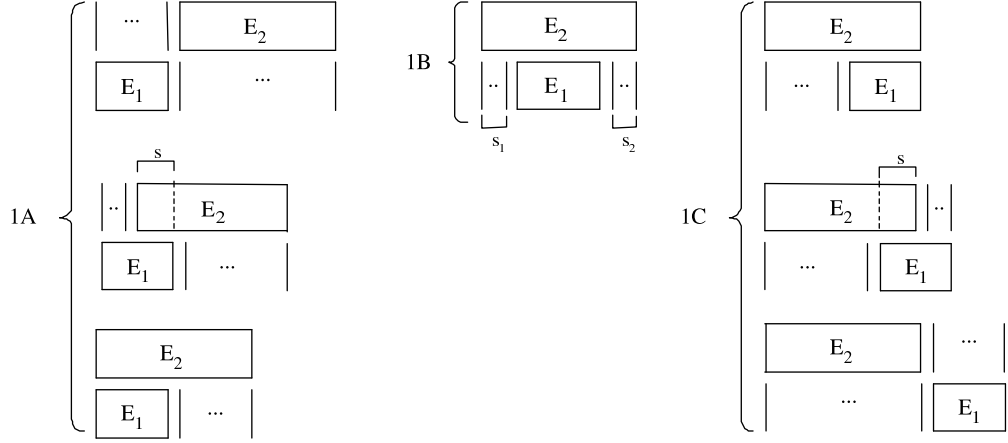}
  \caption{$\mathcal{MM}^{red}_{31}$}\label{figMm31}
\end{figure}

\noindent \underline{{\bf ALGORITHM A (FOR GENERATING THE FINITE
SET $\mathcal{MM}^{red}_{31}$: })} Using the notation of
Def.\ref{MM31}, the \underline{{\bf reduced}} movie-moves
$$\mathcal{M}=\mathcal{MM}(E_1,E_2:m_1,n_1,m_2,n_2) \eqno(*)$$ of type 31,
are now defined to be all those obtained by the following
procedure:

In  the first place, choose in all possible 1596 ways, two
elementary transitions $E_1,E_2$ (not necessarily distinct), such
that neither $E_1$ nor $E_2$ is of type $ET8$. Let us set
$$t_1:=target(E_1),\; t_2=source(E_2)\;.$$

At this point, we must subdivide into 3 disjoint cases: \\
\underline{{\bf Case 1: $t_1 < t_2$}}\hspace{.25in}
\underline{{\bf Case 2:
$t_1 = t_2$}}\hspace{.25in} \underline{{\bf Case 3: $t_1 > t_2$}}\\
In each case, we construct $t_1+t_2+1$ quadruples
$(m_1,n_1,m_2,n_2)$ as described below, insert into $(*)$ to
obtain a set---call it $\mathcal{MM}^{red}_{31}(E_1,E_2)$---
consisting of $t_1+t_2+1$ movie-moves of type 31. Finally, we set
$$\mathcal{MM}^{red}_{31}:=\bigcup \{\mathcal{MM}^{red}_{31}(E_1,E_2)
 | E_1,E_2 \mbox{ elementary transitions neither of type }ET8\}$$
\smallskip

In more detail:\\
 Suppose first we are in Case 1, so that $t_1<t_2$. We subdivide
 into three parametrized SUB-CASES, each providing a value for the quadruple
$(m_1,n_1,m_2,n_2)$, as follows:\\
SUB-CASE 1A(s) (parametrized by $s$ with $0\leq s\leq t_1$) \\
Here $(m_1,n_1,m_2,n_2)=(0,t_2-s,t_1-s,0)$\\
SUB-CASE 1B$(s_1,s_2)$ (parametrized by $s_1,s_2$ with
$$0<s_1,0<s_2,s_1+s_2=t_2-t_1\;)$$---here
 $(m_1,n_1,m_2,n_2)=(s_1,s_2,0,0)$.\\
 SUB-CASE 1C(s) (parametrized by $s$ with $0\leq s\leq
 t_1$)\\
Here $(m_1,n_1,m_2,n_2)=(t_2-s,0,0,t_1-s)$)\\
\smallskip

The number of quadruples $(m_1,n_1,m_2,n_2)$ produced in this way
is
$$(t_1+1)+(t_2-1+1)+(t_1+1)=t_1+t_2+1$$
as asserted above; inserting these quadruples into $(*)$, we
obtain the $t_1+t_2+1$ reduced movie-moves of type 31 in
$\mathcal{MM}^{red}_{31}(E_1,E_2)$.

We next consider Case 2, where $t_1=t_2$. Here we have three
SUB-CASEs, as follows: \\
SUB-CASE 2A (parametrized by $s$ with $0\leq s<t1$)\\
$(m_1,n_1,m_2,n_2)=(0,t_1-s,t_1-s,0)$\\
SUB-CASE 2B $(m_1,n_1,m_2,n_2)=(0,0,0,0)$ \\
SUB-CASE 2C{s} (parametrized by $s$ with $0\leq s<t_1$)\\
$(m_1,n_1,m_2,n_2)=(t_1-s,0,0,t_1-s)$\\
Inserting these values into $(*)$, we obtain
$\mathcal{MM}^{red}_{31}(E_1,E_2)$ in the present case.

There remains Case 3, where $t_1>t_2$. Here we have the three
SUB-CASEs:\\
 SUB-CASE 3A (parametrized by $s$ with $0\leq s \leq t_2$)\\
$(m_1,n_1,m_2,n_2)=(0,t_2-s,t_1-s,0)$\\
SUB-CASE 3B (parametrized by $s_1,s_2$ with
$$0\leq s_1,0\leq s_2,s_1+s_2=t_1-t_2\;)$$
$(m_1,n_1,m_2,n_2)=(0,0,s_1,s_2)$\\
SUB-CASE 3C (parametrized by $s$ with $0\leq s \leq t_2$)\\
$(m_1,n_1,m_2,n_2)=(t_2-s,0,0,t_1-s)$ \\
Inserting these values
into $(*)$, we obtain $\mathcal{MM}^{red}_{31}(E_1,E_2)$ in the
present case.
\smallskip

The motivation for  the preceding algorithm, is the need to ensure
that the following proposition be true---the algorithm may appear
more natural, upon carefully examining the following proof.
\begin{Prop}\label{prop:Assoc31}
Let $\mathcal{A}$ be a  semi-normal amplitude-assignment; then
$\mathcal{A}$ respects all movie-moves of type 31, if and only if
$\mathcal{A}$ respects all {\bf reduced} movie-moves of type 31
(i.e., all elements in the finite set $\MMM^{red}_{31}).$
\end{Prop}
\smallskip

\noindent \underline{{\bf PROOF:}}\\
Let us make the following three assumptions:
\begin{enumerate}
\item $\mathcal{A}$ is a semi-normal amplitude-assignment.
\item $\mathcal{M}$ is a movie-move of type 31
\item $\mathcal{A}$ respects every reduced movie-move of type 31.
\end{enumerate}
It suffices to deduce from these assumptions, that $\mathcal{A}$
respects $\mathcal{M}$.
\smallskip

In the sense of Def.\ref{def:strictDerived}, $\mathcal{M}$ is
strictly derived from a movie-move $\mathcal{M}'$ of the form
$$\mathcal{M}'=\mathcal{M}(E_1,E_2;m_1,n_1,m_2,n_2)$$
where $E_1,E_2$ are elementary transitions, and the m's and n's
are natural numbers. It follows from Prop.\ref{prop:respDerived},
that it suffices to deduce from our three assumptions, that
$\mathcal{A}$ respects $\mathcal{M}'$. Moreover, it follows from
Prop.\ref{prop:31A}(b), that we may assume without loss of
generality, that neither $E_1$ nor $E_2$ is of type ET8.
\smallskip

It follows from the preceding paragraph, that we may, without loss
of generality, strengthen Assumption 2 as follows:\\
 \hspace*{.06in} $2'$.\hspace{.07in}$\mathcal{M}$ has the form
$$\mathcal{M}=\mathcal{M}(E_1,E_2;m_1,n_1,m_2,n_2)$$
with neither $E_1$ nor $E_2$ of type ET8.
\medskip

It suffices to prove, (as we shall in the remainder of our
argument), that given Assumptions 1, $2'$ and 3, it follows that
there exists a reduced movie-move $\mathcal{M}^{red}$, such that
\begin{equation}\label{eqn:key}
\mathcal{A} \mbox{ respects } \mathcal{M}^{red} \Rightarrow
\mathcal{A} \mbox{ respects } \mathcal{M}
\end{equation}
\medskip

\begin{figure}
  \centering
  \includegraphics[scale=0.75]{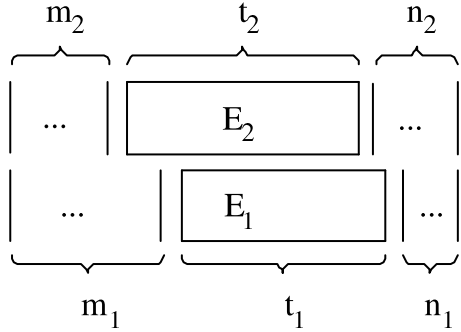}
  \caption{Schematic Construction of $\mathcal{M}$}\label{figMmproofA}
\end{figure}

 Now let us set
$$t_1=target(E_1),\; t_2=source(E_2)$$
Note that assumption $2'$ then implies
\begin{equation}\label{eqn:mnt}
m_1+t_1+n_1=m_2+t_2+n_2
\end{equation}
 In order  to
understand the rationale of the somewhat elaborate case-divisions
which follow, the reader may find it helpful to examine Figure
\ref{figMmproofA}; the case-divisions are according to the
possible
relations of the portions of this Figure.\\

\begin{figure}
  \centering
  \includegraphics[scale=0.6]{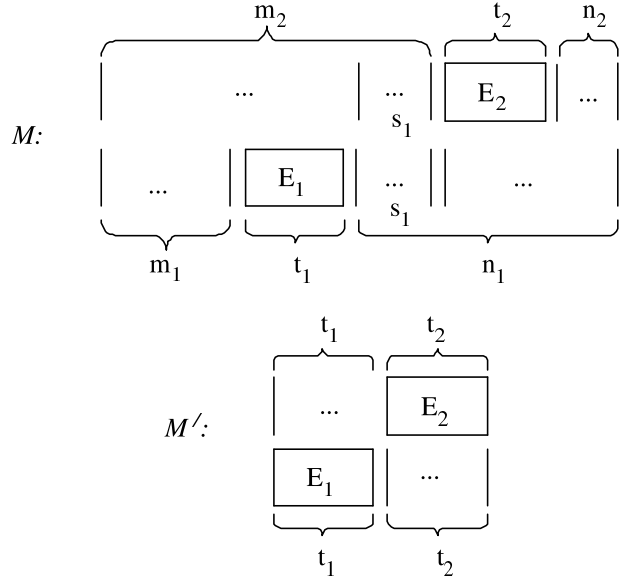}
  \caption{Sub-case 1a}\label{figMm31ProofA1}
\end{figure}

 We begin by again considering (as
in Algorithm A) three cases, according as
$t_1$ is $<,=$ or $> t_2$.\\
\underline{{\bf CASE 1 } $\;t_1<t_2$:}\\
We subdivide this case further, into five disjoint sub-cases, as
follows:\\
\underline{{\bf SUB-CASE 1a   } $\;m_2\geq m_1+t_1:$}\\
(In terms of Figure \ref{figMmproofA}, this means that all of
$E_1$ is to the left of all of $E_2$: cf. Fig.
\ref{figMm31ProofA1}). We must construct a reduced movie-move
$\mathcal{M}'$ such that (\ref{eqn:key}) holds. It is now claimed
that this is the case for
$$\mathcal{M}'=\mathcal{M}(E_1,E_2;0,t_2,t_1,0)$$
Indeed, examination of Algorithm A shows that $\mathcal{M}'$ is
reduced (it falls under SUB-CASE 1A(0) of the Algorithm). It
remains to verify (\ref{eqn:key}):\\

Thus, assume that $\mathcal{A}$ respects $\mathcal{M}'$. In (a) of
Prop. \ref{prop:31A}, with $m_1,n_1,m_2,n_2$ replaced by their
values $0,t_2,t_1,0$ for $\mathcal{M}'$, the first term of
Eqn.\ref{eqn:resp31} becomes
$$(<<t_1>>\otimes <U_2>)\circ (<V_1\otimes <<t_2>>)=<V_1>\otimes
<U_2>$$ ---and similarly for the other 5
 terms. We thus obtain
\begin{equation}\label{eqn:1a}
T_1'\otimes U_1'-T_2'\otimes U_2'=T_3'\otimes U_3'+T_4'\otimes
U_4'-T_5'\otimes U_5'-T_6'\otimes U_6'
\end{equation}
where the  $\KK$-linear transformations
$$T_i':V^{\otimes t_1}\rightarrow V^{\otimes t_1}\mbox{ and }
U_i':V^{\otimes t_2}\rightarrow V^{\otimes t_2}\hspace{.2in}(1\leq
i\leq6)$$ are given by
$$T_1'=<V_1>,T_2=<U_1>,T_3=<E_1>_{\mathcal{A}},T_4=<V_1>,T_5=<U_1>,T_6=<E_1>_{\mathcal{A}}$$
and
$$U_1'=<U_2>,U_2'=<V_2>,U_3'=<U_2>,U_4'=<E_2>_{\mathcal{A}},U_5'=<E_2>_{\mathcal{A}},U_6'=<V_2>$$

On the other hand,setting $s=m_2-m_1-t_2$, (a) of
Prop.\ref{prop:31A}, together with examination of Fig.
\ref{figMm31ProofA1} in this case, shows that $\mathcal{A}$
respects $\mathcal{M}$, if and only if
\begin{equation}\label{eqn:1a1}
T^{(1)}+T^{(2)}=T^{(3)}+T^{(4)}-T^{(5)}-T^{(6)}
\end{equation}
where the $\KK$-linear transformations $T^{(i)}$ are given by:
\begin{eqnarray*}
T^{(1)} &=& (<<m_1>>\otimes <<t_2>>\otimes
<<s_1>>\otimes<U_2>\otimes <<m_2>>)\circ \\
&& (<<m_1>>\otimes <V_1>\otimes <<s_1>>\otimes <<t_2>>\otimes
<<n_2>>)\\
&=& <<m_1>>\otimes <V_1>\otimes <<s_1>>\otimes <U_2>\otimes
<<n_2>>
\end{eqnarray*}
---and similarly,

\begin{eqnarray*}
T^{(2)} &=& <<m_1>>\otimes <U_1>\otimes <<s_1>>\otimes
<V_2>\otimes <<n_2>> \\
T^{(3)} & =&  <<m_1>>\otimes <E_1>_{\mathcal{A}}\otimes
<<s_1>>\otimes <U_2>\otimes
<<n_2>>\\
T^{(4)} &=& <<m_1>>\otimes <V_1>\otimes <<s_1>>\otimes
<E_2>_{\mathcal{A}}\otimes
<<n_2>>\\
T^{(5)} &=& <<m_1>>\otimes <U_1>\otimes <<s_1>>\otimes
<E_2>_{\mathcal{A}}\otimes
<<n_2>>\\
T^{(6)} &=& <<<m_1>>\otimes <E_1>_{\mathcal{A}}\otimes
<<s_1>>\otimes <V_2>\otimes <<n_2>>
\end{eqnarray*}
Clearly, eqn.(\ref{eqn:1a}) implies eqn.(\ref{eqn:1a1}), which
completes the proof in this first Sub-case 1a.\\
\underline{{\bf SUB-CASE 1b   } $\;m_1+t_1>m_2, n_2+t_2>n_1,
m_1\leq m_2
\mbox{ and }n_2\leq n_1$}\\
Here the argument is much simpler than in the preceding Sub-case:\\
Using eqn.\ref{eqn:mnt}, let us set
$$s:=m_1+t_1-m_2 \mbox{ and }n_2+t_2-n_1\;.$$

\begin{figure}
  \centering
  \includegraphics[scale=0.52]{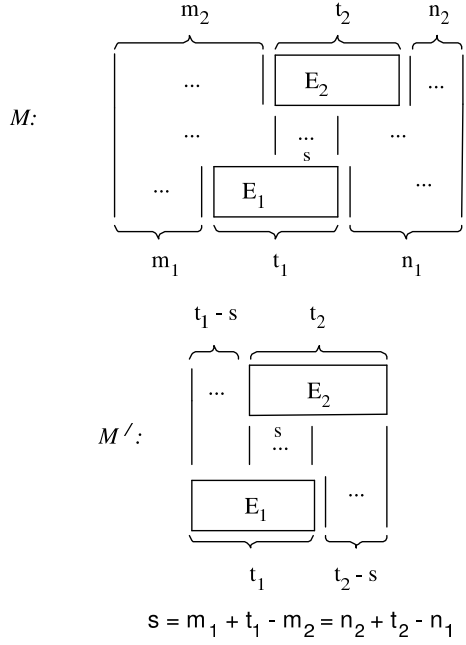}
  \caption{Sub-case 1b}\label{figMm31ProofAb}
\end{figure}

 Examination of Fig. \ref{figMm31ProofAb} shows that here
$$\mathcal{M}=m_1\bullet \mathcal{M}'\bullet n_2$$
where we set
 $$\mathcal{M}':=\mathcal{M}(0,t_2-s,t_1-s,0)$$
The movie-move $\mathcal{M}'$ of type 31, is in fact reduced---
falling under Subcase 1A(s) of Algorithm A. By Assumption 3,
$\mathcal {A}$ respects $\mathcal{M}'$, hence, by (b) of
Prop.\ref{red:comp}, $\mathcal{A}$ respects $\mathcal{M}$---which
proves our proposition holds in the Sub-case 1b.
\smallskip

\begin{figure}
  \centering
  \includegraphics[scale=0.52]{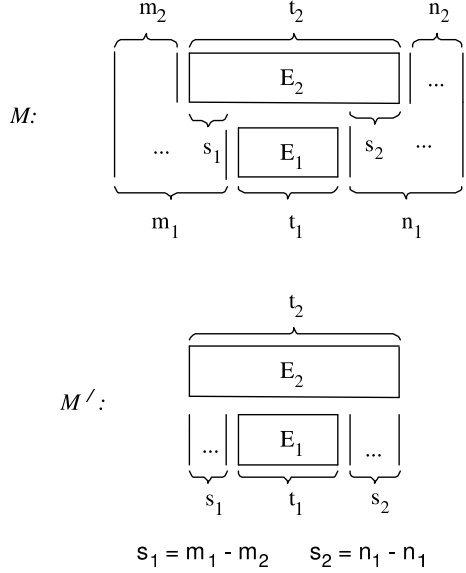}
  \caption{Sub-case 1c}\label{figMm31ProofAc}
\end{figure}

\underline{{\bf SUB-CASE 1c   } $\;m_1>m_2\mbox{ and }n_1> n2$}\\
We then set
$$s_1:=m_1-m_2\mbox{ and }s_2:=n_1-n_2$$
so that $s_1>0, s_2>0$ and (because of eqn.(\ref{eqn:mnt}))
$$s_1+s_2=t_2-t_1\;.$$
Hence the movie-move
$$\mathcal{M}':=\mathcal{M}(E_1,E_2;s_1,s_2,0,0)$$
of type 31, is reduced. (It falls under Sub-case 1B($s_1,s_2$) of
Algorithm A.) Examination of Fig. \ref{figMm31ProofAc} shows that
$$\mathcal{M}=m_1\bullet \mathcal{M}'\bullet n_2\;,$$
and we conclude by precisely the same argument used in the
preceding Sub-case.

\smallskip

The  remaining two Sub-cases of Case 1, are:\\
 \underline{{\bf SUB-CASE 1d   } $\; m_2+t_2>m_1, n_2+t_1>n_1, m_2\leq m_1
 \mbox{ and }n_2\leq n_1$}\\
and\\
\underline{{\bf SUB-CASE 1e   } $\;n_2\geq n_1+t_1$}.\\
These are obtained from Sub-cases 1b and 1a (in that order) by
interchanging left and right (i.e., by the symmetry $m$.) The
related reduced movie-moves are, respectively,
$$\mathcal{M}(E_1,E_2;t_2-s,0,0,t_1-s)\mbox{ where }s=t_1+n_2-n_1=t_2+m_2-m_1$$
and
$$\mathcal{M}(E_1,E_2;t_2,0,0,t_1)\;.$$
\smallskip

Case 3, where $t_1>t_2$, is precisely similar to Case 1 (Also, it
can be obtained from Case 1, by applying the symmetry $t$.)
Finally, Case 2, where $t_1=t_2$, is like Cases 1 and 3.
\subsection{Strongly Normal Amplitude-Assignments}\label{SS:normal}
In the following, we shall be studying certain linear equations
over the ground-field
$$\KK:=\QQ(q)\;.$$

Let $m,n$ be natural numbers. Recall that, in the category
$\Sigma$, defined in \S \ref{SS:Stills}, we denote by
$Hom_{\Sigma}(m,n)$ the set of stills $S$ such that
$$source(S)=m,\;target(S)=n\;, $$
i.e., which have m lines coming in at the bottom, and n lines
going out at the top.  Also, we shall denote by $Map(m,n)$ the set
of all $\KK$-linear maps from $V_1^{\otimes m}$
 to $V_1^{\otimes n}$, where
$$V_1:=V\otimes _{\ZZ[q,q^{-1}]}\KK$$
  If $S$ is a still in
$Hom_{\Sigma}(m,n)$, then by a slight abuse of notation,its
Kauffman amplitude $<S>$ will be identified with the corresponding
element
$$<S>\otimes \KK$$
 in $Map(m,n)$.

\begin{figure}
  \centering
  \includegraphics[scale=0.5]{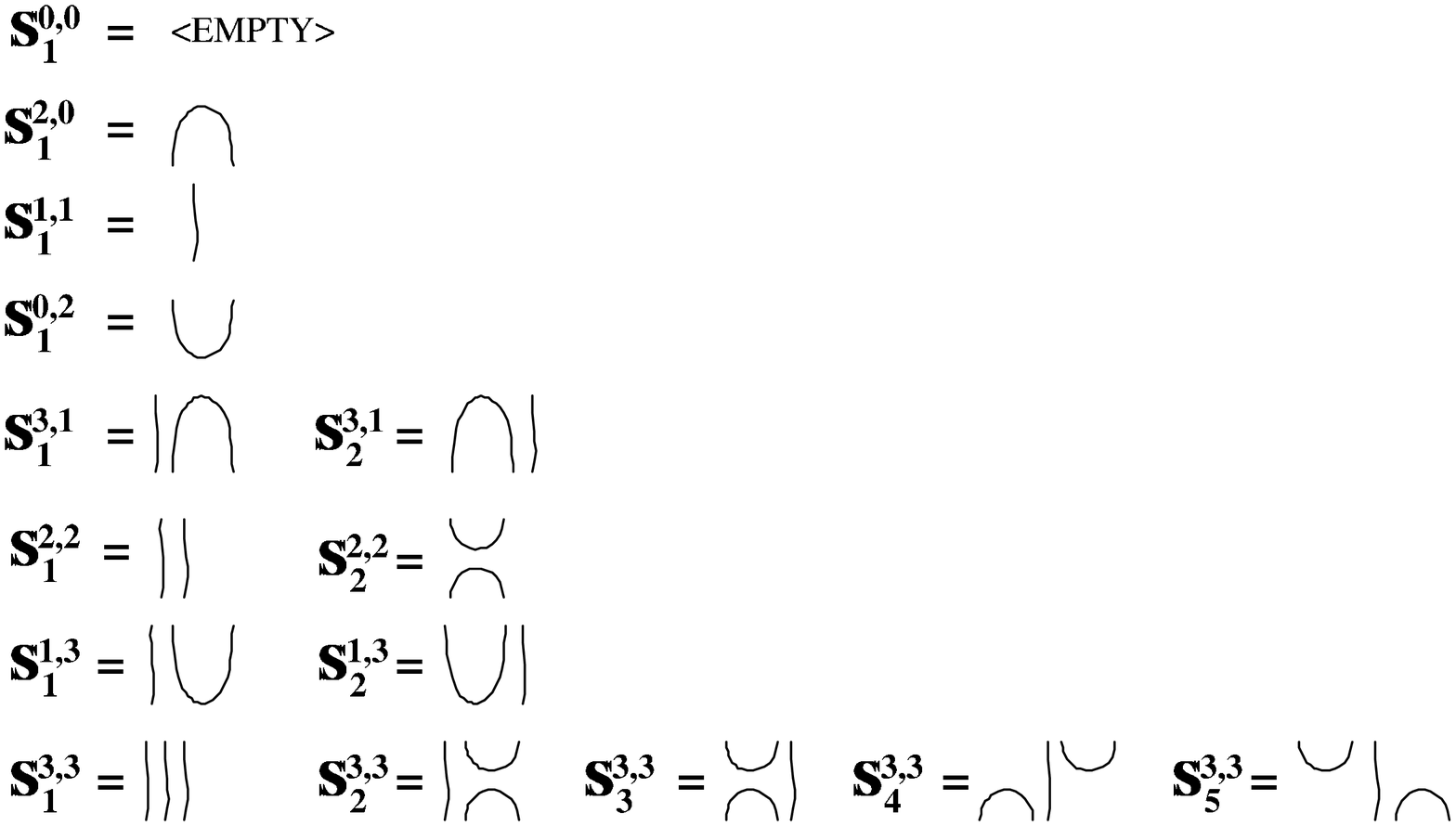}
  \caption{Basic Stills}\label{fig:still_basis}
\end{figure}

For each elementary transition $\mathcal{T}$, of type between 1
and 7, the ordered pair
 $$(in(\mathcal{T}), out(\mathcal{T}))$$
is one of the ordered pairs in the set
$$\mathcal{P}=\{(0,0),(2,0),(1,1),(0,2),(3,1),(2,2),(1,3),(3,3)\}$$
\begin{Def}
To each pair $(m,n)$ in $\mathcal{P}$, we associate a positive
integer $N(m,n)$, and a collection of $N(m,n)$ stills in
$Hom_{\Sigma}(m,n)$, as follows.
(See Fig. \ref{fig:still_basis})\\
\noindent \underline{if $(m,n)=(0,0):$} $N(0,0)=1,\; S_1^{0,0}=1_0$ \\
\underline{if $(m,n)=(2,0)$:} $N(2,0)=1,\; S_1^{2,0}=[0,Cap,0]$\\
\underline{if $(m,n)=(1,1):$} $N(1,1)=1,\; S_1^{1,1}=1_1$\\
\underline{if $(m,n)=(0,2)$:} $N(0,2)=1,\; S_1^{0,2}=[0,Cup,0]$\\
\underline{if $(m,n)=(3,1):$} $N(3,1)=2,\; S_1^{3,1}=[1,Cap,0],\; S_2^{3,1}=[0,Cap,1]$ \\
\underline{if $(m,n)=(2,2,):$} $N(2,2,)=2,\; S_1^{2,2}=1_2,\;
S_2^{2,2}=
[0,Cap,0][0,Cup,0]$ \\
\underline{if $(m,n)=(1,3,):$} $N(1,3)=2,\; S_1^{1,3}=[1,Cup,0],\;S_2^{1,3}=[0,Cup,1]$ \\
\underline{if $(m,n)=(3,3):$} $N(3,3)=5,\; S_1^{3,3}=1_3,\;S_2^{3,3}=[1,Cap,0][1,Cup,0],\\
S_3^{3,3}=[0,Cap,1][0,Cup,1],\;S_4^{3,3}=[0,Cap,1][1,Cup,0],\;S_5^{3,3}=[1,Cap,0][0,Cup,1].$
\end{Def}

\begin{Def}\label{def:normal}
Let $\mathcal{A}$ be an amplitude-assignment. Recall that
 $\mathcal{A}$ is called \underline{{\bf semi-normal}}
when, for every elementary transition $E:S=>T$ of type either 0 or
8,
$$<E>_{\mathcal{A}}=<S>=<T>\;.$$
$\mathcal{A}$ will be called \underline{{\bf strongly normal}} if
it is semi-normal, and if, in addition,
 for every elementary transition
E of type between 1 and 7 inclusive, the linear transformation
$$<E>_{\mathcal{A}}:V_1^{in(E)}=>V_1^{out(E)}$$ satisfies the following condition:\\
Let $m=in(E)$ and $n=out(E)$ (so that $(m,n)$ is an element of $\mathcal{P}$);\\
then $<E>_{\mathcal{A}}$ is a $\KK$-linear combination of the
$N(m,n)$ Kauffman amplitudes
\begin{equation}\label{eqn:basicStills}
\{ <S^{m,n}_j>\;|\;1\leq j\leq N(m,n) \}\;.
\end{equation}
\end{Def}

\bigskip
By imposing this Ansatz of strong normality\footnote{The maps
(\ref{eqn:basicStills}) are in fact $sl_q(2)$-linear. As
motivation for this definition of `strongly normal', we note  that
it is equivalent to the requirement that each $<E>_{\MA}$ be
$sl_q(2)$-linear. We omit the proof, since it is not needed for
any later results, and would add considerably to the length of the
paper, without producing any further invariants.}, the problem of
constructing an amplitude-assignment which is $\MU$-balanced,
becomes, as we shall now see, the feasible problem of solving a
managable number of linear equations (over two thousand, for
instance, if $\MU=\{31\}$) in a manageable number 102 of unknowns.
(The price paid for this desirable reduction,
 is  that we lose track in this paper, of whatever 2-knot
invariants arise via the present constructions from non-strongly
normal amplitude-assignments.)

\medskip

 We may now coordinatize strongly normal amplitude-assignments
$\mathcal{A}$ via 102 parameters in $\KK$, as follows:

If $\mathcal{A}$ is a strongly normal amplitude-assignment, and if
$E$ is an elementary transition of type between 1 and 7
inclusive---say with $m=in(E)$ and $n=out(E)$
--- then, as just proved, there exist uniquely $N(m,n)$ elements
\begin{equation}\label{eqn:parms}
\{PE_j(\mathcal{A})\;|\;1\leq j\leq N(m,n)\}
\end{equation}
in the ground-field $\KK$, such that
\begin{equation}\label{eqn:fund}
<E>_{\mathcal{A}}=\Sigma_{j=1}^{N(m,n)} PE_j(\mathcal{A})<S^{m,n}_j> \;.
\end{equation}
In this  way, we associate to every strongly normal
amplitude-assignment $\mathcal{A}$, a collection
$$\{P^E_j(\mathcal{A}) | E\mbox{ an elementary
transition with }1\leq type(E)\leq 7,\; 1\leq j\leq N(in(E),out(E))\}$$
consisting of
\smallskip

\noindent $\Sigma_{i=1}^7\Sigma_{(E\mbox{ an elementary transition of type } i)}
N(in(E),out(E))$ \\
\vspace{.15in} $=8\times 1+4\times 2+12\times 5+8\times 2+4\times
1+2\times 1+2\times 2 =102$

 \noindent elements of  $\KK$, which we
shall call the $\underline{{\bf normal \; coordinates}}$ of
$\mathcal{A}$.
\smallskip

By means of these normal coordinates, the set of all strongly
normal amplitude-assignments $\mathcal{A}$, is put in bijection
with the free module of rank 102 over $\KK$. The condition that
such $\mathcal{A}$ be $\MU$-balanced, is , as will next be proved,
equivalent to a collection of linear equations over $\KK$ in its
normal coordinates which we shall call the \underline{{\bf
associated equations}}. The final task in the program we have
called Program A, is the computation of the associated equations,
and their complete solution over $\KK$.
\medskip

Let us denote by $\mathcal{PARM}$ the set of 102 formal symbols
$$\{ PE_j\;|\;E\in \mathcal{ET},\;1\leq j\leq N(in(E),out(E))\}\;,$$
and by $\KK.\mathcal{PARM}$ the free $\KK$-module on the set
$\mathcal{PARM}$. Thus,
the `normal coordinates' described above, give the bijection\\
$\Phi:\{\mbox{strongly normal amplitude-assignments
}\mathcal{A}\}\simeq
\KK.\mathcal{PARM}$,\\
$\mathcal{A}\mapsto \Sigma_{E\in \mathcal{ET}}\Sigma_{j=1}^{N((in(E),out(E))}
PE_j(\mathcal{A})PE_j$\\
(where $\mathcal{A}$ may be reconstructed from its image $\Phi(\mathcal{A})$ via
Def.\ref{def:seminormalAA} and Eqn.(\ref{eqn:fund}).)
\subsection{Associated Equations}\label{SS:Assoc}
We next turn to this question:\\
\underline{\bf{Given a flicker F and a strongly normal
amplitude-assignment $\mathcal{A}$,
how do}}\\
\underline{\bf{we compute $<F>_\mathcal{A}$ in terms of
the normal coordinates of $\mathcal{A}$?}}\\
We subdivide the set of flickers into 3 disjoint subsets, which
yield somewhat different answers
to this question. \\
\underline{{\bf CASE 1: $F$ is of type 0 or 8}}\\
With the notation of Def. \ref{defFlicker},
$$F=(E,B,m,n,T)$$
where $E:S'=>S''$ is an elementary transition of type 0 or 8,
$B$ and $T$ are stills (either
or both of which may be empty), $m$ and $n$ are
natural numbers(which may be 0), such that the two following relations are satisfied:
$$target(B)=m+in(E)+n $$
$$source(T)=m+out(E)+n $$
By assumption, $\mathcal{A}$ is strongly normal and hence
semi-normal, so
$$<E>_{\mathcal{A}}=<S'>=<S''>$$
Hence,
$$<F>_\mathcal{A}=<T>\circ (Id_{V_1^{\otimes m}}\otimes <S'>\otimes Id_{V_1^{\otimes n}})
\circ <B> $$
and so finally
$$<F>_\mathcal{A}=<S_F> \eqno{(A_1)}$$
where
$$S_F=T\circ (m\bullet S'\bullet n)\circ B\;. \eqno{(A_2)} $$
(Note that in this case, $<F>_\mathcal{A}=<S_F>$ is independent of $\mathcal{A}$.)\\
\underline{{\bf CASE 2: $F$ is of type $i$ with $1\leq i\leq 7$ }}\\
Here again, $F=(E,B,m,n,T)$
where now $E$ is an elementary transition of type $i$. Let
us assume that $p=in(E)$ and $q=out(E)$; then by Eqn. \ref{eqn:fund},
$$<E>_{\mathcal{A}}=\sum_{j=1}^{N(p,q)}PE_j(\mathcal{A})<S^{p,q}_j>\;.$$
By Def. \ref{Def:ampFlicker},
$$<F>_\mathcal{A}=<T>\circ (Id_{V_1^{\otimes m}}\otimes <E>_\mathcal{A}
\otimes Id_{V_1^{\otimes n}})\circ <B> $$ Combining the two
preceding equations, we get
$$<F>_\mathcal{A}=\sum_{j=1}^{N(p,q)}PE_j(\mathcal{A})<S^{p,q}_j(F)>\;,\eqno(B_1)$$
where the stills $S^{p,q}_j(F)$ (which we note are independent of $\mathcal{A}$)
are given by
$$S^{p,q}_j(F)=T\circ (m\bullet S^{p,q}_j\bullet n)\circ B\;. \eqno{(B_2)} $$
\underline{{\bf CASE 3: $F=1_S$}}\\
Here, by definition,
$$<F>_\mathcal{A}=<S>\;. \eqno{(C)}$$
\medskip

We next consider the analogous question, with flickers replaced  by movies:
\begin{Prop}\label{prop:coordMovieAmp}
Given a movie $M$, with
$$m=in(M),\;n=out(M) $$
there exist the following (depending on $M$, but not on $\mathcal{A}$):\\
\begin{itemize}
\item Natural numbers $p,q$
\item $\epsilon_1,\cdots,\epsilon_p,\epsilon_1',\cdots,\epsilon_q' $ in $\{1,-1\}$
\item Stills $S_1,\cdots,S_p,T_1,\cdots,T_q$ (not necessarily
all distinct) with origin $m$ and terminus $n$
\item Elements (not necessarily all distinct) $(E^1,k^{(1)}),\cdots,(E^q,k^{(q)})$
in $\mathcal{PARM}$
\end{itemize}
\bigskip

such that, for every strongly normal amplitude-assignment
$\mathcal {A}$,\footnote{ If $p$ or $q$ is 0, we use the usual
convention that an empty sum is 0.}
\begin{equation}\label{eqn:coordMovieAmp}
<M>_{\mathcal {A}}=\sum_{i=1}^p \epsilon_i<S_i>+
\sum_{j=1}^q \epsilon_j'(PE^j)_{k^{(j)}}(\mathcal{A}) <T_j>\;.
\end{equation}
\end{Prop}
\smallskip

\noindent \underline{\bf {PROOF:}}\\
 Let us assume the flickers in $M$
are given by
$$M=F_1\cdots F_p\; ,$$
so that the $p+1$ stills in $M$ are given by
$$S_1=source(F_1),\; S_i=source(F_i)=target(F_{i-1})\mbox{ for }1\leq i\leq p,\;
S_{p+1}=target(F_p)\;,$$
Then Def. \ref{def:ampMovie} implies
\begin{equation}\label{eqn:forMovie}
<M>_{\mathcal {A}}=\sum_{i=1}^{p+1}<S_i>-\sum_{i=1}^p<F_i>_{\mathcal{A}}\;.
\end{equation}
Consider the set $\M{E}$ of all $\KK$-linear maps
$$V_1^{\otimes m}\mapsto V_1^{\otimes n}$$
that can be expressed in the form given by the right side of
(\ref{eqn:coordMovieAmp}). If $F$ and $G$ lie in $\M{E}$ then also
$\pm F \pm G $ can be expressed in this form; since each term in
the right side of (\ref{eqn:forMovie}) lies in $\M{E}$, (using
equations
(A), (B), (C) in the preceding discussion), so does $<M>_{\mathcal {A}}$,\\
Q.E.D.
\medskip

\begin{Prop}\label{prop:Assoc}
Let $\mathcal{M}=(M_1,M_2)$ be a grammatical pair of movies; then
there exists a finite collection $($say, of cardinality
$\nu(\mathcal{M}))$
 of inhomogeneous equations in the 102 variables $PE_j$:
\begin{equation}\label{eqn:Assoc}
\sum_{PE_j\in \mathcal{PARM}}C(E,j,k)PE_j=D(E,j,k)\;\;\;\;\;(1\leq k\leq \nu(\mathcal{M})
\end{equation}
having all coefficients $C$ and $D$ in $\KK$---with this property:\\

For every strongly normal amplitude-assignment $\mathcal{A}$, the
system of $\nu(\mathcal{M})$ equations (\ref{eqn:Assoc}) are
necessary and sufficient
for $\mathcal{A}$ to respect $\mathcal{M}$---more precisely,\\
$$(<M_1>_{\mathcal{A}}=<M_2>_{\mathcal{A}})\;\;\Longleftrightarrow\;\;$$
$$\sum_{PE_j\in \mathcal{PARM}}C(E,j,k)PE_j(\mathcal{A})=D(E,j,k)\mbox{ for }
1\leq k\leq \nu(\mathcal{M})\;.$$
\end{Prop}
\smallskip

\noindent Note: We shall call the system (\ref{eqn:Assoc}) of
simultaneous linear equations over $\KK$, the \emph{equations
associated to} $\mathcal{M}$.
\smallskip

\noindent \underline{\bf {PROOF of Prop. \ref{prop:Assoc}:}}\\
The assumption that $\mathcal{M}$ is grammatical, implies (cf. Prop.\ref{Prop:in} )
that there exist natural numbers $m=in(\mathcal{M})$, $n=out(\mathcal{M})$ such that
\begin{equation}\label{eqn:gramm}
m=in(M_1)=in(M_2)\;,\;n=out(M_1)=out(M_2)\;.
\end{equation}
$<M_1>_{\mathcal{A}}$ and $<M_2>_{\mathcal{A}}$ may (by the preceding proposition)
be written in the form (\ref{eqn:coordMovieAmp}), and hence so can
$<M_1>_{}-<M_2>_{\mathcal{A}}$--- i.e. there exist
$$p,q,\epsilon_1,\cdots,\epsilon_p,\epsilon_1',\cdots,\epsilon_q',
S_1,\cdots,S_p,T_1,\cdots,T_q,(E^1,k^{(1)}),\cdots,(E^q,k^{(q)})$$
as in the preceding proposition, (depending on $\mathcal{M}$ but
not on $\mathcal{A}$) such that, for every strongly normal
amplitude-assignment $\mathcal{A}$ ,
$$<M_1>_{\mathcal{A}}-<M_2>_{\mathcal{A}}= \sum_{i=1}^p \epsilon_i<S_i>+
\sum_{j=1}^q \epsilon_j'(PE^j)_{k^{(j)}}(\mathcal{A}) <T_j>\;.$$
Hence, for every strongly normal amplitude-assignment
$\mathcal{A}$, $\mathcal{A}$ respects $\mathcal{M}$, if and only
if
\begin{equation}\label{eqn:resp1}
\sum_{i=1}^p \epsilon_i<S_i>+\sum_{j=1}^q \epsilon_j'(PE^j)_{k^{(j)}}(\mathcal{A}) <T_j>=0
\end{equation}
Now, the Kauffman amplitudes
$$<S_i>,<T_j> \eqno(*)$$
in (\ref{eqn:resp1}), are all $\KK$-linear maps from $V_1^{\otimes
m}$ to $V_1^{\otimes n}$. If we choose $\KK$-bases for
$V_1^{\otimes m}$ and $V_1^{\otimes n}$, the amplitudes in (*) are
all represented by $ 2^n \times 2^m$ matrices, so that the
equation (\ref{eqn:resp1}) between matrices, is equivalent (as
asserted) to a collection of $2^{m+n}$ simultaneous linear
equations over $\KK$ of the form (\ref{eqn:Assoc}).\\
In more detail:\\
Let $\mathcal{B}_m$ denote the set of $2^m$ ordered $m$-tuples of 0's and 1's:
$$\mathcal{B}_m =\{(\epsilon_0,\cdots,\epsilon_m):\mbox{ all }\epsilon's \in \{0,1,\}\}$$
$V_1^{\otimes m}$ has the $\KK$-basis (indexed by $\mathcal{B}_m$)
$$B_m=\{e_\alpha =e_{\epsilon_1}\otimes \cdots \otimes  e_{\epsilon_m}\; |
\; \alpha=(\epsilon_1,\cdots,\epsilon_m)\in \mathcal{B}_m \}  $$
and $V_1^{\otimes n}$ has the similar basis $B_n$. With respect to
this choice of bases for $V_1^{\otimes m}$ and $V_1^{\otimes n}$,
let us assume that $<S_i>,<T_j>$ are represented, respectively, by
the $2^n \times 2^m$ matrices
$$\mathcal{S}^{(i)}_{\beta, \alpha}, \mathcal{T}^{(j)}_{\beta,\alpha}\hspace{.4cm}
(\alpha\in\mathcal{B}_m,\;\beta \in\mathcal{B}_n)$$ over
$\KK$---i.e., that, for all $\alpha$ in $\mathcal{B}_m$,
$$<S_i>(e_{\alpha})=\sum_{\beta\in \mathcal{B}_n}\mathcal{S}^{(i)}_{\beta, \alpha}
e_{\beta}\;\;\; \mbox{   for } 1\leq i\leq p$$
and
$$<T_j>(e_{\alpha})=\sum_{\beta\in \mathcal{B}_n}\mathcal{T}^{(j)}_{\beta, \alpha}
e_{\beta}\;\;\; \mbox{   for } 1\leq j\leq q$$
Then the necessary and sufficient condition that $\mathcal{A}$ respect $\mathcal{M}$,
 i.e., that (\ref{eqn:resp1}) hold, is that, for all $\alpha$ in $\mathcal{B}_m $
 and all $\beta$ in $\mathcal{B}_m $,
\begin{equation}\label{eqn:resp2}
\sum_{i=1}^p \epsilon_i \mathcal{S}^{(i)}_{\beta, \alpha}+
\sum_{j=1}^q \epsilon_j'(PE^j)_{k^{(j)}}(\mathcal{A}) \mathcal{T}^{(j)}_{\beta,\alpha}=0
\end{equation}
\smallskip
Since each of the $2^{m+n}$ equations (\ref{eqn:resp2}) is
essentially in the format (\ref{eqn:Assoc})
the proof of Prop. \ref{prop:Assoc} is complete.
\smallskip

\noindent {\bf REMARKS:}\\
\noindent {\bf REMARK ONE:} For the purposes of the present paper,
the proof just presented should be  thought of, not so much as
demonstrating a  mathematical \underline{{\bf existence theorem}},
but rather as sketching an \underline{{\bf algorithm}} for
computing the system of equations associated to a movie-move
$\mathcal{M}$. This algorithm  is incorporated in 2KnotsLib.\\
\noindent {\bf REMARK TWO:} Let us note that the preceding proof
gives a {\bf specific} collection (\ref{eqn:Assoc}) of equations
associated to $\mathcal{M}$, which we shall denote by
$Assoc(\mathcal{M})$. The cardinality of $Assoc(\mathcal{M})$ is
$\nu(\mathcal{M})=2^{in(\mathcal{M})+out(\mathcal{M})}$.
\bigskip

We are finally ready to construct ---using our class libraries---a
driver program for computing $\MU$-balanced amplitude-assignments.
This program will only consider {\bf strongly normal}
amplitude-assignments; the present authors do not know if there
exist $\MU$-balanced amplitude-assignments which are not strongly
normal, and if there are, our present considerations seem to give
no information concerning them. By contrast, we shall explicitly
compute {\bf all} $\MU$-balanced strongly normal
amplitude-assignments.
\smallskip

We shall from now on assume that $\mathcal{A}$ is a strongly
normal amplitude-assignment. By combining Propositions
\ref{finiteTests}, \ref{prop:Assoc31} and \ref{prop:Assoc}, we
obtain an explicitly computable collection of linear equations
over $\KK$, for  the strongly normal coordinates of $\mathcal{A}$,
whose solutions give the set of all such $\mathcal{A}$ which are
$\MU$-balanced.
In more detail:\\
We divide into two cases, according as $31 \in \mathcal{U}$ or
not.\\
\underline{Consider first the case that $31 \in
\mathcal{U}$:}\\
Since strongly normal implies semi-normal, we may
apply (b) of Prop. \ref{finiteTests}, which says that
$\mathcal{A}$ is $\MU$-balanced, if and only if it respects each
of the finitely many movie-moves in the set
\begin{equation}\label{eqn:assoc1}
\mathcal{S}:=\bigcup_{i\in \MU^c}\;((\mathcal{SYM } \circ
\mathcal{C })\cap \MMM_i)
\end{equation}
 (Recall that $\mathcal{C}$ denotes the set of 40
BL-reduced movie-moves whose type lies in $B=\{i\;|\;1\leq i\leq\
14\} \cup \{21\} \cup \{i\; |\; 23\leq i\leq 30\}$.)
\smallskip

Let $\mathcal{M}=(M_1,M_2)$ be a movie-move in $\mathcal{S}$.
Movie-moves are grammatical, so we may apply Prop.
\ref{prop:Assoc} to $\mathcal{M}$: there is a finite collection
$Assoc(\mathcal{M})$, consisting of
$$\nu(\mathcal{M})=2^{in(\mathcal{M})+out(\mathcal{M})}$$
inhomogeneous linear equations in the 102 normal coordinates of
$\mathcal{A}$ with coefficients in $\KK$:
$$\sum_{PE_j\in \mathcal{PARM}}C(E,j,k)PE_j(\mathcal{A})=D(E,j,k)\mbox{ for }
1\leq k\leq \nu(\mathcal{M})\;.$$
which hold simultaneously, if and only if $\mathcal{A}$ respects $\mathcal{M}$.
\medskip

Let us define the set $Assoc(\MU)$ of \underline{{\bf
$\MU$-associated equations}}, to be the union of the sets of
associated equations to all movie-moves in (\ref{eqn:assoc1}).

Then, combining the two preceding observations, we see that $\mathcal{A}$
 is $\MU$-balanced, if and only if its 102 normal coordinates satisfy
each of the linear equations in the finite set $Assoc(\MU)$ of
associated equations---provided that $31\in \MU$.\\
 \underline{Consider next the case $31\notin \MU$:}\\
 The computation precisely is the same as that in the preceding case, except
 that we must replace $\mathcal{S}$ (given by (\ref{eqn:assoc1}))
 by
 $$\mathcal{S}_1:=(\bigcup_{i\in \MU^c\verb"\"\{31\}}\;((\mathcal{SYM } \circ
\mathcal{C })\cap \MMM_i))\bigcup \MMM_{31}^{red}
$$
and then define $Assoc(\MU)$ to be the union of the sets of
associated equations to all movie-moves in $\mathcal{S}_1$.

\subsection{Proof that Movie-move 31 is Essential}
As an application of the machinery developed above, we now sketch
a proof of the result, announced in the introduction, that the
movie-moves of type 31 are `essential', i.e. cannot be deduced
from the remaining movie-moves, in the sense of Def.
\ref{def:essentialMM}.
\smallskip

The proof to be presented here, involves explicitly solving the
system Assoc($\M{U})$ of equations associated to $\M{U}$, in the
two special cases $\M{U}=\{31\}$ and $\M{U}=\emptyset$.

 What is
the cardinality of the set $Assoc(\{31\})$? Both versions of the
2KnotsLib (augmented in both cases by the driver Program A)
furnish the same answer to this question, namely 12288. There are
not `really' that many equations to solve: we can reduce this
number
substantially, by the two simple expedients of:\\
a) removing wherever it occurs in $Assoc(\{31\})$, the equation $0=0$,\\
and,\\
b) if any equation occurs more than once in $Assoc(\{31\})$,
retaining only one copy (If one equation is a scalar multiple of a
second--- but not identical with it---they are to count as
distinct in this new tally.)
\medskip

Let us denote by $Assoc'(\{31\})$ the smaller (but equivalent)
collection of equations thus obtained. 2KnotsLib tells us the
cardinality of $Assoc'(\{31\})$ is 2856. The most general solution
to this set $Assoc'(\{31\})$ of 2856 inhomogeneous linear
equations over the field $K$ in 102 unknowns is furnished by both
versions of the driver Program B, with the following results
(which would seem difficult to establish without use of  the
computer, at least using only the methods of the present paper):

 The rank of this system of 2856 equations in
102 unknowns is 98. The family of all solutions represents an
affine subspace of $F^{102}$, of affine \footnote{Recall that, if
S is an affine subspace of $K^n$, then if $S$ is empty, it is said
to have affine rank $-1$, while if there exists $v\in S$, then
$\{v'-v | v'\in S \}$ is a vector subspace of $K^n$, whose
dimension over $K$ (which is independent of the choice of $v$ in
$S$) is called the affine rank of $S$.}rank $102-98=4$. (The
reader who wishes to see explicitly the resulting
$\{31\}$-balanced amplitude-assignments, is referred to \S 4.1
below)

 \noindent \underline{{\bf NOTE:}} Since there are many more
equations than unknowns, it seems like more good luck than one
should have a right to expect,
 that this method of constructing $\{31\}$-regular isotopy invariants for
2-knots works. Perhaps the credit for this minor miracle should go
to Def.~\ref{def:ampMovie}, i.e. the preceding paragraph furnishes
evidence for the suitability of formula (\ref{eqn:ampMovie}) for
the amplitude of a movie. In another direction: When equations
have more solutions than one might anticipate, one suspects the
explanation involves the action of some symmetry---here perhaps
the action of the quantum group $sl_q(2)$.
\medskip

More generally, for any proper subset $\mathcal{U}$, we may define
rk($\M{U}$) to be the affine rank of the affine subspace of
$K^{102}$, defined by the system of equations $Assoc(\M{U})$. The
preceding discussion shows that $rk(\{31\})=4$. Another
application of the two 2KnotsLib class libraries shows that
$rk(\emptyset)$=1 (cf \S 4.1 below, based on solving a reduced
system $Assoc'(\emptyset)$ of 10,208 equations in 102 unknowns).
\begin{Prop}\label{prop:essential}
Let $1\leq i\leq 31$; then rk$(\{i\})\geq 1$, and a sufficient
condition for $i$ to be essential, is that rk$(\{i\})> 1$.
\end{Prop}\label{31}
\smallskip

\noindent \underline{Sketch of proof:}Let $S,S_0$ be the affine
subspaces of $K^{102}$ defined respectively by $Assoc(\{i\})$ and
$Assoc(\emptyset)$. Since
$$Assoc(\{i\})\subseteq Assoc(\emptyset)$$
we have $S \supseteq S_0$ and hence rk$(\{i\})\geq 1$. \\
Now assume that $i$ is not essential; then every $\{i\}$-balanced
amplitude-assignment is  $\emptyset$-balanced, which implies
$S=S_0$, and hence rk($\{i\})=1$. Q.E.D.
\smallskip

\noindent \underline{{\bf Corollary:}}\emph{Movie-moves of type 31
are essential} (i.e. the Carter-Rieger-Saito Movie-move Theorem
would no longer be correct if they were omitted from the full list
of movie-moves.)
\smallskip

\noindent \underline{{\bf AN OPEN QUESTION}}\\
 Our program shows that not every
$i$ between 1 and 31 satisfies the condition
$$rk(\{i\})> 1$$ of Prop.~\ref{prop:essential}---for instance, for
type 6 (a collection of movie-moves related to the
Knizhnik-Zamolodchikov equations, as Reidemeister III is to the
Yang-Baxter equations)
both versions of the 2Knots class library show that rk($\{6\}$)=1.\\
This raises the question, which unfortunately we must leave open
at this time, whether the converse to Prop.~\ref{prop:essential}
holds:\\
\underline{{\bf CONJECTURE 1:}} Let $1\leq i \leq 31$; then
movie-moves of type $i$ are essential, if and only if
$$rk(\{i\})>1\;.$$

\section{Some Examples}\label{S:examples}
In \S3.9 we discussed the general solution to the set of equations
$Assoc(\M{U})$ for the two cases $\M{U} = \{31\}$ and $\MU =
\emptyset$ . In the present section, we continue  this discussion,
and compute the associated $\M{U}$-regular isotopy invariants for
seven specific examples of 2-knots---namely, an unknotted sphere,
Klein bottle and torus, two knotted spheres and the `1-twist and
2-twist spun trefoil' described on p.36 of [CS2].
\smallskip

These seven examples are given in the {\bf sf-notation} described
in section \ref{SS:notation}. This is the input format for both
the Java and C++ programs. However the two programs use slightly
different syntax for comments and end markers for input movies.
The Java program expects movies to be terminated by a '.' and
comment lines to begin with '$\#$'; the C++ program (whose
notation will be used in the present section)uses '$\%$' to begin
comment lines and '$\#$' to terminate each input movie. See the
class libraries' respective documentations for details.
\medskip

\subsection{$\mathcal{U}$-regular isotopy for $\M{U}=\{31\}$ and
for $\M{U}=\emptyset$}\label{SS:Basis}

\noindent
\underline{{\bf  The Case $\MU = \{31\}$}} \\
As we saw in \S3.9, the $\{31\}$-regular isotopy invariants are
furnished by the solutions to a set ( furnished by our driver
program A) of 2856 inhomogenous linear equations over the
ground-field $\KK= \QQ(q)$ in  102 unknowns.

This system of 2856 equations in the 102 normal coodinates has
affine rank 4 (cf. footnote (18).). The family of all solutions is
parametrized below, where $t_1, t_2, t_3,$ and $t_4$ can have
arbitrary values from $\KK=\QQ(q)$. (See the ReadMes accompanying
our two class libraries, for discussion of how this
parametrization is obtained.)

\begin{tabbing}
$P_{ETxxxR_1}$ \= $=$ \= $-t_1 +  1 - q^{-3}$ xx \= $P_{ETxxxR_1}$
\= $=$ \= $-t_1 +  1 - q^{-3}$ x \= $P_{ETxxxR_1}$ \= $=$ \= $-t_1
+  1 - q^{-3}$ \kill $P_{ET1I_1}  $ \> $=$ \> $  t_1 $ \>
$P_{ET3mbt_3}  $ \> $=$ \> $  q^{-1} $ \> $P_{ET3tmb_1}  $ \> $=$
\> $  q^{-3} $
 \\ [.05ex]
$P_{ET1R_1}  $ \> $=$ \> $  -t_1 +  1 - q^{-3} $ \> $P_{ET3mbt_4}
$ \> $=$ \> $  q^{-3} $ \> $P_{ET3tmb_2}  $ \> $=$ \> $  q^{-1} $
 \\ [.05ex]
$P_{ET1t_1}  $ \> $=$ \> $  (q^3)t_1 - q^3 +  1 $ \> $P_{ET3mbt_5}
$ \> $=$ \> $  q $ \> $P_{ET3tmb_3}  $ \> $=$ \> $  q^{-1} $
 \\ [.05ex]
$P_{ET1tR_1}  $ \> $=$ \> $  (-q^3)t_1 $ \> $P_{ET3mbtR_1}  $ \>
$=$ \> $  q $ \> $P_{ET3tmb_4}  $ \> $=$ \> $  q $
 \\ [.05ex]
$P_{ET1f_1}  $ \> $=$ \> $  (q^3)t_1 - q^3 +  1 $ \>
$P_{ET3mbtR_2}  $ \> $=$ \> $  q^{-1} $ \> $P_{ET3tmb_5}  $ \> $=$
\> $  q $
 \\ [.05ex]
$P_{ET1fR_1}  $ \> $=$ \> $  (-q^3)t_1 $ \> $P_{ET3mbtR_3}  $ \>
$=$ \> $  q^{-1} $ \> $P_{ET3tmbR_1}  $ \> $=$ \> $  q^{-3} $
 \\ [.05ex]
$P_{ET1ft_1}  $ \> $=$ \> $  t_1 $ \> $P_{ET3mbtR_4}  $ \> $=$ \>
$  q^{-3} $ \> $P_{ET3tmbR_2}  $ \> $=$ \> $  q^{-1} $
 \\ [.05ex]
$P_{ET1ftR_1}  $ \> $=$ \> $  -t_1 +  1 - q^{-3} $ \>
$P_{ET3mbtR_5}  $ \> $=$ \> $  q $ \> $P_{ET3tmbR_3}  $ \> $=$ \>
$  q^{-1} $
 \\ [.05ex]
$P_{ET2I_1}  $ \> $=$ \> $ 1$ \> $P_{ET3mtb_1}  $ \> $=$ \> $
q^{-1} $ \> $P_{ET3tmbR_4}  $ \> $=$ \> $  q $
 \\ [.05ex]
$P_{ET2I_2}  $ \> $=$ \> $ 0$ \> $P_{ET3mtb_2}  $ \> $=$ \> $  q $
\> $P_{ET3tmbR_5}  $ \> $=$ \> $  q $
 \\ [.05ex]
$P_{ET2R_1}  $ \> $=$ \> $ 1$ \> $P_{ET3mtb_3}  $ \> $=$ \> $  q $
\> $P_{ET4I_1}  $ \> $=$ \> $  q $
 \\ [.05ex]
$P_{ET2R_2}  $ \> $=$ \> $ 0$ \> $P_{ET3mtb_4}  $ \> $=$ \> $  q^3
$ \> $P_{ET4I_2}  $ \> $=$ \> $  q^{-1} $
 \\ [.05ex]
$P_{ET2f_1}  $ \> $=$ \> $ 1$ \> $P_{ET3mtb_5}  $ \> $=$ \> $
q^{-1} $ \> $P_{ET4R_1}  $ \> $=$ \> $  q $
 \\ [.05ex]
$P_{ET2f_2}  $ \> $=$ \> $ 0$ \> $P_{ET3mtbR_1}  $ \> $=$ \> $
q^{-1} $ \> $P_{ET4R_2}  $ \> $=$ \> $  q^{-1} $
 \\ [.05ex]
$P_{ET2fR_1}  $ \> $=$ \> $ 1$ \> $P_{ET3mtbR_2}  $ \> $=$ \> $  q
$ \> $P_{ET4t_1}  $ \> $=$ \> $  q^{-1} $
 \\ [.05ex]
$P_{ET2fR_2}  $ \> $=$ \> $ 0$ \> $P_{ET3mtbR_3}  $ \> $=$ \> $  q
$ \> $P_{ET4t_2}  $ \> $=$ \> $  q $
 \\ [.05ex]
$P_{ET3bmt_1}  $ \> $=$ \> $  q^3 $ \> $P_{ET3mtbR_4}  $ \> $=$ \>
$  q^3 $ \> $P_{ET4tR_1}  $ \> $=$ \> $  q^{-1} $
 \\ [.05ex]
$P_{ET3bmt_2}  $ \> $=$ \> $  q $ \> $P_{ET3mtbR_5}  $ \> $=$ \> $
q^{-1} $ \> $P_{ET4tR_2}  $ \> $=$ \> $  q $
 \\ [.05ex]
$P_{ET3bmt_3}  $ \> $=$ \> $  q $ \> $P_{ET3tbm_1}  $ \> $=$ \> $
q^{-1} $ \> $P_{ET4f_1}  $ \> $=$ \> $  q^{-1} $
 \\ [.05ex]
$P_{ET3bmt_4}  $ \> $=$ \> $  q^{-1} $ \> $P_{ET3tbm_2}  $ \> $=$
\> $  q $ \> $P_{ET4f_2}  $ \> $=$ \> $  q $
 \\ [.05ex]
$P_{ET3bmt_5}  $ \> $=$ \> $  q^{-1} $ \> $P_{ET3tbm_3}  $ \> $=$
\> $  q $ \> $P_{ET4fR_1}  $ \> $=$ \> $  q^{-1} $
 \\ [.05ex]
$P_{ET3bmtR_1}  $ \> $=$ \> $  q^3 $ \> $P_{ET3tbm_4}  $ \> $=$ \>
$  q^{-1} $ \> $P_{ET4fR_2}  $ \> $=$ \> $  q $
 \\ [.05ex]
$P_{ET3bmtR_2}  $ \> $=$ \> $  q $ \> $P_{ET3tbm_5}  $ \> $=$ \> $
q^3 $ \> $P_{ET4ft_1}  $ \> $=$ \> $  q $
 \\ [.05ex]
$P_{ET3bmtR_3}  $ \> $=$ \> $  q $ \> $P_{ET3tbmR_1}  $ \> $=$ \>
$  q^{-1} $ \> $P_{ET4ft_2}  $ \> $=$ \> $  q^{-1} $
 \\ [.05ex]
$P_{ET3bmtR_4}  $ \> $=$ \> $  q^{-1} $ \> $P_{ET3tbmR_2}  $ \>
$=$ \> $  q $ \> $P_{ET4ftR_1}  $ \> $=$ \> $  q $
 \\ [.05ex]
$P_{ET3bmtR_5}  $ \> $=$ \> $  q^{-1} $ \> $P_{ET3tbmR_3}  $ \>
$=$ \> $  q $ \> $P_{ET4ftR_2}  $ \> $=$ \> $  q^{-1} $
 \\ [.05ex]
$P_{ET3btm_1}  $ \> $=$ \> $  q $ \> $P_{ET3tbmR_4}  $ \> $=$ \> $
q^{-1} $ \>

 \\ [.05ex]
$P_{ET3btm_2}  $ \> $=$ \> $  q^{-1} $ \> $P_{ET3tbmR_5}  $ \> $=$
\> $  q^3 $ \>

 \\ [.05ex]
$P_{ET3btm_3}  $ \> $=$ \> $  q^{-1} $ \> $P_{ET5I_1}  $ \> $=$ \>
$  -t_4 + (q^2 + q^{-2})(t_3 - 1)  - t_2  +  2 $
 \\ [.05ex]
$P_{ET3btm_4}  $ \> $=$ \> $  q $ \> $P_{ET5R_1}  $ \> $=$ \> $
t_4 + (-q^2 - q^{-2})(t_3 - 1) + t_2 $
 \\ [.05ex]
$P_{ET3btm_5}  $ \> $=$ \> $  q^{-3} $ \> $P_{ET5t_1}  $ \> $=$ \>
$  -t_4 + (q^2 + q^{-2})(t_3 - 1) - t_2  + 2 $
 \\ [.05ex]
$P_{ET3btmR_1}  $ \> $=$ \> $  q $ \> $P_{ET5tR_1}  $ \> $=$ \> $
t_4 + (-q^2 - q^{-2})(t_3 - 1) + t_2 $
 \\ [.05ex]
$P_{ET3btmR_2}  $ \> $=$ \> $  q^{-1} $ \> $P_{ET6I_1}  $ \> $=$
\> $  (q^2  - 1 + q^{-2})(t_4 + t_3 - 2) - t_2 $
 \\ [.05ex]
$P_{ET3btmR_3}  $ \> $=$ \> $  q^{-1} $ \> $P_{ET6R_1}  $ \> $=$
\> $  t_2 $
 \\ [.05ex]
$P_{ET3btmR_4}  $ \> $=$ \> $  q $ \> $P_{ET7I_1}  $ \> $=$ \> $
(q^2 + q^{-2})(t_3 - 1) - t_2 +  1 $
 \\ [.05ex]
$P_{ET3btmR_5}  $ \> $=$ \> $  q^{-3} $ \> $P_{ET7I_2}  $ \> $=$
\> $  t_3 $
 \\ [.05ex]
$P_{ET3mbt_1}  $ \> $=$ \> $  q $ \> $P_{ET7R_1}  $ \> $=$ \> $
t_4 + (-q^2 +  1 - q^{-2})(t_3 - 1) + t_2 $
 \\ [.05ex]
$P_{ET3mbt_2}  $ \> $=$ \> $  q^{-1} $ \> $P_{ET7R_2}  $ \> $=$ \>
$  t_4 $
 \\ [.05ex]
\end{tabbing}
We will denote by $\mathcal{A}_0$ the $\{31\}$-balanced normal
amplitude assignment obtained from our solution by setting each
$t_i = 0$. In addition, for each $i = 1, 2, 3, 4$, we  define a
$\{31\}$-balanced normal amplitude assignment $\mathcal{A}_i$ to
be the normal coordinate assignment obtained from  our solutions,
where $t_i = 1$
and $t_j = 0$ for each $j \ne i$. \\

For a movie $M$,
$${<}M{>}_{i}:={<}M{>}_{\mathcal{A}_i}\;\;\;\;(0\leq i\leq 4)$$

\noindent will denote the amplitudes of $M$ for the
$\{31\}$-balanced amplitude assignments $\mathcal{A}_{i}$.
\bigskip

\noindent \underline{{\bf The Case $\MU = \emptyset$}}: \\

Consideration of the set $\M{M}^{red}_{31}$ alone generates 3194
$\emptyset$-associated equations after removing zero equations and
repetitions.

With the extra associated equations provided by movie move 31, the
solutions to the $\emptyset$-associated equations in the 102
normal coordinates affine rank 1. The general solution in this
case is given below where $t$ can have arbitrary value in the
ground-field, $\KK$.
\bigskip

\begin{tabbing} $P_{ETxxxR_1}$ \= $=$ \= $-t_1 +  1 - q^{-3}$ xxxx
\= $P_{ETxxxR_1}$ \= $=$ \= $-t_1 +  1 - q^{-3}$ \= $P_{ETxxxR_1}$
\= $=$ \= $-t_1 +  1 - q^{-3}$ \kill $P_{ET1I_1}  $ \> $=$ \> $  (
1 + q^{-3})t - q^{-3} $ \> $P_{ET3mbt_3}  $ \> $=$ \> $  q^{-1} $
\> $P_{ET3tmb_1}  $ \> $=$ \> $  q^{-3} $
 \\ [.05ex]
$P_{ET1R_1}  $ \> $=$ \> $  ( -1 - q^{-3})t +  1 $ \>
$P_{ET3mbt_4}  $ \> $=$ \> $  q^{-3} $ \> $P_{ET3tmb_2}  $ \> $=$
\> $  q^{-1} $
 \\ [.05ex]
$P_{ET1t_1}  $ \> $=$ \> $  (q^3 +  1)t - q^3 $ \> $P_{ET3mbt_5}
$ \> $=$ \> $  q $ \> $P_{ET3tmb_3}  $ \> $=$ \> $  q^{-1} $
 \\ [.05ex]
$P_{ET1tR_1}  $ \> $=$ \> $  (-q^3  - 1)t +  1 $ \> $P_{ET3mbtR_1}
$ \> $=$ \> $  q $ \> $P_{ET3tmb_4}  $ \> $=$ \> $  q $
 \\ [.05ex]
$P_{ET1f_1}  $ \> $=$ \> $  (q^3 +  1)t - q^3 $ \> $P_{ET3mbtR_2}
$ \> $=$ \> $  q^{-1} $ \> $P_{ET3tmb_5}  $ \> $=$ \> $  q $
 \\ [.05ex]
$P_{ET1fR_1}  $ \> $=$ \> $  (-q^3  - 1)t +  1 $ \> $P_{ET3mbtR_3}
$ \> $=$ \> $  q^{-1} $ \> $P_{ET3tmbR_1}  $ \> $=$ \> $  q^{-3} $
 \\ [.05ex]
$P_{ET1ft_1}  $ \> $=$ \> $  ( 1 + q^{-3})t - q^{-3} $ \>
$P_{ET3mbtR_4}  $ \> $=$ \> $  q^{-3} $ \> $P_{ET3tmbR_2}  $ \>
$=$ \> $  q^{-1} $
 \\ [.05ex]
$P_{ET1ftR_1}  $ \> $=$ \> $  ( -1 - q^{-3})t +  1 $ \>
$P_{ET3mbtR_5}  $ \> $=$ \> $  q $ \> $P_{ET3tmbR_3}  $ \> $=$ \>
$  q^{-1} $
 \\ [.05ex]
$P_{ET2I_1}  $ \> $=$ \> $ 1$ \> $P_{ET3mtb_1}  $ \> $=$ \> $
q^{-1} $ \> $P_{ET3tmbR_4}  $ \> $=$ \> $  q $
 \\ [.05ex]
$P_{ET2I_2}  $ \> $=$ \> $ 0$ \> $P_{ET3mtb_2}  $ \> $=$ \> $  q $
\> $P_{ET3tmbR_5}  $ \> $=$ \> $  q $
 \\ [.05ex]
$P_{ET2R_1}  $ \> $=$ \> $ 1$ \> $P_{ET3mtb_3}  $ \> $=$ \> $  q $
\> $P_{ET4I_1}  $ \> $=$ \> $  q $
 \\ [.05ex]
$P_{ET2R_2}  $ \> $=$ \> $ 0$ \> $P_{ET3mtb_4}  $ \> $=$ \> $  q^3
$ \> $P_{ET4I_2}  $ \> $=$ \> $  q^{-1} $
 \\ [.05ex]
$P_{ET2f_1}  $ \> $=$ \> $ 1$ \> $P_{ET3mtb_5}  $ \> $=$ \> $
q^{-1} $ \> $P_{ET4R_1}  $ \> $=$ \> $  q $
 \\ [.05ex]
$P_{ET2f_2}  $ \> $=$ \> $ 0$ \> $P_{ET3mtbR_1}  $ \> $=$ \> $
q^{-1} $ \> $P_{ET4R_2}  $ \> $=$ \> $  q^{-1} $
 \\ [.05ex]
$P_{ET2fR_1}  $ \> $=$ \> $ 1$ \> $P_{ET3mtbR_2}  $ \> $=$ \> $  q
$ \> $P_{ET4t_1}  $ \> $=$ \> $  q^{-1} $
 \\ [.05ex]
$P_{ET2fR_2}  $ \> $=$ \> $ 0$ \> $P_{ET3mtbR_3}  $ \> $=$ \> $  q
$ \> $P_{ET4t_2}  $ \> $=$ \> $  q $
 \\ [.05ex]
$P_{ET3bmt_1}  $ \> $=$ \> $  q^3 $ \> $P_{ET3mtbR_4}  $ \> $=$ \>
$  q^3 $ \> $P_{ET4tR_1}  $ \> $=$ \> $  q^{-1} $
 \\ [.05ex]
$P_{ET3bmt_2}  $ \> $=$ \> $  q $ \> $P_{ET3mtbR_5}  $ \> $=$ \> $
q^{-1} $ \> $P_{ET4tR_2}  $ \> $=$ \> $  q $
 \\ [.05ex]
$P_{ET3bmt_3}  $ \> $=$ \> $  q $ \> $P_{ET3tbm_1}  $ \> $=$ \> $
q^{-1} $ \> $P_{ET4f_1}  $ \> $=$ \> $  q^{-1} $
 \\ [.05ex]
$P_{ET3bmt_4}  $ \> $=$ \> $  q^{-1} $ \> $P_{ET3tbm_2}  $ \> $=$
\> $  q $ \> $P_{ET4f_2}  $ \> $=$ \> $  q $
 \\ [.05ex]
$P_{ET3bmt_5}  $ \> $=$ \> $  q^{-1} $ \> $P_{ET3tbm_3}  $ \> $=$
\> $  q $ \> $P_{ET4fR_1}  $ \> $=$ \> $  q^{-1} $
 \\ [.05ex]
$P_{ET3bmtR_1}  $ \> $=$ \> $  q^3 $ \> $P_{ET3tbm_4}  $ \> $=$ \>
$  q^{-1} $ \> $P_{ET4fR_2}  $ \> $=$ \> $  q $
 \\ [.05ex]
$P_{ET3bmtR_2}  $ \> $=$ \> $  q $ \> $P_{ET3tbm_5}  $ \> $=$ \> $
q^3 $ \> $P_{ET4ft_1}  $ \> $=$ \> $  q $
 \\ [.05ex]
$P_{ET3bmtR_3}  $ \> $=$ \> $  q $ \> $P_{ET3tbmR_1}  $ \> $=$ \>
$  q^{-1} $ \> $P_{ET4ft_2}  $ \> $=$ \> $  q^{-1} $
 \\ [.05ex]
$P_{ET3bmtR_4}  $ \> $=$ \> $  q^{-1} $ \> $P_{ET3tbmR_2}  $ \>
$=$ \> $  q $ \> $P_{ET4ftR_1}  $ \> $=$ \> $  q $
 \\ [.05ex]
$P_{ET3bmtR_5}  $ \> $=$ \> $  q^{-1} $ \> $P_{ET3tbmR_3}  $ \>
$=$ \> $  q $ \> $P_{ET4ftR_2}  $ \> $=$ \> $  q^{-1} $
 \\ [.05ex]
$P_{ET3btm_1}  $ \> $=$ \> $  q $ \> $P_{ET3tbmR_4}  $ \> $=$ \> $
q^{-1} $ \>

 \\ [.05ex]
$P_{ET3btm_2}  $ \> $=$ \> $  q^{-1} $ \> $P_{ET3tbmR_5}  $ \> $=$
\> $  q^3 $ \>

 \\ [.05ex]
$P_{ET3btm_3}  $ \> $=$ \> $  q^{-1} $ \> $P_{ET5I_1}  $ \> $=$ \>
$ 1$
 \\ [.05ex]
$P_{ET3btm_4}  $ \> $=$ \> $  q $ \> $P_{ET5R_1}  $ \> $=$ \> $ 1$
 \\ [.05ex]
$P_{ET3btm_5}  $ \> $=$ \> $  q^{-3} $ \> $P_{ET5t_1}  $ \> $=$ \>
$ 1$
 \\ [.05ex]
$P_{ET3btmR_1}  $ \> $=$ \> $  q $ \> $P_{ET5tR_1}  $ \> $=$ \> $
1$
 \\ [.05ex]
$P_{ET3btmR_2}  $ \> $=$ \> $  q^{-1} $ \> $P_{ET6I_1}  $ \> $=$
\> $  (q^2 +  1 + q^{-2})t - q^2 - q^{-2} $
 \\ [.05ex]
$P_{ET3btmR_3}  $ \> $=$ \> $  q^{-1} $ \> $P_{ET6R_1}  $ \> $=$
\> $  (-q^2  - 1 - q^{-2})t +  1 $
 \\ [.05ex]
$P_{ET3btmR_4}  $ \> $=$ \> $  q $ \> $P_{ET7I_1}  $ \> $=$ \> $
t $
 \\ [.05ex]
$P_{ET3btmR_5}  $ \> $=$ \> $  q^{-3} $ \> $P_{ET7I_2}  $ \> $=$
\> $  -t +  1 $
 \\ [.05ex]
$P_{ET3mbt_1}  $ \> $=$ \> $  q $ \> $P_{ET7R_1}  $ \> $=$ \> $
-t +  1 $
 \\ [.05ex]
$P_{ET3mbt_2}  $ \> $=$ \> $  q^{-1} $ \> $P_{ET7R_2}  $ \> $=$ \>
$  t $
 \\ [.05ex]
\end{tabbing}

We will denote by $\mathcal{A}_a$ the $\emptyset$-balanced normal
amplitude assignment obtained from the preceding solution by
setting $t = 0$, and by $\M{A}_b$  the  $\emptyset$-balanced
normal amplitude assignment obtained by setting $t = 1$. For a
movie $M$, we then set
$${<}M{>}_{a}:={<}M{>}_{\M{A}_a}$$
and
$${<}M{>}_{b}:={<}M{>}_{\M{A}_b}$$
\smallskip

\noindent \underline{Note that, for any 2-tangle $K$,
${<}M{>}_{a}$
has the same value for all movies $M$}\\
\underline{ representing K, --- and the same is true for
${<}\;{>}_{b}$.} (This is one of the few places in the present
paper, where use is made of Th. \ref{th:main}.) However, the
computations below provide evidence that the resulting isotopy
invariants ${<}\;{>}_{a}$ and ${<}\;{>}_{b}$ are of a rather
trivial nature ---that is, these
computations support the following conjecture:\\

\noindent \underline{{\bf Conjecture 2:}} For every 2-knot $K$,
and every CRS movie $M$ representing $K$,
$${<}M{>}_{a}={<}M{>}_{b}=\mbox{ the number of connected components
of }K\;.$$
\subsection{Example: An Unknotted Sphere}
\par\noindent
The first 2-knot is a simple unknotted 2-sphere.
\begin{verbatim}
% UnknottedSphere
ss => sf[0,Cup,0][0,Cap,0]sf => ff #
\end{verbatim}
with amplitudes:

Amplitudes for ($\MU = \{31\}$):
\begin{eqnarray*}
{<}\mbox{UnknottedSphere}{>}_{0} & = & q^{-2}+q^{2} \\
{<}\mbox{UnknottedSphere}{>}_{1} & = & q^{-2}+q^{2} \\
{<}\mbox{UnknottedSphere}{>}_{2} & = & q^{-2}+q^{2} \\
{<}\mbox{UnknottedSphere}{>}_{3} & = & 1 \\
{<}\mbox{UnknottedSphere}{>}_{4} & = & 1
\end{eqnarray*}

Amplitudes for ($\MU = \emptyset$):
\begin{eqnarray*}
{<}\mbox{UnknottedSphere}{>}_{a} & = & 1 \\
{<}\mbox{UnknottedSphere}{>}_{b} & = & 1
\end{eqnarray*}

\subsection{Example: A Knotted Sphere}
\par\noindent
The next 2-knot is a knotted sphere communicated to us by Carter.
\begin{verbatim}
% KnottedSphereA
ss => ssf[0,Cup,0][0,Cap,0]f =>
f[0,Cup,0][0,Cap,0]f[0,Cup,0]ss[0,Cap,0] =>
[0,Cup,0][0,Cap,0][0,Cup,0]f[1,Cup,1]ss[0,Cap,2]f[0,Cap,0] =>
[0,Cup,0][0,Cap,0][0,Cup,0][1,Cup,1]ssf[1,Cap,1][1,Cup,1]f
   [0,Cap,2][0,Cap,0]
=> [0,Cup,0][0,Cap,0][0,Cup,0][1,Cup,1]f[2,NW,0][2,NE,0]f
   [1,Cap,1]s[1,Cup,1][0,Cap,2]s[0,Cap,0]
=> [0,Cup,0][0,Cap,0][0,Cup,0][1,Cup,1][2,NW,0]s[2,NE,0]
   [1,Cap,1]sff[0,Cap,0]
=> [0,Cup,0][0,Cap,0][0,Cup,0]s[1,Cup,1][2,NW,0]sf[1,NW,1]
   [2,Cap,0]f[0,Cap,0]
=> [0,Cup,0][0,Cap,0][0,Cup,0]f[2,Cup,0]s[1,NE,1]f[1,NW,1]s
   [2,Cap,0][0,Cap,0]
=> [0,Cup,0][0,Cap,0]s[0,Cup,0][2,Cup,0]sff[2,Cap,0][0,Cap,0] =>
[0,Cup,0][0,Cap,0]f[0,Cup,0]s[0,Cup,2]f[2,Cap,0]s[0,Cap,0] =>
[0,Cup,0]s[0,Cap,0][0,Cup,0]sf[0,Cap,0][0,Cup,0]f[0,Cap,0] =>
[0,Cup,0]f[2,Cup,0]s[0,Cap,2]f[0,Cap,0]s[0,Cup,0][0,Cap,0] =>
[0,Cup,0][2,Cup,0]ssf[2,Cap,0][0,Cap,0]f[0,Cup,0][0,Cap,0] =>
[0,Cup,0]s[2,Cup,0]f[1,NE,1]s[1,NW,1]f[2,Cap,0][0,Cap,0]
   [0,Cup,0][0,Cap,0]
=> [0,Cup,0]f[1,Cup,1][2,NW,0]sf[1,NW,1][2,Cap,0]s[0,Cap,0]
   [0,Cup,0][0,Cap,0]
=> [0,Cup,0][1,Cup,1]s[2,NW,0]f[2,NE,0]s[1,Cap,1]f[0,Cap,0]
   [0,Cup,0][0,Cap,0]
=> [0,Cup,0][1,Cup,1]ssff[1,Cap,1][0,Cap,0][0,Cup,0][0,Cap,0] =>
[0,Cup,0][1,Cup,1]f[2,NE,0]s[2,NW,0]f[1,Cap,1]s[0,Cap,0]
   [0,Cup,0][0,Cap,0]
=> [0,Cup,0]s[1,Cup,1][2,NE,0]sf[1,NE,1][2,Cap,0]f[0,Cap,0]
   [0,Cup,0][0,Cap,0]
=> [0,Cup,0]f[2,Cup,0]s[1,NW,1]f[1,NE,1]s[2,Cap,0][0,Cap,0]
   [0,Cup,0][0,Cap,0]
=> [0,Cup,0][2,Cup,0]sff[2,Cap,0][0,Cap,0]s[0,Cup,0][0,Cap,0] =>
[0,Cup,0]s[2,Cup,0]f[0,Cap,2]s[0,Cap,0]f[0,Cup,0][0,Cap,0] =>
[0,Cup,0]f[0,Cap,0][0,Cup,0]sf[0,Cap,0][0,Cup,0]s[0,Cap,0] =>
[0,Cup,0][0,Cap,0]s[0,Cup,0]f[0,Cup,2]s[2,Cap,0]f[0,Cap,0] =>
[0,Cup,0][0,Cap,0]f[0,Cup,0][2,Cup,0]ssf[2,Cap,0][0,Cap,0] =>
[0,Cup,0][0,Cap,0][0,Cup,0]s[2,Cup,0]f[1,NE,1]s[1,NW,1]f
   [2,Cap,0][0,Cap,0]
=> [0,Cup,0][0,Cap,0][0,Cup,0]f[1,Cup,1][2,NW,0]sf[1,NW,1]
   [2,Cap,0]s[0,Cap,0]
=> [0,Cup,0][0,Cap,0][0,Cup,0][1,Cup,1]s[2,NW,0]f[2,NE,0]s
   [1,Cap,1]f[0,Cap,0]
=> [0,Cup,0][0,Cap,0][0,Cup,0][1,Cup,1]ssff[1,Cap,1][0,Cap,0] =>
[0,Cup,0][0,Cap,0][0,Cup,0][1,Cup,1]f[2,NE,0]s[2,NW,0]f
   [1,Cap,1]s[0,Cap,0]
=> [0,Cup,0][0,Cap,0][0,Cup,0]s[1,Cup,1][2,NE,0]sf[1,NE,1]
   [2,Cap,0]f[0,Cap,0]
=> [0,Cup,0][0,Cap,0][0,Cup,0]f[2,Cup,0]s[1,NW,1]f[1,NE,1]s
   [2,Cap,0][0,Cap,0]
=> [0,Cup,0][0,Cap,0]s[0,Cup,0][2,Cup,0]sff[2,Cap,0][0,Cap,0] =>
[0,Cup,0][0,Cap,0]f[0,Cup,0]s[0,Cup,2]f[2,Cap,0]s[0,Cap,0] =>
[0,Cup,0]s[0,Cap,0][0,Cup,0]sf[0,Cap,0][0,Cup,0]f[0,Cap,0] =>
[0,Cup,0]f[2,Cup,0]s[0,Cap,2]f[0,Cap,0]s[0,Cup,0][0,Cap,0] =>
[0,Cup,0][2,Cup,0]ssf[2,Cap,0][0,Cap,0]f[0,Cup,0][0,Cap,0] =>
[0,Cup,0]s[2,Cup,0]f[1,NW,1]s[1,NE,1]f[2,Cap,0][0,Cap,0]
   [0,Cup,0][0,Cap,0]
=> [0,Cup,0]f[1,Cup,1][2,NE,0]sf[1,NE,1][2,Cap,0]s[0,Cap,0]
   [0,Cup,0][0,Cap,0]
=> [0,Cup,0][1,Cup,1]s[2,NE,0]f[2,NW,0]s[1,Cap,1]f[0,Cap,0]
   [0,Cup,0][0,Cap,0]
=> [0,Cup,0][1,Cup,1]ssff[1,Cap,1][0,Cap,0][0,Cup,0][0,Cap,0] =>
[0,Cup,0][1,Cup,1]f[2,NW,0]s[2,NE,0]f[1,Cap,1]s[0,Cap,0]
   [0,Cup,0][0,Cap,0]
=> [0,Cup,0]s[1,Cup,1][2,NW,0]sf[1,NW,1][2,Cap,0]f[0,Cap,0]
   [0,Cup,0][0,Cap,0]
=> [0,Cup,0]f[2,Cup,0]s[1,NE,1]f[1,NW,1]s[2,Cap,0][0,Cap,0]
   [0,Cup,0][0,Cap,0]
=> s[0,Cup,0][2,Cup,0]sff[2,Cap,0][0,Cap,0][0,Cup,0][0,Cap,0] =>
f[0,Cup,0]s[0,Cup,2]f[2,Cap,0]s[0,Cap,0][0,Cup,0][0,Cap,0] =>
[0,Cup,0]sf[0,Cap,0][0,Cup,0]sf[0,Cap,0][0,Cup,0][0,Cap,0] =>
[0,Cup,0]ff[0,Cap,0]s[0,Cup,0][0,Cap,0]s => s[0,Cup,0][0,Cap,0]sff
=> ff #
\end{verbatim}
with amplitudes:

Amplitudes for ($\MU = \{31\}$):
\begin{eqnarray*}
{<}\mbox{KnottedSphereA}{>}_{0} & = & -q^{-8}-q^{-6}-3q^{-4}-q^{-2}-4-q^{2}-3q^{4}-q^{6}-q^{8} \\
{<}\mbox{KnottedSphereA}{>}_{1} & = & -q^{-8}-q^{-6}-3q^{-4}-q^{-2}-4-q^{2}-3q^{4}-q^{6}-q^{8} \\
{<}\mbox{KnottedSphereA}{>}_{2} & = & -q^{-8}-2q^{-6}-4q^{-4}-4q^{-2}-6-4q^{2}-4q^{4}-2q^{6}-q^{8} \\
{<}\mbox{KnottedSphereA}{>}_{3} & = & q^{-6}+q^{-4}+3q^{-2}+3+3q^{2}+q^{4}+q^{6} \\
{<}\mbox{KnottedSphereA}{>}_{4} & = &
-q^{-8}-q^{-6}-4q^{-4}-3q^{-2}-5-3q^{2}-4q^{4}-q^{6}-q^{8}
\end{eqnarray*}

Amplitudes for ($\MU = \emptyset$):
\begin{eqnarray*}
{<}\mbox{KnottedSphereA}{>}_{a} & = & 1 \\
{<}\mbox{KnottedSphereA}{>}_{b} & = & 1
\end{eqnarray*}

\subsection{Example: Another knotted sphere}
\par\noindent
The next 2-knot is a knotted sphere from (Colin C. Adams, The Knot
Book, Fig.10.13):
\begin{verbatim}
% KnottedSphereB
ss => f[0,Cup,0]ss[0,Cap,0]f =>
[0,Cup,0]f[0,Cup,2]ss[0,Cap,2]f[0,Cap,0] =>
[0,Cup,0][0,Cup,2]f[1,NE,1]ss[1,NW,1]f[0,Cap,2][0,Cap,0] =>
[0,Cup,0][0,Cup,2][1,NE,1]f[0,NW,2]ss[0,NE,2]f[1,NW,1]
   [0,Cap,2][0,Cap,0]
=> [0,Cup,0][0,Cup,2][1,NE,1][0,NW,2]f[2,NW,0]ss[2,NE,0]f
   [0,NE,2][1,NW,1][0,Cap,2][0,Cap,0]
=> [0,Cup,0][0,Cup,2][1,NE,1][0,NW,2][2,NW,0]f[2,Cup,2]
   [1,Cap,3]ssf[2,NE,0][0,NE,2][1,NW,1][0,Cap,2][0,Cap,0]
=> [0,Cup,0][0,Cup,2][1,NE,1][0,NW,2][2,NW,0][2,Cup,2][1,Cap,3]
   ssf[1,Cup,3][2,Cap,2]f[2,NE,0][0,NE,2][1,NW,1][0,Cap,2][0,Cap,0]
=> [0,Cup,0][0,Cup,2][1,NE,1][0,NW,2][2,NW,0][2,Cup,2]
   s[1,Cap,3]f[1,Cap,1]s[1,Cup,1]f[1,Cup,3][2,Cap,2]
   [2,NE,0][0,NE,2][1,NW,1][0,Cap,2][0,Cap,0]
=> [0,Cup,0][0,Cup,2][1,NE,1][0,NW,2][2,NW,0]s[2,Cup,2]f
   [3,Cap,1]s[1,Cap,1]f[1,Cup,1][1,Cup,3][2,Cap,2][2,NE,0]
   [0,NE,2][1,NW,1][0,Cap,2][0,Cap,0]
=> [0,Cup,0][0,Cup,2][1,NE,1][0,NW,2][2,NW,0]ff[1,Cap,1]
   s[1,Cup,1][1,Cup,3]s[2,Cap,2][2,NE,0][0,NE,2][1,NW,1]
   [0,Cap,2][0,Cap,0]
=> [0,Cup,0][0,Cup,2][1,NE,1][0,NW,2][2,NW,0][1,Cap,1]
   f[1,Cup,1]s[3,Cup,1]f[2,Cap,2]s[2,NE,0][0,NE,2][1,NW,1]
   [0,Cap,2][0,Cap,0]
=> [0,Cup,0][0,Cup,2][1,NE,1][0,NW,2][2,NW,0]s[1,Cap,1][1,Cup,1]s
   ff[2,NE,0][0,NE,2][1,NW,1][0,Cap,2][0,Cap,0]
=> [0,Cup,0][0,Cup,2][1,NE,1][0,NW,2]s[2,NW,0]ff[2,NE,0]s
   [0,NE,2][1,NW,1][0,Cap,2][0,Cap,0]
=> [0,Cup,0][0,Cup,2][1,NE,1]s[0,NW,2]ff[0,NE,2]s[1,NW,1]
   [0,Cap,2][0,Cap,0]
=> [0,Cup,0][0,Cup,2]s[1,NE,1]ff[1,NW,1]s[0,Cap,2][0,Cap,0] =>
[0,Cup,0]s[0,Cup,2]ff[0,Cap,2]s[0,Cap,0] => s[0,Cup,0]ff[0,Cap,0]s
=> ff #
\end{verbatim}
with amplitudes:

Amplitudes for ($\MU = \{31\}$):
\begin{eqnarray*}
{<}\mbox{KnottedSphereB}{>}_{0} & = & 2q^{-16}-q^{-14}+2q^{-12}-2q^{-6}-2q^{-4}-2q^{-2}-4 \\
& & \quad{}-2q^{2}-2q^{4}-2q^{6}+2q^{12}-q^{14}+2q^{16} \\
{<}\mbox{KnottedSphereB}{>}_{1} & = & 2q^{-16}-q^{-14}+2q^{-12}-2q^{-6}-2q^{-4}-2q^{-2}-4 \\
& & \quad{}-2q^{2}-2q^{4}-2q^{6}+2q^{12}-q^{14}+2q^{16} \\
{<}\mbox{KnottedSphereB}{>}_{2} & = & 2q^{-16}+q^{-14}+2q^{-12}-2q^{-6}-4q^{-4}-6q^{-2}-8 \\
& & \quad{}-6q^{2}-4q^{4}-2q^{6}+2q^{12}+q^{14}+2q^{16} \\
{<}\mbox{KnottedSphereB}{>}_{3} & = & -2q^{-14}+2q^{-4}+4q^{-2}+5+4q^{2}+2q^{4}-2q^{14} \\
{<}\mbox{KnottedSphereB}{>}_{4} & = & 2q^{-16}+2q^{-12}-2q^{-6}-4q^{-4}-6q^{-2}-7 \\
& & \quad{}-6q^{2}-4q^{4}-2q^{6}+2q^{12}+2q^{16}
\end{eqnarray*}

Amplitudes for ($\MU = \emptyset$):
\begin{eqnarray*}
{<}\mbox{KnottedSphereB}{>}_{a} & = & 1 \\
{<}\mbox{KnottedSphereB}{>}_{b} & = & 1
\end{eqnarray*}

\subsection{Example: Klein Bottle}
\par\noindent
The next 2-knot is a Klein bottle
\begin{verbatim}
% Klein bottle
ss => f[0,Cup,0]ss[0,Cap,0]f =>
[0,Cup,0]ssf[0,Cup,2][1,Cap,1]f[0,Cap,0] =>
[0,Cup,0]f[1,Cup,1]s[0,Cap,2]f[0,Cup,2]s[1,Cap,1][0,Cap,0] =>
[0,Cup,0][1,Cup,1]sfsf[1,Cap,1][0,Cap,0] =>
[0,Cup,0]s[1,Cup,1]f[2,NE,0]s[2,NW,0]f[1,Cap,1][0,Cap,0] =>
[0,Cup,0]f[2,Cup,0][1,NW,1]fs[2,NW,0][1,Cap,1]s[0,Cap,0] =>
[0,Cup,0][2,Cup,0]s[1,NW,1]f[1,NE,1]s[2,Cap,0]f[0,Cap,0] =>
[0,Cup,0][2,Cup,0]sfsf[2,Cap,0][0,Cap,0] =>
[0,Cup,0]s[2,Cup,0]f[1,Cap,1]s[1,Cup,1]f[2,Cap,0][0,Cap,0] =>
[0,Cup,0]ffs[1,Cup,1][2,Cap,0]s[0,Cap,0] => s[0,Cup,0]ff[0,Cap,0]s
=> ff #
\end{verbatim}
with amplitudes:

Amplitudes for ($\MU = \{31\}$):
\begin{eqnarray*}
{<}\mbox{KleinBottle}{>}_{0} & = & q^{-4}+2+q^{4} \\
{<}\mbox{KleinBottle}{>}_{1} & = & q^{-4}+2+q^{4} \\
{<}\mbox{KleinBottle}{>}_{2} & = & q^{-4}+2+q^{4} \\
{<}\mbox{KleinBottle}{>}_{3} & = & 1 \\
{<}\mbox{KleinBottle}{>}_{4} & = & 1
\end{eqnarray*}

Amplitudes for ($\MU = \emptyset$):
\begin{eqnarray*}
{<}\mbox{KleinBottle}{>}_{a} & = & 1 \\
{<}\mbox{KleinBottle}{>}_{b} & = & 1
\end{eqnarray*}

\subsection{Example: Simple Torus}
\par\noindent
The next 2-knot is a simple (unknotted) torus
\begin{verbatim}
% SimpleTorus
ss => f[0,Cup,0]ss[0,Cap,0]f =>
[0,Cup,0]fs[0,Cap,0][0,Cup,0]fs[0,Cap,0] => s[0,Cup,0]ff[0,Cap,0]s
=> ff #
\end{verbatim}
with amplitudes:

Amplitudes for ($\MU = \{31\}$):
\begin{eqnarray*}
{<}\mbox{SimpleTorus}{>}_{0} & = & q^{-4}+2+q^{4} \\
{<}\mbox{SimpleTorus}{>}_{1} & = & q^{-4}+2+q^{4} \\
{<}\mbox{SimpleTorus}{>}_{2} & = & q^{-4}+2+q^{4} \\
{<}\mbox{SimpleTorus}{>}_{3} & = & 1 \\
{<}\mbox{SimpleTorus}{>}_{4} & = & 1
\end{eqnarray*}

Amplitudes for ($\MU = \emptyset$):
\begin{eqnarray*}
{<}\mbox{SimpleTorus}{>}_{a} & = & 1 \\
{<}\mbox{SimpleTorus}{>}_{b} & = & 1
\end{eqnarray*}

\subsection{Example: 1-Twist Spun Trefoil}
\par\noindent
The next knot is a spun trefoil from~\cite{CS2}
\begin{verbatim}
% 1-twist spun trefoil: C&S p. 36
ss => f[0,Cup,0]ss[0,Cap,0]f =>
[0,Cup,0]f[0,Cup,2]ss[0,Cap,2]f[0,Cap,0] =>
[0,Cup,0][0,Cup,2]f[1,NE,1]ss[1,NW,1]f[0,Cap,2][0,Cap,0] =>
[0,Cup,0][0,Cup,2][1,NE,1]f[1,NE,1]ss[1,NW,1]f[1,NW,1];
   [0,Cap,2][0,Cap,0]
=> [0,Cup,0][0,Cup,2][1,NE,1][1,NE,1]f[1,NE,1]ss[1,NW,1]f
   [1,NW,1][1,NW,1][0,Cap,2][0,Cap,0]
=> [0,Cup,0][0,Cup,2][1,NE,1][1,NE,1][1,NE,1]f[2,Cap,0]s
   [2,Cup,0]f[1,NW,1]s[1,NW,1][1,NW,1][0,Cap,2][0,Cap,0]
=> [0,Cup,0][0,Cup,2][1,NE,1][1,NE,1][1,NE,1][2,Cap,0]sf
   [1,Cup,1]s[2,NE,0]f[1,NW,1][1,NW,1][0,Cap,2][0,Cap,0]
=> [0,Cup,0][0,Cup,2][1,NE,1][1,NE,1][1,NE,1][2,Cap,0]f
   [1,Cup,1]s[1,NE,1]f[2,NE,0][1,NW,1]s[1,NW,1];
   [0,Cap,2][0,Cap,0]
=> [0,Cup,0][0,Cup,2][1,NE,1][1,NE,1][1,NE,1][2,Cap,0];
   [1,Cup,1]f[2,NW,0][1,NE,1][2,NE,0]f[1,NW,1]s
   [0,Cap,2][0,Cap,0]s
=> [0,Cup,0][0,Cup,2][1,NE,1][1,NE,1][1,NE,1][2,Cap,0];
   [1,Cup,1][2,NW,0]s[1,NE,1][2,NE,0][1,NW,1]sf
   [2,Cap,0][0,Cap,0]f
=> [0,Cup,0][0,Cup,2][1,NE,1][1,NE,1][1,NE,1][2,Cap,0]s
   [1,Cup,1][2,NW,0]sf[2,NW,0][1,NE,1][2,NE,0]f
   [2,Cap,0][0,Cap,0]
=> [0,Cup,0][0,Cup,2][1,NE,1][1,NE,1][1,NE,1][2,Cap,0]f
   [2,Cup,0][1,NE,1]f[2,NW,0][1,NE,1]s[2,NE,0];
   [2,Cap,0]s[0,Cap,0]
=> [0,Cup,0][0,Cup,2][1,NE,1][1,NE,1][1,NE,1][2,Cap,0]s
   [2,Cup,0]s[1,NE,1][2,NW,0][1,NE,1]f[2,Cap,0]f[0,Cap,0]
=> [0,Cup,0][0,Cup,2][1,NE,1][1,NE,1][1,NE,1][2,Cap,0]f
   [2,Cup,0]s[2,NE,0]f[1,NE,1][2,NW,0]s[1,NE,1];
   [2,Cap,0][0,Cap,0]
=> [0,Cup,0][0,Cup,2][1,NE,1][1,NE,1][1,NE,1][2,Cap,0];
   [2,Cup,0]f[1,NW,1][2,NE,0][1,NE,1]sf[1,NE,1];
   [2,Cap,0]s[0,Cap,0]
=> [0,Cup,0][0,Cup,2][1,NE,1][1,NE,1][1,NE,1][2,Cap,0];
   [2,Cup,0][1,NW,1]s[2,NE,0][1,NE,1]f[2,NW,0]s
   [1,Cap,1]f[0,Cap,0]
=> [0,Cup,0][0,Cup,2][1,NE,1][1,NE,1][1,NE,1][2,Cap,0];
   [2,Cup,0][1,NW,1]f[1,NW,1][2,NE,0]s[1,NE,1]f
   [1,Cap,1]s[0,Cap,0]
=> [0,Cup,0][0,Cup,2][1,NE,1][1,NE,1][1,NE,1][2,Cap,0];
   [2,Cup,0][1,NW,1][1,NW,1]s[2,NE,0]f[1,Cap,1]sf[0,Cap,0]
=> [0,Cup,0][0,Cup,2][1,NE,1][1,NE,1][1,NE,1][2,Cap,0];
   [2,Cup,0][1,NW,1][1,NW,1]f[1,NW,1]s[2,Cap,0]f[0,Cap,0]s
=> [0,Cup,0][0,Cup,2][1,NE,1][1,NE,1][1,NE,1]s[2,Cap,0];
   [2,Cup,0]s[1,NW,1][1,NW,1][1,NW,1]f[0,Cap,2][0,Cap,0]f
=> [0,Cup,0][0,Cup,2][1,NE,1][1,NE,1]s[1,NE,1]ff[1,NW,1]s
   [1,NW,1][1,NW,1][0,Cap,2][0,Cap,0]
=> [0,Cup,0][0,Cup,2][1,NE,1]s[1,NE,1]ff[1,NW,1]s[1,NW,1];
   [0,Cap,2][0,Cap,0]
=> [0,Cup,0][0,Cup,2]s[1,NE,1]ff[1,NW,1]s[0,Cap,2][0,Cap,0] =>
[0,Cup,0]s[0,Cup,2]ff[0,Cap,2]s[0,Cap,0] => s[0,Cup,0]ff[0,Cap,0]s
=> ff #
\end{verbatim}
with amplitudes:

Amplitudes for ($\MU = \{31\}$):
\begin{eqnarray*}
{<}\mbox{1--TwistSpunTrefoil}{>}_{0} & = & q^{-14}+q^{14} \\
{<}\mbox{1--TwistSpunTrefoil}{>}_{1} & = & q^{-14}+q^{14} \\
{<}\mbox{1--TwistSpunTrefoil}{>}_{2} & = & q^{-14}+q^{14} \\
{<}\mbox{1--TwistSpunTrefoil}{>}_{3} & = & 1 \\
{<}\mbox{1--TwistSpunTrefoil}{>}_{4} & = & 1
\end{eqnarray*}

Amplitudes for ($\MU = \emptyset$):
\begin{eqnarray*}
{<}\mbox{1--TwistSpunTrefoil}{>}_{a} & = & 1 \\
{<}\mbox{1--TwistSpunTrefoil}{>}_{b} & = & 1
\end{eqnarray*}

\subsection{Example: 2-twist spun trefoil}
\par\noindent
The next knot is another spun trefoil from~\cite{CS2}
\begin{verbatim}
% 2-twist spun trefoil: C&S p. 36
ss => f[0,Cup,0]ss[0,Cap,0]f =>
[0,Cup,0]f[0,Cup,2]ss[0,Cap,2]f[0,Cap,0] =>
[0,Cup,0][0,Cup,2]f[1,NE,1]ss[1,NW,1]f[0,Cap,2][0,Cap,0] =>
[0,Cup,0][0,Cup,2][1,NE,1]f[1,NE,1]ss[1,NW,1]f[1,NW,1];
   [0,Cap,2][0,Cap,0]
=> [0,Cup,0][0,Cup,2][1,NE,1][1,NE,1]f[1,NE,1]ss[1,NW,1]f
   [1,NW,1][1,NW,1][0,Cap,2][0,Cap,0]
=> [0,Cup,0][0,Cup,2][1,NE,1][1,NE,1][1,NE,1]f[2,Cap,0]s
   [2,Cup,0]f[1,NW,1]s[1,NW,1][1,NW,1][0,Cap,2][0,Cap,0]
=> [0,Cup,0][0,Cup,2][1,NE,1][1,NE,1][1,NE,1][2,Cap,0]sf
   [1,Cup,1]s[2,NE,0]f[1,NW,1][1,NW,1][0,Cap,2][0,Cap,0]
=> [0,Cup,0][0,Cup,2][1,NE,1][1,NE,1][1,NE,1][2,Cap,0]f
   [1,Cup,1]s[1,NE,1]f[2,NE,0][1,NW,1]s[1,NW,1];
   [0,Cap,2][0,Cap,0]
=> [0,Cup,0][0,Cup,2][1,NE,1][1,NE,1][1,NE,1][2,Cap,0];
   [1,Cup,1]f[2,NW,0][1,NE,1][2,NE,0]f[1,NW,1]s
   [0,Cap,2][0,Cap,0]s
=> [0,Cup,0][0,Cup,2][1,NE,1][1,NE,1][1,NE,1][2,Cap,0];
   [1,Cup,1][2,NW,0]s[1,NE,1][2,NE,0][1,NW,1]sf
   [2,Cap,0][0,Cap,0]f
=> [0,Cup,0][0,Cup,2][1,NE,1][1,NE,1][1,NE,1][2,Cap,0]s
   [1,Cup,1][2,NW,0]sf[2,NW,0][1,NE,1][2,NE,0]f
   [2,Cap,0][0,Cap,0]
=> [0,Cup,0][0,Cup,2][1,NE,1][1,NE,1][1,NE,1][2,Cap,0]f
   [2,Cup,0][1,NE,1]f[2,NW,0][1,NE,1]s[2,NE,0];
   [2,Cap,0]s[0,Cap,0]
=> [0,Cup,0][0,Cup,2][1,NE,1][1,NE,1][1,NE,1][2,Cap,0]s
   [2,Cup,0]s[1,NE,1][2,NW,0][1,NE,1]f[2,Cap,0]f[0,Cap,0]
=> [0,Cup,0][0,Cup,2][1,NE,1][1,NE,1][1,NE,1][2,Cap,0]f
   [2,Cup,0]s[2,NE,0]f[1,NE,1][2,NW,0]s[1,NE,1];
   [2,Cap,0][0,Cap,0]
=> [0,Cup,0][0,Cup,2][1,NE,1][1,NE,1][1,NE,1][2,Cap,0];
   [2,Cup,0]f[1,NW,1][2,NE,0][1,NE,1]sf[1,NE,1];
   [2,Cap,0]s[0,Cap,0]
=> [0,Cup,0][0,Cup,2][1,NE,1][1,NE,1][1,NE,1][2,Cap,0];
   [2,Cup,0][1,NW,1]s[2,NE,0][1,NE,1]f[2,NW,0]s
   [1,Cap,1]f[0,Cap,0]
=> [0,Cup,0][0,Cup,2][1,NE,1][1,NE,1][1,NE,1][2,Cap,0];
   [2,Cup,0][1,NW,1]f[1,NW,1][2,NE,0]s[1,NE,1]f
   [1,Cap,1]s[0,Cap,0]
=> [0,Cup,0][0,Cup,2][1,NE,1][1,NE,1][1,NE,1][2,Cap,0];
   [2,Cup,0][1,NW,1][1,NW,1]s[2,NE,0]f
   [1,Cap,1]sf[0,Cap,0]
=> [0,Cup,0][0,Cup,2][1,NE,1][1,NE,1][1,NE,1][2,Cap,0];
   [2,Cup,0][1,NW,1][1,NW,1]f[1,NW,1]s
   [2,Cap,0]f[0,Cap,0]s
=> [0,Cup,0][0,Cup,2][1,NE,1][1,NE,1][1,NE,1][2,Cap,0]s
   [2,Cup,0][1,NW,1]s[1,NW,1][1,NW,1]f[0,Cap,2][0,Cap,0]f
=> [0,Cup,0][0,Cup,2][1,NE,1][1,NE,1][1,NE,1][2,Cap,0]sf
   [1,Cup,1]s[2,NE,0]f[1,NW,1][1,NW,1][0,Cap,2][0,Cap,0]
=> [0,Cup,0][0,Cup,2][1,NE,1][1,NE,1][1,NE,1][2,Cap,0]f
   [1,Cup,1]s[1,NE,1]f[2,NE,0][1,NW,1]s[1,NW,1];
   [0,Cap,2][0,Cap,0]
=> [0,Cup,0][0,Cup,2][1,NE,1][1,NE,1][1,NE,1][2,Cap,0];
   [1,Cup,1]f[2,NW,0][1,NE,1][2,NE,0]f[1,NW,1]s
   [0,Cap,2][0,Cap,0]s
=> [0,Cup,0][0,Cup,2][1,NE,1][1,NE,1][1,NE,1][2,Cap,0];
   [1,Cup,1][2,NW,0]s[1,NE,1][2,NE,0][1,NW,1]sf
   [2,Cap,0][0,Cap,0]f
=> [0,Cup,0][0,Cup,2][1,NE,1][1,NE,1][1,NE,1][2,Cap,0]s
   [1,Cup,1][2,NW,0]sf[2,NW,0][1,NE,1][2,NE,0]f
   [2,Cap,0][0,Cap,0]
=> [0,Cup,0][0,Cup,2][1,NE,1][1,NE,1][1,NE,1][2,Cap,0]f
   [2,Cup,0][1,NE,1]f[2,NW,0][1,NE,1]s[2,NE,0];
   [2,Cap,0]s[0,Cap,0]
=> [0,Cup,0][0,Cup,2][1,NE,1][1,NE,1][1,NE,1][2,Cap,0]s
   [2,Cup,0]s[1,NE,1][2,NW,0][1,NE,1]f[2,Cap,0]f[0,Cap,0]
=> [0,Cup,0][0,Cup,2][1,NE,1][1,NE,1][1,NE,1][2,Cap,0]f
   [2,Cup,0]s[2,NE,0]f[1,NE,1][2,NW,0]s[1,NE,1];
   [2,Cap,0][0,Cap,0]
=> [0,Cup,0][0,Cup,2][1,NE,1][1,NE,1][1,NE,1][2,Cap,0];
   [2,Cup,0]f[1,NW,1][2,NE,0][1,NE,1]sf[1,NE,1];
   [2,Cap,0]s[0,Cap,0]
=> [0,Cup,0][0,Cup,2][1,NE,1][1,NE,1][1,NE,1][2,Cap,0];
   [2,Cup,0][1,NW,1]s[2,NE,0][1,NE,1]f[2,NW,0]s
   [1,Cap,1]f[0,Cap,0]
=> [0,Cup,0][0,Cup,2][1,NE,1][1,NE,1][1,NE,1][2,Cap,0];
   [2,Cup,0][1,NW,1]f[1,NW,1][2,NE,0]s[1,NE,1]f
   [1,Cap,1]s[0,Cap,0]
=> [0,Cup,0][0,Cup,2][1,NE,1][1,NE,1][1,NE,1][2,Cap,0];
   [2,Cup,0][1,NW,1][1,NW,1]s[2,NE,0]f[1,Cap,1]sf[0,Cap,0]
=> [0,Cup,0][0,Cup,2][1,NE,1][1,NE,1][1,NE,1][2,Cap,0];
   [2,Cup,0][1,NW,1][1,NW,1]f[1,NW,1]s[2,Cap,0]f[0,Cap,0]s
=> [0,Cup,0][0,Cup,2][1,NE,1][1,NE,1][1,NE,1]s[2,Cap,0];
   [2,Cup,0]s[1,NW,1][1,NW,1][1,NW,1]f[0,Cap,2][0,Cap,0]f
=> [0,Cup,0][0,Cup,2][1,NE,1][1,NE,1]s[1,NE,1]ff[1,NW,1]s
   [1,NW,1][1,NW,1][0,Cap,2][0,Cap,0]
=> [0,Cup,0][0,Cup,2][1,NE,1]s[1,NE,1]ff[1,NW,1]s[1,NW,1];
   [0,Cap,2][0,Cap,0]
=> [0,Cup,0][0,Cup,2]s[1,NE,1]ff[1,NW,1]s[0,Cap,2][0,Cap,0] =>
[0,Cup,0]s[0,Cup,2]ff[0,Cap,2]s[0,Cap,0] => s[0,Cup,0]ff[0,Cap,0]s
=> ff #
\end{verbatim}
with amplitudes:

Amplitudes for ($\MU = \{31\}$):
\begin{eqnarray*}
{<}\mbox{2--TwistSpunTrefoil}{>}_{0} & = & q^{-14}+q^{14} \\
{<}\mbox{2--TwistSpunTrefoil}{>}_{1} & = & q^{-14}+q^{14} \\
{<}\mbox{2--TwistSpunTrefoil}{>}_{2} & = & q^{-14}+q^{14} \\
{<}\mbox{2--TwistSpunTrefoil}{>}_{3} & = & 1 \\
{<}\mbox{2--TwistSpunTrefoil}{>}_{4} & = & 1
\end{eqnarray*}

Amplitudes for ($\MU = \emptyset$):
\begin{eqnarray*}
{<}\mbox{2--TwistSpunTrefoil}{>}_{a} & = & 1 \\
{<}\mbox{2--TwistSpunTrefoil}{>}_{b} & = & 1
\end{eqnarray*}

\end{document}